\documentclass[11pt,letterpaper]{amsart}
\ifdefined\pdfminorversion\pdfminorversion=7\fi
\usepackage{amsmath,amssymb,amsthm}
\usepackage{graphicx}
\usepackage[percent]{overpic}
\newcommand{\R}{\mathbb{R}}
\newcommand{\Z}{\mathbb{Z}}
\newcommand{\Q}{\mathbb{Q}}
\newcommand{\C}{\mathbb{C}}

\newcommand{\ZN}[1]{\Z_{#1}}
\newcommand{\ex}[1]{e^{2\pi i #1}}
\newcommand{\ps}[1]{p_{#1}}
\newcommand{\es}[1]{e_{#1}}
\newcommand{\zroot}[1]{\zeta_{#1}}
\newcommand{\muS}{\mu_S}
\newcommand{\mask}[1]{#1(X)}

\theoremstyle{plain}
\newtheorem{theorem}{Theorem}[section]
\newtheorem{proposition}[theorem]{Proposition}
\newtheorem{lemma}[theorem]{Lemma}
\newtheorem{corollary}[theorem]{Corollary}
\theoremstyle{definition}
\newtheorem{definition}[theorem]{Definition}
\theoremstyle{remark}
\newtheorem{remark}[theorem]{Remark}

\usepackage[hidelinks,pdftitle={Spectral unions of four intervals tile the real line},pdfauthor={Alexander Okezue Bell}]{hyperref}
\title{Spectral unions of four intervals tile the real line}
\author{Alexander Okezue Bell}
\thanks{Department of Mathematics, Stanford University, 450 Jane Stanford Way, Building 380, Stanford, CA 94305-2125, USA. Corresponding author: Alexander Okezue Bell. Email: \href{mailto:okezue@stanford.edu}{\texttt{okezue@stanford.edu}}.}
\date{}
\begin{document}
\begin{abstract}
We prove that every spectral union of four bounded intervals tiles the real line by translations. The argument combines periodicity of spectra, cofactor identities for the secular determinant, and a common tiling complement for the finite spectral fibers. We also establish a rigidity theorem for six complex numbers with prescribed vanishing power sums.
\end{abstract}
\maketitle
\section{Introduction}
\label{sec:intro}

Fuglede proved that a bounded measurable subset of Euclidean space tiles by a lattice iff it is spectral with the dual lattice as a spectrum, and he asked whether spectrality and translational tiling remain equivalent when the lattice restriction is dropped \cite{fuglede1974}. In dimension one the question is settled for a union of two intervals, by \L aba \cite{laba2001}, and for a union of three. Bose, Anil Kumar, Krishnan and Madan proved that a union of three intervals which tiles the real line is spectral \cite{bose-etal2010}. Shi proved the converse, and with it Fuglede's conjecture for unions of at most three intervals \cite[Theorem~1.1]{shi2026}.

\begin{theorem}
\label{thm:main}
Let \(\Omega \subset \R\) be a bounded measurable set which, up to a null set, is a union of four bounded intervals. If \(\Omega\) is spectral, then \(\Omega\) tiles \(\R\) by translations.
\end{theorem}

A set is spectral when there is a set of frequencies for which the normalized exponentials form an orthonormal basis of \(L^2(\Omega)\), and it tiles by translations when countably many translates partition the line up to null sets. The two definitions are stated precisely in Section~\ref{sec:fiber}. The proof of Theorem~\ref{thm:main} is given in Section~\ref{sec:assembly}, after the reductions it uses.

Shi's theorem comes with a geometric description \cite[Theorem~1.2]{shi2026}. If
\[
\Omega = [a_1,b_1) \cup [a_2,b_2) \cup [a_3,b_3)
\]
is spectral, of total length \(L\), then after a relabeling of the components at least one of three alternatives holds. The set is a fundamental domain of the lattice \(L\Z\). Or the three lengths are equal to \(h = L/3\), and \(\Omega\) is a three-fold lattice multitile by \(h\Z\). Or one length is the sum of the other two and equals \(h = L/2\), the two shorter components occupy complementary arcs of the circle \(\R/h\Z\), and \(\Omega\) is a two-fold lattice multitile by \(h\Z\). In the second and third alternatives a separate argument turns the lattice multitile into a translational tiling. The route to this trichotomy runs through a self-adjoint realization of the derivative on the union of intervals. Its domain is encoded by a boundary unitary, the spectrum is the zero set of the associated secular determinant, and cofactor identities on that determinant constrain the unitary and the lengths.

The same operator, determinant and rank-one cofactor identity are available for any finite number of intervals. For more than three components, Shi's classification no longer closes the argument \cite[Section~10]{shi2026}. For a larger number of lengths the cofactor identities involve many principal minors and many additive relations among subset sums, and the exhaustive list that closes the argument for three lengths does not directly extend. Separately, a zero sum of two unimodular numbers is antipodal, and a zero sum of three is a rotated equilateral triangle, which is the rigidity his fiber lemma uses for two-fold and three-fold lattice multitiles. A zero sum of four unimodular numbers has continuous degrees of freedom, so that lemma does not extend by the same reasoning. Theorem~\ref{thm:main} is the four-interval case of the missing implication. Its proof does not attempt a complete classification of the boundary unitary. It combines the cofactor identities with the periodicity of the spectrum and with the finite spectral pairs carried by the fibers of the period.

\begin{figure}[ht]
\centering
\begin{overpic}[width=0.9\textwidth]{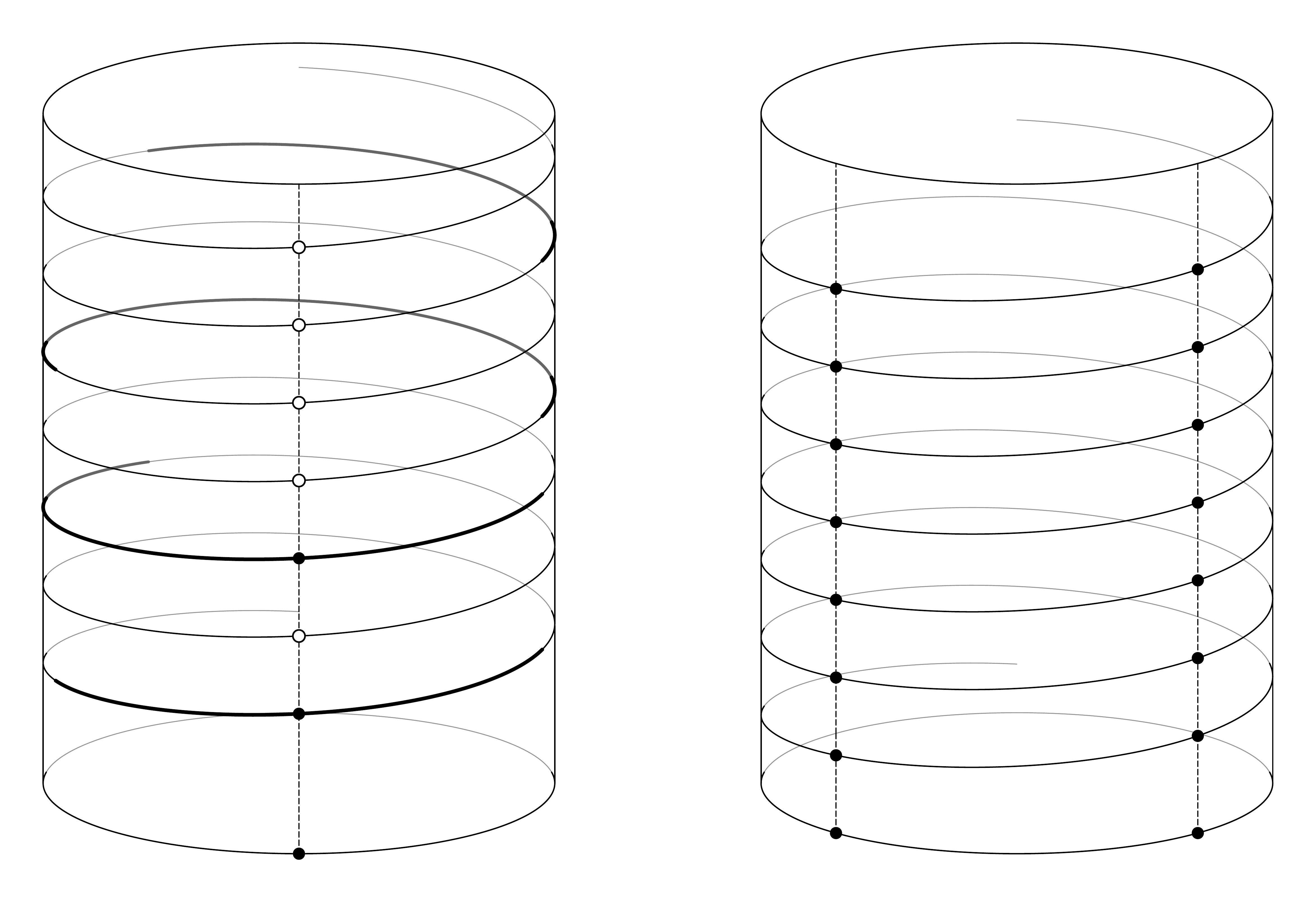}
\put(2.04,29.61){\makebox(0,0)[r]{$\Omega$}}
\put(22.72,0.96){\makebox(0,0){$x$}}
\put(97.96,35.82){\makebox(0,0)[l]{$\Lambda$}}
\put(77.28,0.96){\makebox(0,0){$S$}}
\end{overpic}
\caption{On the left, the line is wound on a cylinder of circumference \(h=1\), with the four intervals of \(\Omega=[0.6,1)\cup[2.4,3)\cup[4,4.6)\cup[6,6.4)\) drawn heavy and the parts behind the cylinder in gray. The vertical line over \(x\) meets the winding at the points \(x+mh\), and the filled ones form the fiber \(A_x=\{0,2\}\). On the right, the spectrum \(\Lambda=\{0,\tfrac14\}+\Z\) is wound on a cylinder of circumference \(h^{-1}\). It lies on \(k=2\) vertical lines, whose positions on the base circle, measured in turns, are the points of \(S\).}
\label{fig:fiber}
\end{figure}

A spectrum of a bounded measurable set on the line is periodic \cite{ik2013}. If \(h^{-1}\) is the minimal period, intersecting \(\Omega\) with the translates of \(h\Z\) turns the continuous problem into a family of finite ones. Each such translate meets \(\Omega\) in finitely many points, and the integers \(m\) for which \(x+mh\in\Omega\) form the fiber at \(x\). Spectrality of \(\Omega\) with a spectrum made of \(k\) cosets of \(h^{-1}\Z\) is equivalent to almost every fiber being a spectrum of one and the same set \(S\) of \(k\) points on the circle, and then \(\lvert\Omega\rvert=kh\). Figure~\ref{fig:fiber} shows one such cut. The proof bounds \(k\), and then turns the common structure of the fibers into a tiling.

The bound on \(k\) comes from the boundary unitary. The secular determinant is an exponential polynomial whose frequencies are the subset sums of the four lengths, and its zero set is the spectrum. Section~\ref{sec:boundary} proves these facts for any finite number of intervals, and derives from the periodicity of the spectrum the uniqueness principle for exponential polynomials that vanish on it. The cofactor identity then says that, at each value of a subset sum, the coefficient-weighted sum of the subset indicators is proportional to the all-ones vector. A nonempty proper subset can therefore carry a nonzero coefficient only if its length sum is shared with other subsets, in such a way that their indicators and the all-ones vector are linearly dependent. For four lengths these values, the levels, are classified in Theorem~\ref{thm:levelclass}. This is where the argument departs from the three-interval case: it classifies the levels, and the boundary unitary itself is analyzed only for two proportions of the lengths. Unless the lengths are proportional to \((1,1,1,2)\) or to \((1,1,2,2)\), every level is a multiple of \(\lvert\Omega\rvert/k_0\) for some \(k_0\le 4\). The secular determinant is then a polynomial in \(\ex{\lvert\Omega\rvert\zeta/k_0}\), its zero set is invariant under translation by \(k_0/\lvert\Omega\rvert\), and \(k\) divides \(k_0\). Lengths proportional to \((1,1,1,2)\) force a monomial boundary unitary and \(k=1\). Lengths proportional to \((1,1,2,2)\) give \(k\le 3\), except in one configuration: up to scaling, \(\Omega\) is a union of integer runs of lengths \(2\), \(2\), \(1\) and \(1\), the fiber is constant, and \(k=6\).

When \(k\le 4\), the passage from the fibers to a tiling follows from Dutkay and Haussermann's common-complement theorem \cite{dutkay-haussermann2015}, which applies to every \(k\le 5\). Section~\ref{sec:four} gives an explicit construction in the notation used here. Shi's fiber lemma \cite[Lemma~7.1]{shi2026} covers \(k\in\{2,3\}\). Since a zero sum of two or three unimodular numbers is rigid, every fiber reduces, modulo an integer determined by \(S\), to a coset of one cyclic subgroup, and an interval of residues is a common tiling complement. For \(k=4\) the zero sums of four unimodular numbers do have a continuous parameter, but each of them splits into two antipodal pairs, and a dephased Hadamard matrix of order four, whose first row and column consist of ones, has a column with entries \(\pm 1\). These two facts force every fiber, modulo an integer determined by \(S\), to be a union of two cosets of a subgroup of order two, and one set \(B\subset\Z\), depending only on \(S\), is again a tiling complement of every fiber. Then \(\Omega\) tiles \(\R\) by \(hB\). For five cosets, Haagerup's classification of complex Hadamard matrices of order five \cite{haagerup1997} gives the same coset structure as for two and three, although Theorem~\ref{thm:main} does not need that case. Kolountzakis and Matolcsi used these classifications to prove that a spectral set of at most five elements in a finite abelian group tiles \cite{km2006}. On the line, a single complement must serve all the fibers at once. We record the resulting statement below and prove it as Theorem~\ref{thm:fiberfew}.

\begin{theorem}
\label{thm:introcomplement}
Let \(E\subset\R\) be a bounded measurable set with a spectrum \(\Gamma+h^{-1}\Z\), where \(h>0\) and \(\Gamma\) consists of \(k\le 5\) elements that are pairwise distinct modulo \(h^{-1}\Z\). Then \(E\) tiles \(\R\) by a translation set \(hB\), where \(B\subset\Z\) depends only on \(h\Gamma\subset\R/\Z\).
\end{theorem}

The six-point fiber of the exceptional configuration does not occur. If a finite set \(S\subset\R/\Z\) admits an integer spectrum, then \(S-S\subset\Q/\Z\), by a theorem of Fu and Song \cite{fu-song2026}, recorded here as Theorem~\ref{thm:rational}. Once the differences are rational, an integer spectrum containing two consecutive integers forces six roots of unity summing to zero into one of two rigid shapes, two rotated equilateral triangles or three antipodal pairs, and either shape gives \(S\) a translational symmetry, which the minimal period forbids.

The same fiber leads to a second statement that is not about tiling, and that may be of independent interest. If \(\{0,1,x,x+1,y,z\}\) is a spectrum of \(S\), then the six numbers \(\ex{s}\), \(s\in S\), have vanishing power sums at the exponents \(1\), \(x-1\), \(x\) and \(x+1\). These four relations alone force a nontrivial translation symmetry of \(S\) and the rationality of \(S-S\).

\begin{theorem}
\label{thm:introsix}
Let \(x\ge 2\) be an integer, and let \(\omega_1,\ldots,\omega_6\) be distinct complex numbers of modulus \(1\) such that \(\sum_{s=1}^6\omega_s^n=0\) for \(n=1\), \(x-1\), \(x\) and \(x+1\). Then the six numbers form three antipodal pairs or two rotated equilateral triangles, and every ratio \(\omega_s/\omega_t\) is a root of unity.
\end{theorem}

Section~\ref{sec:sixnumbers} proves more. Six nonzero complex numbers with product \(1\) whose power sums vanish at \(\pm 1\), \(\pm(x-1)\), \(\pm x\) and \(\pm(x+1)\) have \(\es{3}=0\) or \(\es{2}=\es{4}=0\), and all of them are roots of unity, of order at most \(6(x+1)\). Here \(\es{j}\) denotes the \(j\)th elementary symmetric function. The proof is local-global. For algebraic solutions, non-archimedean estimates make all six numbers units at every finite place, while a two-sided extension of Fu and Song's modulus rigidity theorem makes all their conjugates unimodular. Kronecker's theorem then makes them roots of unity, the classification of vanishing sums gives the stated alternatives, and the Nullstellensatz carries the resulting coefficient identities to complex solutions. More generally, the modulus theorem forces \(d\) nonzero complex numbers to have equal modulus when their power sums vanish on the nonzero differences of a set of \(N\) integers and \(d<2N\); this threshold is sharp.

A separate construction gives a statement about weak tilings. Section~\ref{sec:weak} exhibits a weak tiling of the complement of a four-interval set by a positive measure which is not a convex combination of proper tilings.

The paper is organized as follows. Section~\ref{sec:fiber} fixes the definitions, recalls the periodicity theorem, and proves the fiber reformulation. Section~\ref{sec:boundary} constructs the boundary unitary and the secular determinant of a spectral union of intervals. Section~\ref{sec:cosets} bounds the number of points in a fiber. Section~\ref{sec:four} constructs a common tiling complement for lattice multitiles with at most five cosets. Section~\ref{sec:six} excludes the six-point fiber. Section~\ref{sec:assembly} assembles these statements into the proof of Theorem~\ref{thm:main}. Section~\ref{sec:sixnumbers} proves the theorem on six unimodular numbers. Section~\ref{sec:weak} gives the weak-tiling counterexample. Section~\ref{sec:remarks} records what the method gives for an arbitrary finite number of intervals.

\section{Periodicity and the fiber reformulation}
\label{sec:fiber}

Throughout, \(\Omega \subset \R\) is a bounded measurable set of Lebesgue measure \(L = \lvert\Omega\rvert > 0\). The inner product on \(L^2(\Omega)\) is \(\langle f, g \rangle = \int_\Omega f\,\overline{g}\). The Fourier transform of its indicator is \(\widehat{1_\Omega}(\xi) = \int_\Omega \ex{-\xi t}\,dt\).

\begin{definition}
\label{def:spectrum}
For \(\lambda \in \R\) the normalized exponential is the function
\[
e_\lambda(x) = L^{-1/2} \ex{\lambda x}, \qquad x \in \Omega.
\]
A set \(\Lambda \subset \R\) is a \emph{spectrum} of \(\Omega\), and \(\Omega\) is \emph{spectral}, when \(\{e_\lambda : \lambda \in \Lambda\}\) is an orthonormal basis of \(L^2(\Omega)\). Orthogonality is equivalent to the vanishing of \(\widehat{1_\Omega}\) on \((\Lambda - \Lambda) \setminus \{0\}\), since \(\langle e_\lambda, e_{\lambda'} \rangle = L^{-1} \widehat{1_\Omega}(\lambda' - \lambda)\). The set \(\Omega\) \emph{tiles \(\R\) by translations} when there is a countable set \(T \subset \R\) such that the translates \(\Omega + t\), \(t \in T\), are pairwise disjoint up to null sets and their union is \(\R\) up to a null set.
\end{definition}

\begin{definition}
\label{def:fiber}
Let \(h > 0\), and let \(\Gamma \subset \R\) be a finite set whose elements are pairwise distinct modulo \(h^{-1}\Z\). Write \(k = \lvert\Gamma\rvert\), \(\Lambda = \Gamma + h^{-1}\Z\) and \(S = h\Gamma \subset \R/\Z\). Then \(\lvert S\rvert = k\) and \(\Lambda = h^{-1}(S + \Z)\). For \(x \in \R\) the \emph{fiber} of \(\Omega\) at \(x\) is
\[
A_x = \{ m \in \Z : x + mh \in \Omega \}.
\]
The set \(A_x\) is finite. Indeed, if \(\Omega \subset [-R, R]\), then \(x + mh \in \Omega\) forces \(\lvert m\rvert \le (R + \lvert x\rvert)/h\), and for \(x \in [0, h)\) the bound is \(\lvert m\rvert \le (R + h)/h\). The same bound shows that \(x \mapsto \lvert A_x\rvert\) and, for each \(\sigma \in \R/\Z\), the character sum
\[
\widehat{A_x}(\sigma) = \sum_{m \in A_x} \ex{\sigma m}
\]
are bounded measurable functions of \(x \in [0, h)\): each is a finite sum of translates of \(1_\Omega\). Translation by one period shifts the fiber by one integer, \(A_{x+h} = A_x - 1\), and therefore \(\widehat{A_{x+h}}(\sigma) = \ex{-\sigma}\,\widehat{A_x}(\sigma)\).
\end{definition}

\begin{definition}
\label{def:spectralpair}
Let \(S \subset \R/\Z\) be finite, of cardinality \(k\), and let \(A \subset \Z\) be finite. The set \(A\) is a \emph{spectrum} of \(S\) when \(\lvert A\rvert = k\) and \(\widehat{A}(\sigma) = 0\) \(\forall \sigma \in (S - S) \setminus \{0\}\). Equivalently, the matrix \(H = \bigl(\ex{sa}\bigr)_{s \in S,\, a \in A}\) is a complex Hadamard matrix of order \(k\): every entry has modulus \(1\), and \(HH^* = kI\). The entries of \(HH^*\) are the sums \(\widehat{A}(s - s')\). Since a left inverse of a square matrix is a two-sided inverse, one then has \(H^*H = kI\) as well, and the counting measure \(\mu_S = \sum_{s \in S} \delta_s\), with the positive-sign transform \(\widehat{\mu_S}(n) = \sum_{s \in S}\ex{ns}\) for \(n \in \Z\), satisfies \(\widehat{\mu_S}(a - a') = 0\) whenever \(a \ne a'\) in \(A\).
\end{definition}

The next statement is the periodicity theorem of Iosevich and Kolountzakis \cite[Theorem~1.5]{ik2013}, scaled from measure one to measure \(L\), together with the elementary description of the period group.

\begin{theorem}
\label{thm:periodic}
Let \(\Lambda\) be a spectrum of a bounded measurable set \(\Omega \subset \R\) of measure \(L\). Every positive period of \(\Lambda\) lies in \(L^{-1}\Z_{>0}\). The set of all periods is \(h^{-1}\Z\) for a unique \(h > 0\), so that \(h^{-1}\) is the minimal positive period. Writing \(\Gamma = \Lambda \cap [0, h^{-1})\), \(k = \lvert\Gamma\rvert\) and \(S = h\Gamma \subset \R/\Z\) as in Definition~\ref{def:fiber}, one has \(\Lambda = \Gamma + h^{-1}\Z\) and \(\lvert S\rvert = k\), and \(S\) has no nontrivial translation symmetry.
\end{theorem}

\begin{proof}
Let \(\Omega_1 = L^{-1}\Omega\) and \(\Lambda_1 = L\Lambda\). The change of variables \(x = Ly\) gives \(\lvert\Omega_1\rvert = 1\). If \(\lambda \ne \lambda'\) lie in \(\Lambda\), then \(\int_\Omega \ex{(\lambda - \lambda')x}\,dx = 0\), so
\[
\int_{\Omega_1} \ex{(L\lambda - L\lambda') y}\,dy = L^{-1} \int_\Omega \ex{(\lambda - \lambda') x}\,dx = 0,
\]
while the integral equals \(1\) when \(\lambda = \lambda'\). Thus the unnormalized exponentials at the frequencies \(\Lambda_1\) are orthonormal in \(L^2(\Omega_1)\). They are complete. If \(f \in L^2(\Omega_1)\) is orthogonal to every such exponential, the substitution \(x = Ly\) gives
\[
\int_\Omega f(x/L)\,\overline{e_\lambda(x)}\,dx = L^{1/2} \int_{\Omega_1} f(y)\,\ex{-(L\lambda) y}\,dy = 0\quad\forall\lambda \in \Lambda.
\]
Completeness of \(\{e_\lambda\}\) forces \(x \mapsto f(x/L)\) to vanish almost everywhere on \(\Omega\), and therefore \(f\) vanishes almost everywhere on \(\Omega_1\). By the theorem of Iosevich and Kolountzakis \cite[Theorem~1.5]{ik2013}, \(\Lambda_1\) is periodic and every positive period of \(\Lambda_1\) is a positive integer. Every positive period of \(\Lambda\) is therefore of the form \(T/L\) with \(T \in \Z_{>0}\).

The Fourier transform \(\widehat{1_\Omega}\) is continuous, \(\widehat{1_\Omega}(0) = L\), and \(\widehat{1_\Omega}\) vanishes on \((\Lambda - \Lambda) \setminus \{0\}\). Hence there is a \(\delta > 0\) such that any two distinct points of \(\Lambda\) differ by at least \(\delta\). The period group \(P = \{ t \in \R : \Lambda + t = \Lambda \}\) is a nonzero subgroup of \(L^{-1}\Z\), by the periodicity and integrality of periods just cited. Thus \(P = h^{-1}\Z\) for a unique \(h > 0\), and its positive generator \(h^{-1}\) lies in \(L^{-1}\Z_{>0}\).

Uniform discreteness makes \(\Gamma = \Lambda \cap [0, h^{-1})\) finite, and \(\Lambda = \Gamma + h^{-1}\Z\). If \(S + \alpha = S\) for some \(\alpha \in \R/\Z\) which is not the identity, any lift of \(\alpha\) satisfies \(\Lambda + \alpha/h = \Lambda\), so \(\alpha/h \in h^{-1}\Z\) and \(\alpha \in \Z\), a contradiction. Thus \(S\) has no nontrivial translation symmetry.
\end{proof}

After the change of variables \(t=x+mh\), frequencies in one coset of \(h^{-1}\Z\) have the same dependence on \(m\). The following theorem identifies spectrality with the corresponding basis property on almost every fiber.

\begin{theorem}
\label{thm:fiber}
Let \(\Omega \subset \R\) be bounded and measurable, of measure \(L\), and let \(h\), \(\Gamma\), \(k\), \(\Lambda\) and \(S\) be as in Definition~\ref{def:fiber}. The following are equivalent.
\begin{enumerate}
\item[\rm (i)] \(\Lambda\) is a spectrum of \(\Omega\).
\item[\rm (ii)] For almost every \(x \in [0, h)\), the fiber \(A_x\) is a spectrum of \(S\).
\end{enumerate}
When these hold, \(L = kh\).
\end{theorem}

\begin{proof}
Set
\[
\widetilde{\Omega} = \{ (x, m) \in [0, h) \times \Z : x + mh \in \Omega \}.
\]
Equip \(\widetilde{\Omega}\) with the restriction of Lebesgue measure times counting measure, and write \(L^2(\widetilde{\Omega})\) for the corresponding space of square-integrable functions. Define \(V \colon L^2(\Omega) \to L^2(\widetilde{\Omega})\) by \((Vf)(x, m) = f(x + mh)\). Every real number has a unique representation \(t = x + mh\) with \(x \in [0, h)\) and \(m \in \Z\), so the map \((x, m) \mapsto x + mh\) is a bijection of \(\widetilde{\Omega}\) onto \(\Omega\). On each slice the map is a translation, and Lebesgue measure is translation invariant, so \(V\) is a linear isometry. It is surjective: given \(F \in L^2(\widetilde{\Omega})\), the function \(f\) on \(\Omega\) defined by \(f(x + mh) = F(x, m)\) is measurable, square-integrable, and satisfies \(Vf = F\). Thus \(V\) is unitary.

For \(\gamma \in \R\), \(q \in \Z\) and \(\lambda = \gamma + q/h\),
\begin{equation}
\label{eq:fibermap}
(Ve_\lambda)(x, m) = L^{-1/2} \ex{(\gamma + q/h)x} \ex{h\gamma m},
\end{equation}
because \(\ex{qm} = 1\). The field of vectors \(v_\gamma(x) = \bigl(\ex{h\gamma m}\bigr)_{m \in A_x}\) is therefore independent of \(q\). Its inner product on the fiber is
\[
\langle v_\gamma(x), v_{\gamma'}(x) \rangle = \widehat{A_x}\bigl(h(\gamma - \gamma')\bigr).
\]
For \(\gamma,\gamma' \in \Gamma\) with \(\gamma \neq \gamma'\), the frequency \(h(\gamma - \gamma')\) is a nonzero point of \(S - S\), and every nonzero point of \(S - S\) arises in this way. Applying the change of variables to the integrable function \(e_{\gamma+q/h}\,\overline{e_{\gamma'+q'/h}}\) gives the inner-product formula
\begin{equation}
\label{eq:fiberinner}
\langle e_{\gamma+q/h}, e_{\gamma'+q'/h} \rangle = L^{-1} \int_0^h \ex{(\gamma - \gamma')x} \ex{(q - q')x/h} \widehat{A_x}\bigl(h(\gamma - \gamma')\bigr)\,dx.
\end{equation}
The same change of variables with the character \(t \mapsto \ex{-qt/h}\) gives
\begin{equation}
\label{eq:poisson}
\widehat{1_\Omega}(q/h) = \int_0^h \ex{-qx/h} \lvert A_x\rvert\,dx, \qquad q \in \Z,
\end{equation}
and in particular \(\int_0^h \lvert A_x\rvert\,dx = L\). One standard fact is used repeatedly: if \(\varphi \in L^1([0, h))\) and \(\int_0^h \varphi(x) \ex{-qx/h}\,dx = 0\) \(\forall q \in \Z\), then \(\varphi = 0\) almost everywhere.

Assume (i). Fix \(\gamma \in \Gamma\). The exponentials \(e_{\gamma+q/h}\), \(q \in \Z\), are pairwise orthogonal, so \(\widehat{1_\Omega}\) vanishes on \(h^{-1}\Z \setminus \{0\}\) by the formula in Definition~\ref{def:spectrum}. Combined with \eqref{eq:poisson}, the function \(\lvert A_x\rvert - L/h\) on \([0, h)\) has all Fourier coefficients equal to zero. Hence \(\lvert A_x\rvert = L/h\) for almost every \(x \in [0, h)\). The value is an integer on a set of positive measure, so \(k' := L/h\) is a positive integer.

Next fix \(\gamma \ne \gamma'\) in \(\Gamma\), and set \(\delta = \gamma - \gamma'\) and \(\sigma = h\delta \bmod 1\). For every \(q \in \Z\) the frequencies \(\gamma + q/h\) and \(\gamma'\) are distinct points of \(\Lambda\). Formula \eqref{eq:fiberinner}, with the second integer equal to zero, says that
\[
\int_0^h \ex{qx/h} \bigl(\ex{\delta x}\,\widehat{A_x}(\sigma)\bigr)\,dx = 0\quad\forall q \in \Z.
\]
As \(q\) runs through \(\Z\), these are all the Fourier coefficients of \(x \mapsto \ex{\delta x}\,\widehat{A_x}(\sigma)\). Thus \(\widehat{A_x}(\sigma) = 0\) for almost every \(x\). Finitely many pairs \(\gamma \ne \gamma'\) produce finitely many null sets. Outside a single null set \(N_1\) one has \(\lvert A_x\rvert = k'\) and \(\widehat{A_x}(\sigma) = 0\) \(\forall \sigma \in (S - S) \setminus \{0\}\).

At any such \(x\) the vectors \(v_\gamma(x)\), \(\gamma \in \Gamma\), are pairwise orthogonal and of squared norm \(k' \ge 1\) in \(\C^{A_x}\). They are linearly independent, so \(k \le k'\). Suppose \(k' > k\). The fiber takes only finitely many values on \([0, h)\), so there is a set \(A \subset \Z\) with \(\lvert A\rvert = k'\) for which
\[
X_A = \{ x \in [0, h) \setminus N_1 : A_x = A \}
\]
has positive measure. The corresponding vectors in \(\C^A\) are \(k\) orthogonal vectors in a space of dimension \(k'\), so some \(w \in \C^A \setminus \{0\}\) is orthogonal to all of them. Define \(F \in L^2(\widetilde{\Omega})\) by \(F(x, m) = w_m\) for \(x \in X_A\) and \(m \in A\), and by \(F = 0\) elsewhere. Then \(F\) is not zero. Let \(f = V^{-1}F\). For \(\lambda = \gamma + q/h\), unitarity and \eqref{eq:fibermap} give
\[
\langle f, e_\lambda \rangle = L^{-1/2} \int_{X_A} \ex{-(\gamma + q/h)x} \langle w, v_\gamma \rangle\,dx = 0,
\]
where \(v_\gamma = \bigl(\ex{h\gamma m}\bigr)_{m \in A}\). Thus \(f\) is a nonzero function orthogonal to every exponential of \(\Lambda\), which contradicts completeness. Therefore \(k' = k\), \(L = kh\), and (ii) holds.

Assume (ii), and write \(X_0\) for the full-measure set of \(x \in [0, h)\) at which \(A_x\) is a spectrum of \(S\). Formula \eqref{eq:poisson} at \(q = 0\) gives \(L = \int_0^h \lvert A_x\rvert\,dx = kh\).

Let \(\gamma, \gamma' \in \Gamma\) and \(q, q' \in \Z\). If \(\gamma \ne \gamma'\), then \(\widehat{A_x}\) vanishes at \(h(\gamma - \gamma')\) \(\forall x \in X_0\), so \eqref{eq:fiberinner} gives orthogonality. If \(\gamma = \gamma'\), then \(\widehat{A_x}(0) = k\) on \(X_0\), and \eqref{eq:fiberinner} gives
\[
\langle e_{\gamma+q/h}, e_{\gamma+q'/h} \rangle = L^{-1} \int_0^h \ex{(q - q')x/h}\, k\,dx = \frac{kh}{L} \delta_{q,q'} = \delta_{q,q'}.
\]
Distinct pairs \((\gamma, q)\) give distinct frequencies, because the elements of \(\Gamma\) are distinct modulo \(h^{-1}\Z\). The family \(\{e_\lambda : \lambda \in \Lambda\}\) is therefore orthonormal.

For \(x \in X_0\) the same vectors \(v_\gamma(x)\) are pairwise orthogonal of squared norm \(k\) in the \(k\)-dimensional space \(\C^{A_x}\), so they form a basis. Let \(f \in L^2(\Omega)\) be orthogonal to every \(e_\lambda\), and set \(F = Vf\). For each \(\gamma \in \Gamma\) the function
\[
\varphi_\gamma(x) = \sum_{m \in A_x} F(x, m) \ex{-h\gamma m}
\]
is measurable. Cauchy's inequality gives
\[
\lvert\varphi_\gamma(x)\rvert^2 \le \lvert A_x\rvert \sum_{m \in A_x} \lvert F(x, m)\rvert^2,
\]
so \(\varphi_\gamma \in L^2([0, h))\). Unitarity and \eqref{eq:fibermap} give, \(\forall q \in \Z\),
\[
0 = \langle f, e_{\gamma+q/h} \rangle = L^{-1/2} \int_0^h \ex{-qx/h} \bigl(\ex{-\gamma x}\,\varphi_\gamma(x)\bigr)\,dx.
\]
Hence \(\varphi_\gamma = 0\) almost everywhere. Outside one null set, \(F(x, \cdot)\) is orthogonal to a basis of \(\C^{A_x}\), so \(F(x, \cdot) = 0\). Therefore \(f = 0\), and \(\Lambda\) is a spectrum of \(\Omega\).
\end{proof}

\section{The boundary unitary}
\label{sec:boundary}

Theorem~\ref{thm:fiber} holds for every bounded measurable set, and it places no restriction on \(k\). For a union of intervals the restriction comes from the lengths, through a self-adjoint realization of the derivative. Orthogonality already indicates the mechanism. Let \(\Omega=\bigsqcup_{j=1}^{n}[a_j,b_j)\) be a union of \(n\) pairwise disjoint bounded intervals, with \(\ell_j=b_j-a_j>0\) and \(L=\lvert\Omega\rvert\), and let \(\Lambda\) be a spectrum of \(\Omega\). For \(\lambda\neq\mu\) in \(\Lambda\), integration over each component gives
\[
0=\int_\Omega\ex{(\lambda-\mu)x}\,dx=\frac{1}{2\pi i(\lambda-\mu)}\sum_{j=1}^{n}\Bigl(\ex{\lambda b_j}\,\overline{\ex{\mu b_j}}-\ex{\lambda a_j}\,\overline{\ex{\mu a_j}}\Bigr).
\]
The vectors \((\ex{\lambda a_j})_{j}\), which record an exponential of the spectrum at the left ends of the components, and the vectors \((\ex{\lambda b_j})_{j}\), which record it at the right ends, therefore have the same inner products, and every one of them has squared norm \(n\). The assignment of the first family to the second extends to an isometry between their spans. Lemma~\ref{lem:boundaryunitary} shows that this coupling is a single unitary matrix, determined by \(\Lambda\), and Proposition~\ref{prop:secular} recovers \(\Lambda\) as the zero set of one exponential polynomial built from that matrix and the lengths.

This is the framework of Shi \cite[Section~4]{shi2026}, who treats three intervals and records that it extends to any number of them \cite[Section~10]{shi2026}. The description of self-adjoint realizations of the derivative by unitary boundary conditions goes back to Fuglede \cite{fuglede1974}; see also Dutkay and Jorgensen \cite[Section~3]{dutkay-jorgensen2015}. We recall the construction for \(n\) intervals, in the notation needed here. The uniqueness principle of Lemma~\ref{lem:freqwidth} is derived from the periodicity of \(\Lambda\) rather than from a count of zeros: periodicity reduces it to the fact that a polynomial of degree less than \(k\) cannot vanish at \(k\) points.

Write \(W=\bigoplus_{j=1}^{n}H^1(a_j,b_j)\) for the space of functions on \(\Omega\) whose restriction to each component is absolutely continuous on \([a_j,b_j]\) with derivative in \(L^2\). For \(f\in W\) put \(f(a)=(f(a_j))_{j=1}^{n}\) and \(f(b)=(f(b_j))_{j=1}^{n}\), the values at the two ends of the components, and write \(\langle x,y\rangle=\sum_j x_j\overline{y_j}\) on \(\C^n\). Integration by parts on each component gives Green's formula
\begin{equation}
\label{eq:green}
\langle f',g\rangle+\langle f,g'\rangle=\langle f(b),g(b)\rangle-\langle f(a),g(a)\rangle,\qquad f,g\in W.
\end{equation}
Let \(D\) be the operator that is diagonal in the orthonormal basis \(\{e_\lambda:\lambda\in\Lambda\}\),
\begin{equation}
\label{eq:diagop}
Df=\sum_{\lambda\in\Lambda}\lambda\langle f,e_\lambda\rangle e_\lambda,\qquad \operatorname{dom}D=\Bigl\{f\in L^2(\Omega):\sum_{\lambda\in\Lambda}\lambda^2\lvert\langle f,e_\lambda\rangle\rvert^2<\infty\Bigr\}.
\end{equation}
It is self-adjoint. If \(Df=\zeta f\) for some \(\zeta\in\C\), then \((\zeta-\lambda)\langle f,e_\lambda\rangle=0\) \(\forall \lambda\), so the eigenvalues of \(D\) are the points of \(\Lambda\) and the eigenspace of \(\lambda\) is \(\C e_\lambda\).

\begin{lemma}
\label{lem:boundaryunitary}
Let \(\Omega=\bigsqcup_{j=1}^{n}[a_j,b_j)\) be a union of \(n\) pairwise disjoint bounded intervals, let \(\Lambda\) be a spectrum of \(\Omega\), and let \(D\) be the operator \eqref{eq:diagop}. Then there is a unique unitary matrix \(U\in U(n)\) such that
\[
\operatorname{dom}D=\{f\in W:f(b)=Uf(a)\},
\]
and \(Df=(2\pi i)^{-1}f'\) \(\forall f\in\operatorname{dom}D\).
\end{lemma}

\begin{proof}
Let \(W_0\) be the space of functions in \(W\) that vanish at all \(2n\) endpoints. For \(\lambda\in\Lambda\) the function \(e_\lambda\) lies in \(W\), with \(e_\lambda'=2\pi i\lambda e_\lambda\). For \(\varphi\in W_0\), formula \eqref{eq:green} gives \(\langle(2\pi i)^{-1}\varphi',e_\lambda\rangle=\lambda\langle\varphi,e_\lambda\rangle\). Summing the squares of these coefficients over the basis shows that \(\varphi\in\operatorname{dom}D\) and \(D\varphi=(2\pi i)^{-1}\varphi'\). Thus \(D\) extends the operator \(D_0=(2\pi i)^{-1}d/dx\) with domain \(W_0\), and \(D=D^*\) is a restriction of \(D_0^*\). A function \(g\in L^2(\Omega)\) lies in the domain of \(D_0^*\) exactly when its distributional derivative on each open component lies in \(L^2\), that is, when \(g\in W\), and then \(D_0^*g=(2\pi i)^{-1}g'\). Hence \(\operatorname{dom}D\subset W\), and \(Df=(2\pi i)^{-1}f'\) for \(f\in\operatorname{dom}D\).

Let \(\omega\) be the Hermitian form on \(\C^n\oplus\C^n\) given by \(\omega((x,y),(x',y'))=\langle y,y'\rangle-\langle x,x'\rangle\). It is nondegenerate. Let \(V\) be the subspace \(\{(f(a),f(b)):f\in\operatorname{dom}D\}\). For \(f,g\in W\), formula \eqref{eq:green} reads
\begin{equation}
\label{eq:greenD}
\langle(2\pi i)^{-1}f',g\rangle-\langle f,(2\pi i)^{-1}g'\rangle=(2\pi i)^{-1}\,\omega\bigl((f(a),f(b)),(g(a),g(b))\bigr).
\end{equation}
Since \(D\) is symmetric, \(\omega\) vanishes on \(V\times V\). Conversely, suppose that \(\omega(v,(x,y))=0\) \(\forall v\in V\). Choose \(g\in W\), affine on each component, with \(g(a)=x\) and \(g(b)=y\). By \eqref{eq:greenD}, \(\langle Df,g\rangle=\langle f,(2\pi i)^{-1}g'\rangle\) \(\forall f\in\operatorname{dom}D\), so \(g\in\operatorname{dom}D^*=\operatorname{dom}D\) and \((x,y)\in V\). Thus \(V\) is its own orthogonal space for \(\omega\), and nondegeneracy gives \(\dim V=2n-\dim V\), that is, \(\dim V=n\). For \((x,y)\in V\) one has \(\omega((x,y),(x,y))=0\), so \(\lVert y\rVert=\lVert x\rVert\). The projection \((x,y)\mapsto x\) is therefore injective on \(V\), hence onto \(\C^n\), and \(V\) is the graph of a linear isometry \(U\) of \(\C^n\).

Every \(f\in\operatorname{dom}D\) satisfies \(f(b)=Uf(a)\) by construction. If \(f\in W\) satisfies \(f(b)=Uf(a)\), then \((f(a),f(b))\in V\), and \eqref{eq:greenD} together with the vanishing of \(\omega\) on \(V\times V\) gives \(\langle Dg,f\rangle=\langle g,(2\pi i)^{-1}f'\rangle\) \(\forall g\in\operatorname{dom}D\). Hence \(f\in\operatorname{dom}D^*=\operatorname{dom}D\). The space \(V\) determines \(U\), so \(U\) is unique.
\end{proof}

The matrix \(U\) is the \emph{boundary unitary} of the spectral pair \((\Omega,\Lambda)\). Every \(e_\lambda\) lies in \(\operatorname{dom}D\), so Lemma~\ref{lem:boundaryunitary} gives the boundary relation
\begin{equation}
\label{eq:boundaryrelation}
\bigl(\ex{\lambda b_j}\bigr)_{j=1}^{n}=U\bigl(\ex{\lambda a_j}\bigr)_{j=1}^{n},\qquad\lambda\in\Lambda,
\end{equation}
which is the coupling of left and right endpoint values found at the beginning of the section. The next statement shows that the coupling characterizes the spectrum: a complex number \(\zeta\) lies in \(\Lambda\) as soon as \(U\) carries some nonzero vector \(w\) to \(\bigl(\ex{\zeta\ell_j}w_j\bigr)_{j}\).

\begin{proposition}
\label{prop:secular}
Let \(\Omega\), \(\Lambda\) and \(U\) be as in Lemma~\ref{lem:boundaryunitary}. For \(\zeta\in\C\) put \(D_\ell(\zeta)=\operatorname{diag}(\ex{\zeta\ell_1},\ldots,\ex{\zeta\ell_n})\), \(M(\zeta)=D_\ell(\zeta)-U\), \(F(\zeta)=\det M(\zeta)\), and \(y_\zeta=(\ex{\zeta a_j})_{j=1}^{n}\). Then the following hold.
\begin{enumerate}
\item[\rm (i)] \(F(\zeta)=0\) iff \(\zeta\in\Lambda\). In particular every zero of \(F\) is real.
\item[\rm (ii)] For \(\lambda\in\Lambda\), the kernel of \(M(\lambda)\) is \(\C y_\lambda\), so \(M(\lambda)\) has rank \(n-1\), and
\[
\operatorname{adj}M(\lambda)=\kappa(\lambda)\,y_\lambda\bigl(D_\ell(\lambda)y_\lambda\bigr)^*
\]
for a scalar \(\kappa(\lambda)\neq 0\).
\item[\rm (iii)] For \(\lambda\in\Lambda\), \(F'(\lambda)=2\pi iL\,\kappa(\lambda)\neq 0\). Thus every zero of \(F\) is simple.
\end{enumerate}
\end{proposition}

\begin{proof}
Fix \(\zeta\in\C\). For \(w\in\C^n\) let \(f_w\) be the function equal to \(w_j\ex{\zeta(x-a_j)}\) on \([a_j,b_j)\). It lies in \(W\), with \(f_w(a)=w\), \(f_w(b)=D_\ell(\zeta)w\), and \((2\pi i)^{-1}f_w'=\zeta f_w\). If \(M(\zeta)w=0\), then \(f_w(b)=Uf_w(a)\), so Lemma~\ref{lem:boundaryunitary} puts \(f_w\) in \(\operatorname{dom}D\) with \(Df_w=\zeta f_w\). When \(w\neq 0\) the function \(f_w\) is nonzero, so \(\zeta\) is an eigenvalue of \(D\). Hence \(\zeta\in\Lambda\) and \(f_w\in\C e_\zeta\). The map \(w\mapsto f_w\) is injective, so \(F(\zeta)\neq 0\) when \(\zeta\notin\Lambda\), and \(\ker M(\zeta)\) has dimension at most one when \(\zeta\in\Lambda\). Conversely, for \(\lambda\in\Lambda\) one has \(\ex{\lambda b_j}=\ex{\lambda\ell_j}\ex{\lambda a_j}\), so the boundary relation \eqref{eq:boundaryrelation} reads \(D_\ell(\lambda)y_\lambda=Uy_\lambda\), that is, \(M(\lambda)y_\lambda=0\). This proves (i) and the description of the kernel.

Let \(\lambda\in\Lambda\). The matrix \(D_\ell(\lambda)\) is unitary, and \(M(\lambda)^*=-D_\ell(\lambda)^*M(\lambda)U^*\). Hence \(\ker M(\lambda)^*=U\ker M(\lambda)\), which is the line \(\C D_\ell(\lambda)y_\lambda\) because \(Uy_\lambda=D_\ell(\lambda)y_\lambda\). The adjugate of a matrix of rank \(n-1\) has rank one. From \(M(\lambda)\operatorname{adj}M(\lambda)=\operatorname{adj}M(\lambda)\,M(\lambda)=F(\lambda)I=0\), its columns lie in \(\ker M(\lambda)\) and the conjugate transposes of its rows lie in \(\ker M(\lambda)^*\). This gives (ii).

Jacobi's formula gives
\[
F'(\lambda)=\operatorname{tr}\bigl(\operatorname{adj}M(\lambda)\,M'(\lambda)\bigr),\qquad M'(\lambda)=2\pi i\operatorname{diag}(\ell_1,\ldots,\ell_n)\,D_\ell(\lambda).
\]
By (ii) the trace equals
\[
2\pi i\,\kappa(\lambda)\,\bigl(D_\ell(\lambda)y_\lambda\bigr)^*\operatorname{diag}(\ell_1,\ldots,\ell_n)\,D_\ell(\lambda)y_\lambda=2\pi i\,\kappa(\lambda)\sum_{j=1}^{n}\ell_j,
\]
since every coordinate of \(D_\ell(\lambda)y_\lambda\) has modulus one. This is (iii).
\end{proof}

\begin{lemma}
\label{lem:freqwidth}
Let \(\Lambda\) be a spectrum of a bounded measurable set \(\Omega\subset\R\) of measure \(L\), and let \(G(\zeta)=\sum_{\beta\in\mathcal B}g_\beta\ex{\beta\zeta}\) for a finite set \(\mathcal B\subset\R\) and complex coefficients \(g_\beta\). Suppose that \(\mathcal B\) lies in an interval of length strictly less than \(L\). If \(G\) vanishes on \(\Lambda\), then \(g_\beta=0\) \(\forall \beta\in\mathcal B\).
\end{lemma}

\begin{proof}
By Theorems~\ref{thm:periodic} and~\ref{thm:fiber}, \(\Lambda=\Gamma+h^{-1}\Z\), where the \(k\) elements of \(\Gamma\) are pairwise distinct modulo \(h^{-1}\Z\) and \(L=kh\). Call two frequencies in \(\mathcal B\) equivalent when their difference lies in \(h\Z\), and for a class \(c\) put \(G_c(\zeta)=\sum_{\beta\in c}g_\beta\ex{\beta\zeta}\). Let \(\beta_c\) be the least element of \(c\). For \(\gamma\in\Gamma\) and \(q\in\Z\), every \(\beta\in c\) satisfies \(\ex{\beta q/h}=\ex{\beta_c q/h}\), and therefore
\[
0=G(\gamma+q/h)=\sum_c\bigl(\ex{\beta_c/h}\bigr)^{q}\,G_c(\gamma).
\]
The numbers \(\ex{\beta_c/h}\) are pairwise distinct, so the sequences \(q\mapsto\bigl(\ex{\beta_c/h}\bigr)^q\) are linearly independent, and \(G_c(\gamma)=0\) for every class \(c\) and every \(\gamma\in\Gamma\). Every element of \(c\) has the form \(\beta_c+mh\) with \(m\in\Z\) and \(0\le mh<L=kh\), so \(0\le m\le k-1\). Thus \(G_c(\zeta)=\ex{\beta_c\zeta}P_c(\ex{h\zeta})\) for a polynomial \(P_c\) of degree at most \(k-1\) whose coefficients are the \(g_\beta\) with \(\beta\in c\). The polynomial \(P_c\) vanishes at the \(k\) points \(\ex{h\gamma}\), \(\gamma\in\Gamma\), which are distinct because the elements of \(\Gamma\) are distinct modulo \(h^{-1}\Z\). Hence \(P_c=0\), and every coefficient vanishes.
\end{proof}

The secular determinant of a spectral union of intervals vanishes on \(\Lambda\) and has the frequencies \(0\) and \(L\), so the bound on the width in Lemma~\ref{lem:freqwidth} cannot be relaxed to \(L\).

\section{The number of points in a fiber}
\label{sec:cosets}

Let \(\Omega = \bigsqcup_{j=1}^{4} [a_j,b_j)\) be a union of four pairwise disjoint non-adjacent bounded intervals, with \(\ell_j = b_j - a_j > 0\) and measure \(L\), and let \(\Lambda\) be a spectrum of \(\Omega\), with minimal period \(h^{-1}\) and with \(k\), \(\Gamma\), and \(S\) as in Theorem~\ref{thm:periodic}. Then \(L = kh\) by Theorem~\ref{thm:fiber}. This section bounds \(k\). The bound holds except possibly for constant fibers consisting of four integer runs of lengths \(2,2,1,1\), with \(k = 6\). This exceptional type is excluded by Lemma~\ref{lem:sixpoint}, and the resulting inequality \(k \le 4\) is Corollary~\ref{cor:k4}.

Let \(U \in U(4)\), \(D_\ell\), \(M\), \(F\), \(y_\lambda\), and \(\kappa\) be as in Lemma~\ref{lem:boundaryunitary} and Proposition~\ref{prop:secular}, with \(n = 4\). Thus \(F\) is an exponential polynomial whose zeros are exactly the points of \(\Lambda\), all real and simple, and at each \(\lambda \in \Lambda\) the adjugate of \(M(\lambda)\) is the rank-one matrix \(\kappa(\lambda)\, y_\lambda \bigl(D_\ell(\lambda) y_\lambda\bigr)^*\), with \(F'(\lambda) = 2\pi i L\, \kappa(\lambda)\). The uniqueness principle used below is Lemma~\ref{lem:freqwidth}: an exponential polynomial which vanishes on \(\Lambda\) and whose frequencies lie in an interval of length strictly less than \(L\) vanishes identically.

For \(T \subset \{1,2,3,4\}\) write \(\ell_T = \sum_{j \in T} \ell_j\), with \(\ell_\emptyset = 0\) and \(\ell_{\{1,2,3,4\}} = L\), and write \(z_T(\zeta) = \ex{\zeta \ell_T}\), abbreviating \(z_{\{j\}}\) to \(z_j\). Let \(U[T]\) be the principal minor of \(U\) on the index set \(T\), with the convention \(U[\emptyset] = 1\), and set \(c_T = (-1)^{4-\lvert T\rvert} U[T^c]\). In particular \(c_\emptyset = \det U\) and \(c_{\{1,2,3,4\}} = 1\). Let \(C_j(\zeta)\) be the \((j,j)\) cofactor of \(M(\zeta)\).

\begin{theorem}
\label{thm:cofactor}
The secular determinant and its diagonal cofactors are
\[
F(\zeta) = \sum_T c_T z_T(\zeta), \qquad z_j(\zeta)\, C_j(\zeta) = \sum_{T \ni j} c_T z_T(\zeta),
\]
and \(z_j(\zeta)\, C_j(\zeta) = F'(\zeta)/(2\pi i L)\) for every index \(j\) and every \(\zeta \in \mathbb{C}\). Equivalently, for every real \(\alpha\), the vector \(V_\alpha = \sum_{\ell_T = \alpha} c_T 1_T \in \mathbb{C}^4\) is a scalar multiple of the all-ones vector: the sum \(\sum_{\ell_T = \alpha,\, T \ni j} c_T\) does not depend on \(j\). If a proper nonempty set \(T\) is the unique subset with length sum \(\ell_T\), then \(U[T^c] = 0\). If exactly one other subset \(T'\) has the same length sum, then either \(T' = T^c\), in which case \(\ell_T = L/2\) and \(c_T = c_{T^c}\), or else \(c_T = c_{T'} = 0\).
\end{theorem}

\begin{proof}
Row \(i\) of \(M(\zeta)\) is \(z_i(\zeta)\, e_i^T - u_i\), where \(e_i\) is the \(i\)-th standard basis vector and \(u_i\) is the \(i\)-th row of \(U\). Expanding the determinant by multilinearity in the rows, the summand in which the diagonal term \(z_i e_i^T\) is chosen for \(i \in T\) and the row \(-u_i\) is chosen for \(i \notin T\) equals \(z_T(\zeta)\) times \((-1)^{\lvert T^c\rvert} U[T^c]\). This is \(c_T z_T(\zeta)\), and the first identity follows. The cofactor \(C_j\) is the secular determinant of the principal submatrix with index \(j\) removed. The same expansion, multiplied by \(z_j\), gives the second identity.

At a point \(\lambda \in \Lambda\), the diagonal entry of the adjugate is
\[
C_j(\lambda) = \kappa(\lambda)\, (y_\lambda)_j \overline{z_j(\lambda)\,(y_\lambda)_j}.
\]
Since \(\lvert(y_\lambda)_j\rvert = 1\) and \(\lambda\) is real, this equals \(\kappa(\lambda)\, \overline{z_j(\lambda)}\), and therefore
\[
z_j(\lambda) C_j(\lambda) = \kappa(\lambda) = \frac{F'(\lambda)}{2\pi i L}.
\]
The difference
\[
E_j(\zeta) = z_j(\zeta)\, C_j(\zeta) - \frac{F'(\zeta)}{2\pi i L}
\]
is an exponential polynomial. The coefficient of \(z_{\{1,2,3,4\}}\) in \(z_j C_j\) is \(1\), and in \(F'/(2\pi i L)\) it is \(L/L = 1\). The coefficient of the constant term is zero in both. Every remaining frequency is a length \(\ell_T\) with \(T\) proper and nonempty, hence lies in \([\min_i \ell_i,\, L - \min_i \ell_i]\). The frequency width of \(E_j\) is at most \(L - 2\min_i \ell_i < L\). Since \(E_j\) vanishes on \(\Lambda\), Lemma~\ref{lem:freqwidth} gives \(E_j \equiv 0\).

The four polynomials \(z_j C_j\) are therefore equal to \(F'/(2\pi i L)\). The coefficient of \(\ex{\zeta \alpha}\) in \(z_j C_j\) is the \(j\)-th coordinate of \(V_\alpha\), and in \(F'/(2\pi i L)\) it is \((\alpha/L) \sum_{\ell_T = \alpha} c_T\). Exponentials with distinct real frequencies are linearly independent, so \(V_\alpha\) equals that scalar times the all-ones vector.

Suppose \(\emptyset \neq T \neq \{1,2,3,4\}\) and no other subset has length sum \(\ell_T\). Then \(c_T 1_T\) is a multiple of the all-ones vector. The indicator \(1_T\) is not itself a multiple of that vector, so \(c_T = 0\), that is \(U[T^c] = 0\).

Suppose exactly one other subset \(T'\) has the same length sum, and write \(c_T 1_T + c_{T'} 1_{T'} = \mu \mathbf{1}\). If some index lies in neither set, then \(\mu = 0\). If \(T \cup T' \neq \{1,2,3,4\}\), or if the two sets are nested, the same vanishing follows by reading a coordinate in the symmetric difference, and then both coefficients vanish. If \(T \cup T' = \{1,2,3,4\}\) and \(T \cap T' \neq \emptyset\), a coordinate in the intersection and a coordinate in \(T \setminus T'\) force one coefficient to vanish and then the other. The remaining case is \(T' = T^c\). Both coefficients then equal \(\mu\), and \(\ell_T = \ell_{T^c}\) forces \(\ell_T = L/2\).
\end{proof}

By the complementary-minor identity (\cite{horn-johnson2013}, the Jacobi identity for the minors of \(U^{-1} = U^*\)), one has \(U[T^c] = \det(U)\, \overline{U[T]}\) \(\forall T\). In particular \(U[T] = 0\) iff \(U[T^c] = 0\), and \(\lvert c_T\rvert = \lvert c_{T^c}\rvert\).

\begin{theorem}
\label{thm:monomial}
Let \(U\) be an invertible complex \(n \times n\) matrix whose nonempty proper principal minors all vanish. Then the directed graph with an edge \(i \to j\) whenever \(u_{ij} \neq 0\) is a single \(n\)-cycle: \(u_{ij} = 0\) unless \(j = \sigma(i)\) for one cyclic permutation \(\sigma\), and every entry \(u_{i\sigma(i)}\) is nonzero. If \(U\) is unitary, each of those entries has modulus one, so \(U\) is monomial. In particular, if the sixteen subset sums \(\ell_T\) of a spectral union of four intervals are pairwise distinct, then \(k = 1\) and \(\Lambda = \lambda_0 + L^{-1}\mathbb{Z}\) for some real \(\lambda_0\).
\end{theorem}

\begin{proof}
The case \(n=1\) is immediate, so assume \(n\ge 2\). Write \(G\) for the directed graph of the nonzero entries. A directed cycle of length \(k\) is a closed directed walk through \(k\) distinct vertices; a loop has length one. There is no directed cycle of length \(k\) for any \(k = 1,\ldots,n-1\). Indeed \(U[\{i\}] = u_{ii} = 0\). Fix \(2 \le k \le n-1\) and assume there is no shorter cycle. For a set \(T\) of size \(k\), the minor \(U[T]\) is a sum over permutations of \(T\). A nonzero term comes from a permutation whose functional graph lies in \(G\). Every cycle of that permutation, fixed points included, would be a directed cycle in \(G\) of length at most \(k\), hence of length exactly \(k\), so the permutation is a Hamiltonian cycle of the induced subgraph \(G[T]\). That subgraph has at most one Hamiltonian cycle. Two distinct Hamiltonian cycles would differ by an edge, and that edge, followed by the complementary arc of the first cycle, would close a directed cycle of length between \(2\) and \(k-1\). A single Hamiltonian cycle \(\sigma_0\) would contribute \(U[T] = \operatorname{sgn}(\sigma_0) \prod_{i \in T} u_{i\sigma_0(i)} \neq 0\), contradicting the hypothesis. Thus \(G[T]\) has no \(k\)-cycle.

The determinant of \(U\) is nonzero, and by the previous paragraph only \(n\)-cycles contribute to it, so \(G\) has a Hamiltonian cycle \(H\). An edge of \(G\) outside \(H\) would close with a path of \(H\) to a directed cycle of length between \(2\) and \(n-1\). The edges of \(G\) are therefore exactly the edges of \(H\). If \(U\) is unitary, each row has norm one and a single nonzero entry, so that entry has modulus one.

Now let the sixteen sums \(\ell_T\) be pairwise distinct. Theorem~\ref{thm:cofactor} gives \(U[T] = 0\) for every proper nonempty \(T\). The matrix \(U\) is unitary, hence invertible, so it is monomial on a \(4\)-cycle \(\sigma\). In the expansion of \(\det(D_\ell(\zeta)-U)\), a permutation contributes only if it takes each index either to itself or to its successor under \(\sigma\). Choosing the successor once forces the permutation to be \(\sigma\) throughout, because \(\sigma\) is a single cycle. A \(4\)-cycle is an odd permutation, and the contribution of \(\sigma\) is \(\operatorname{sgn}(\sigma) (-1)^4 \prod_i u_{i\sigma(i)} = -\eta\) with \(\lvert\eta\rvert = 1\). The identity permutation contributes \(z_{\{1,2,3,4\}}\). Thus \(F(\zeta) = \ex{\zeta L} - \eta\), and the zero set is a single coset of \(L^{-1}\mathbb{Z}\). The minimal period \(h^{-1} = k/L\) divides \(L^{-1}\), so \(k = 1\).
\end{proof}

\begin{proposition}
\label{prop:periodsupport}
Let \(k_0\) be a positive integer. If \(\ell_T \in (L/k_0)\mathbb{Z}\) \(\forall T\) with \(c_T \neq 0\), then \(F(\zeta) = q\bigl(\ex{(L/k_0)\zeta}\bigr)\) for a monic polynomial \(q\) of degree \(k_0\) with simple unimodular roots, and \(k_0/L\) is a period of \(\Lambda\). Conversely, write \(F(\zeta) = \sum_\alpha d_\alpha \ex{\alpha\zeta}\) with \(d_\alpha = \sum_{\ell_T = \alpha} c_T\). If \(p > 0\) is a period of \(\Lambda\), then \(\ex{p\alpha} = 1\) \(\forall \alpha\) with \(d_\alpha \neq 0\). In particular \(pL \in \mathbb{Z}\).
\end{proposition}

\begin{proof}
Each surviving frequency is an integer multiple of \(h_0 = L/k_0\), and \(L = k_0 h_0\) is the largest, so \(F(\zeta) = q(\ex{h_0 \zeta})\) with \(q\) monic of degree \(k_0\) and \(q(0) = c_\emptyset = \det U \neq 0\). A root \(\rho\) of \(q\) is nonzero. If \(\lvert\rho\rvert \neq 1\), the solution of \(\ex{h_0 \zeta} = \rho\) is a non-real zero of \(F\), which Proposition~\ref{prop:secular} excludes. If \(\rho\) is a multiple root, then \(F\) has a multiple zero, because \(F'(\zeta) = 2\pi i h_0 \ex{h_0 \zeta}\, q'(\ex{h_0 \zeta})\), and multiple zeros are excluded by the same proposition. Thus \(q\) has \(k_0\) distinct unimodular roots, and \(\Lambda\) is the union of the corresponding \(k_0\) cosets of \(h_0^{-1}\mathbb{Z}\). In particular \(k_0/L\) is a period.

Conversely, if \(\Lambda + p = \Lambda\), the difference
\[
G(\zeta) = F(\zeta+p) - F(\zeta) = \sum_\alpha d_\alpha\bigl(\ex{p\alpha}-1\bigr)\ex{\alpha\zeta}
\]
vanishes on \(\Lambda\). Its constant term is zero. If \(\ex{pL} \neq 1\), the top frequency \(L\) survives in \(G\); otherwise every frequency of \(G\) lies strictly between \(0\) and \(L\). In both cases the frequencies of \(G\) lie in an interval of length strictly less than \(L\), so Lemma~\ref{lem:freqwidth} gives \(d_\alpha(\ex{p\alpha}-1) = 0\) \(\forall \alpha\). Every \(\alpha\) with \(d_\alpha \neq 0\) therefore satisfies \(\ex{p\alpha} = 1\), and \(d_L = 1\) gives \(pL \in \mathbb{Z}\). The coefficients \(c_T\) of a single frequency may cancel, so the conclusion concerns the sums \(d_\alpha\) and not the individual \(c_T\).
\end{proof}

To apply Proposition~\ref{prop:periodsupport}, it suffices to restrict the subset sums that can carry nonzero coefficients. A level is a number \(\alpha \in (0,L)\) for which the indicators \(\{1_T : \ell_T = \alpha\}\), together with the all-ones vector, are linearly dependent in \(\mathbb{C}^4\). By Theorem~\ref{thm:cofactor}, a proper subset can have \(c_T \neq 0\) only when its length sum is such a level: the coefficient vector \((c_T)_{\ell_T = \alpha}\) lies in the space of relations \(\sum \gamma_T 1_T \in \mathbb{C}\mathbf{1}\). Two subsets of equal sum are never nested, because the lengths are positive, so each level family is an antichain. Complementation exchanges the level \(\alpha\) with the level \(L-\alpha\) and preserves dependence. Indeed, if \(\sum \gamma_T 1_T = \mu \mathbf{1}\), then \(\sum \gamma_T 1_{T^c} = \bigl(\sum \gamma_T - \mu\bigr)\mathbf{1}\).

\begin{theorem}
\label{thm:levelclass}
The levels of four positive lengths are contained in
\[
\{L/2,\, L/3,\, 2L/3,\, L/4,\, 3L/4,\, 2L/5,\, 3L/5\}.
\]
Exactly one of the following holds.
\begin{enumerate}
\item The lengths are proportional to a permutation of \((1,1,1,2)\), and the levels are \(2L/5\) and \(3L/5\).
\item The lengths are proportional to a permutation of \((1,1,2,2)\), and the levels are \(L/3\), \(L/2\), and \(2L/3\).
\item The levels are contained in \((L/k_0)\mathbb{Z}\) for some integer \(k_0 \in \{1,2,3,4\}\). All four lengths equal gives \(k_0 = 4\). A partition of the four indices into three blocks of equal length sum, with no subset summing to \(L/2\), gives \(k_0 = 3\). A subset summing to \(L/2\), with no such three-block partition, gives \(k_0 = 2\). The absence of every level gives \(k_0 = 1\).
\end{enumerate}
In the third case, Proposition~\ref{prop:periodsupport} supplies the period \(k_0/L\). The minimal period \(k/L\) divides it, so \(k \le k_0 \le 4\).
\end{theorem}

\begin{proof}
The organizing principle is the size of the level family, after complementation has reduced the claim to levels \(\alpha \le L/2\). A family of size one is independent from the all-ones vector, so it contributes no level. A family of size two is dependent with the all-ones vector iff the two sets are complementary, by the second vanishing rule of Theorem~\ref{thm:cofactor}, and then \(\alpha = L/2\).

Consider a family of three sets \(T_1,T_2,T_3\) and a relation \(\gamma_1 1_{T_1} + \gamma_2 1_{T_2} + \gamma_3 1_{T_3} = \mu \mathbf{1}\) in which not every coefficient is zero. If \(\mu = 0\), no coefficient vanishes, because two distinct nonempty \(0/1\) vectors are linearly independent. Solving for the first indicator gives \(1_{T_1} = a\, 1_{T_2} + b\, 1_{T_3}\) with \(a, b \neq 0\). Comparing coordinates yields \(T_1 \subset T_2 \cup T_3\), and the same inclusion with the indices cycled. The scalar \(a\) equals \(1\) at each point of \(T_1 \cap T_2 \setminus T_3\), and \(b\) equals \(1\) at each point of \(T_1 \cap T_3 \setminus T_2\), while \(a+b = 0\) on \(T_2 \cap T_3 \setminus T_1\) and \(a+b = 1\) on the triple intersection. Since \(T_1\) sits in neither of the other sets, both of those first two regions are nonempty, so \(a = b = 1\). Then \(a+b = 2\), the two regions that would require the value \(0\) or \(1\) are empty, and \(T_2 \cap T_3 = \emptyset\). Combined with \(T_2 \subset T_1 \cup T_3\), this forces \(T_2 \subset T_1\), and likewise \(T_3 \subset T_1\), so \(T_1 = T_2 \sqcup T_3\). The length of \(T_1\) would be \(2\alpha\). Thus \(\mu \neq 0\), and one may take \(\mu = 1\).

A zero coefficient would leave a complementary pair, hence \(\alpha = L/2\), where the family is closed under complementation and therefore has even size. All three coefficients are nonzero. Every coordinate lies in at least one set. The incidence type of \(j \in \{1,2,3,4\}\) is \(\{r \in \{1,2,3\}:j\in T_r\}\); the coefficients indexed by this type sum to \(1\). If three singleton types occur, the coefficients are \((1,1,1)\), no other type agrees with that value, and the three sets partition \(\{1,2,3,4\}\). This is a three-block partition of equal sum, at level \(L/3\). If exactly two singleton types occur, a proper pair type would kill a coefficient and the doubleton formed by those two indices would require the sum \(2\), so every remaining coordinate has the triple type and the third set sits inside the intersection of the first two. If exactly one singleton type occurs, the remaining coordinates share one type, either the complementary pair of indices or all three indices, and in both subcases the other two sets coincide. If no singleton type occurs and the triple type occurs, each pair type kills a coefficient, so every coordinate has the triple type and the three sets coincide. If only pair types occur and one of them is missing, two of the sets are nested or equal. If all three pair types occur, the coefficients are \(\bigl(\tfrac12,\tfrac12,\tfrac12\bigr)\) and the sets are the complements of a three-block partition, so the level is \(2L/3\), already reduced to \(L/3\) by complementation. The only three-set level at most \(L/2\) is therefore the three-block partition at \(L/3\), and the complementary family sits at \(2L/3\).

Consider a family of at least four sets at a level \(\alpha < L/2\). A three-set \(T = \{1,2,3,4\}\setminus\{i\}\) with \(\ell_T = \alpha < L/2\) would require \(\ell_i > L/2 > \alpha\). Every other member of the antichain would then contain \(i\), and its sum would exceed \(\alpha\). A family of size at least two at such a level therefore contains no three-set, and consists of singletons and pairs. Three or more singletons, with no room for a disjoint pair, force all four lengths to equal \(\alpha\), hence \(\alpha = L/4\) and the family of all four singletons; its complement is the family of all four triples, at level \(3L/4\). Exactly two singletons admit only the opposite pair, and the family then has size at most three, already treated. Exactly one singleton \(\{i\}\) forces the three pairs in the complement of \(\{i\}\) into the family, hence three equal lengths \(s\) and \(\ell_i = 2s\), so \(\alpha = 2s = 2L/5\). The lengths are a permutation of \((s,s,s,2s)\), and the complementary level is \(3L/5\). A family of pairs alone, of size at least four and below \(L/2\), is impossible: the six pairs form three complementary couples, and two members of one couple already sum to \(L\), so at most three pairs can share a sum strictly below \(L/2\).

At level \(L/2\) the family is a union of complementary couples. A singleton of length \(L/2\) excludes every pair, and a second singleton would force a zero length, so the family is that singleton and its complement. Otherwise the family consists of complementary couples of pairs. Two such couples force the length pattern \((a,b,b,a)\) up to order, and the third couple joins iff \(a = b\). In every subcase the relation \(1_T + 1_{T^c} = \mathbf{1}\) shows that the level is genuine.

These are all the possibilities, so every level lies in the stated set. If \(2L/5\) is a level, the lengths are a permutation of \((s,s,s,2s)\). Every subset sum is then an integer multiple of \(s\), while \(L/2 = 5s/2\), \(L/3 = 5s/3\), and \(L/4 = 5s/4\) are not, so the only levels are \(2L/5\) and \(3L/5\). This is the first case. If \(L/4\) is a level, all four lengths are equal and the levels are \(L/4\), \(L/2\), and \(3L/4\), which is the third case with \(k_0 = 4\).

Otherwise every level lies in \(\{L/2, L/3, 2L/3\}\). A single denominator \(k_0 \le 4\) fails only when both \(L/2\) and \(L/3\) occur. A three-block partition of equal sum has, up to order, two blocks of one index and one block of two, so the lengths are \((s,s,t,u)\) with \(s = t+u = L/3\). A subset sum equal to \(L/2 = 3s/2\) cannot be a single length, nor the sum of the two short lengths \(t+u = s\), nor the sum \(2s\) of the two long lengths. It must be \(s+t\) or \(s+u\), hence \(t = u = s/2\). The lengths are proportional to \((2,2,1,1)\), the levels are \(L/3\), \(L/2\), and \(2L/3\), and this is the second case. If \(L/2\) and \(L/3\) do not occur together, the third case holds with \(k_0 = 2\), \(k_0 = 3\), or \(k_0 = 1\).

In the third case every subset sum \(\ell_T\) with \(c_T \neq 0\) lies in \((L/k_0)\mathbb{Z}\), because a nonzero proper coefficient can occur only at a level and the levels lie on that lattice, as do \(0\) and \(L\). Proposition~\ref{prop:periodsupport} gives the period \(k_0/L\). The group of periods is \((k/L)\mathbb{Z}\), so \(k\) divides \(k_0\).
\end{proof}

The first case is settled by the kernel description in Proposition~\ref{prop:secular}, and the second by the shape of \(F\).

\begin{theorem}
\label{thm:cell1112}
Suppose the lengths are \((s,s,s,2s)\) up to order. Then \(U\) is monomial on a \(4\)-cycle, \(\Lambda = \lambda_0 + L^{-1}\mathbb{Z}\) for some real \(\lambda_0\), and \(k = 1\).
\end{theorem}

\begin{proof}
Label the intervals so that \(\ell_1 = \ell_2 = \ell_3 = s\) and \(\ell_4 = 2s\), so that \(L = 5s\). The level analysis of Theorem~\ref{thm:levelclass} allows nonzero coefficients for proper nonempty subsets only at \(2s = 2L/5\) and \(3s = 3L/5\). At level \(s\) the three short singletons form a family whose indicators, with a zero in the fourth coordinate, are independent from the all-ones vector, so the corresponding coefficients vanish. Thus the complementary minors vanish as well: \(u_{11} = u_{22} = u_{33} = 0\). At level \(3s\) the family is the triple of short intervals together with the three pairs consisting of one short interval and the long interval. The level identity equates the three pairwise sums of coefficients and forces the three pair-minors on the short intervals to be equal, say to \(\gamma'\), and the diagonal entry \(u_{44}\) to equal \(-2\gamma'\). At level \(2s\) the same reading, or complementation, gives \(U[\{1,4\}] = U[\{2,4\}] = U[\{3,4\}] = \gamma\) and \(U[\{1,2,3\}] = -2\gamma\) for a scalar \(\gamma\). The complementary-minor identity gives \(\lvert\gamma\rvert = \lvert\gamma'\rvert\). If \(\gamma' = 0\), then \(\gamma = 0\), every nonempty proper principal minor vanishes, and Theorem~\ref{thm:monomial} gives the conclusion.

Suppose \(\gamma' \neq 0\). We show that \(U\) then has an eigenvector supported on the three short intervals. Since their lengths are equal, such a vector gives a secular kernel vector with zero fourth coordinate, contrary to Proposition~\ref{prop:secular}. Set \(g = \lvert\gamma'\rvert > 0\). Write \(X\) for the \(3 \times 3\) principal block on the short intervals, \(p\) for the column \((u_{14},u_{24},u_{34})^T\), and \(q\) for the column \((u_{41},u_{42},u_{43})^T\). The surviving minors read \(x_{ii} = 0\), \(x_{ij}x_{ji} = -\gamma'\) for the three short pairs, and \(p_i q_i = -\gamma\). In particular every \(p_i\) and \(q_i\) is nonzero. Set \(P_i = \lvert p_i\rvert^2\) and \(a_{ij} = \lvert x_{ij}\rvert^2\). Unitarity of the rows and columns gives \(\sum_i P_i = 1 - \lvert u_{44}\rvert^2 = 1 - 4g^2\), \(\sum_k a_{ik} = 1 - P_i\), and \(\sum_k a_{ki} = 1 - g^2/P_i\). The off-diagonal block equations \((XX^*)_{ij} = -p_i\overline{p_j}\) and \((X^*X)_{ij} = -\overline{q_i}q_j\), read through the vanishing diagonal of \(X\), give the six product relations
\begin{align*}
a_{13}a_{23} &= P_1P_2, & a_{12}a_{32} &= P_1P_3, & a_{21}a_{31} &= P_2P_3,\\
a_{31}a_{32} &= g^4/(P_1P_2), & a_{21}a_{23} &= g^4/(P_1P_3), & a_{12}a_{13} &= g^4/(P_2P_3),
\end{align*}
together with \(a_{ij}a_{ji} = g^2\). The product of the first three relations and the three evaluations \(a_{ij}a_{ji} = g^2\) yield \(P_1P_2P_3 = g^3\).

Set \(B = a_{12}\). Solving the products produces \(a_{21} = g^2/B\), \(a_{13} = gP_1/B\), \(a_{23} = P_2 B/g\), \(a_{31} = gB/P_1\), and \(a_{32} = P_1P_3/B\). The norm of the first column, multiplied by \(P_1/g\), has the same left-hand side as the norm of the first row and gives \(1 - P_1 = P_1/g - g\), hence \(P_1 = g\). The second row and the third column then have the same left-hand side, so \(P_2 = g^2/P_3\). The third row and the second column reduce to one quadratic equation \(B + gP_3/B = 1 - P_3\), while the first row is \(B + g^2/B = 1 - g\). Subtracting the cleared equations gives \((g-P_3)(B+g) = 0\), hence \(P_3 = g\) and then \(P_2 = g\). The mass relation \(3g = 1 - 4g^2\) has the positive root \(g = 1/4\). Thus \(\lvert p_i\rvert = 1/2\) and \(\lvert u_{44}\rvert = 1/2\).

From \(q_i = -\gamma/p_i\) and \(\lvert p_i\rvert^2 = 1/4\) one gets \(\overline{q_i} = -4\overline{\gamma}\, p_i\). The conjugate of the column \(q\) is a scalar multiple of \(p\). The top-left block of \(U^*U = I\) reads \(X^*X + \overline{q}\, q^T = I_3\), and the top-right block reads \(X^*p = -u_{44}\overline{q}\), so \(X^*p\) is a multiple of \(p\). For every \(v\) orthogonal to \(p\), the correction \(\overline{q}\,(q^T v)\) is orthogonal to \(v\), and therefore \(\lVert Xv\rVert = \lVert v\rVert\). Take the inner product \(\langle u,v\rangle = \sum_i u_i\overline{v_i}\). The same block equation gives \(\langle Xv, p\rangle = \langle v, X^*p\rangle = 0\), so \(X\) preserves the Hermitian orthogonal complement of \(p\). Its restriction to that two-dimensional space is a linear isometry, hence unitary, and has a unimodular eigenvalue \(\mu\) with eigenvector \(v\) orthogonal to \(p\). Since \(\overline{q}\) is a scalar multiple of \(p\) and \(\langle p,v\rangle = 0\), one has \(\langle\overline{q}, v\rangle = 0\). This inner product equals \(\overline{q^T v}\), so \(q^T v = 0\). The block equation \(X^*Xv = v - \overline{q}\,(q^T v)\) therefore gives \(X^*Xv = v\), and \(Xv = \mu v\) gives \(U(v,0)^T = \mu(v,0)^T\).

Choose a real \(\zeta_0\) with \(\ex{s\zeta_0} = \mu\). Then \(D_\ell(\zeta_0) = \operatorname{diag}(\mu,\mu,\mu,\mu^2)\), and \((v,0)^T\) lies in the kernel of \(M(\zeta_0)\). Thus \(F(\zeta_0) = 0\), so \(\zeta_0 \in \Lambda\), and Proposition~\ref{prop:secular} says that this kernel is spanned by the unimodular vector \(y_{\zeta_0}\). A kernel vector with a zero coordinate cannot be a scalar multiple of a unimodular vector. This contradiction forces \(\gamma' = 0\), and the first paragraph applies.
\end{proof}

\begin{proposition}
\label{prop:cell1122}
Suppose the intervals are labeled so that \(\ell_1 = \ell_2 = s\) and \(\ell_3 = \ell_4 = 2s\), and set \(L = 6s\) and \(w = \ex{s\zeta}\). Let \(\beta = u_{12}u_{21}\), \(\gamma_1 = U[\{1,3\}]\), and \(\gamma_2 = U[\{1,4\}]\). Then \(U[\{2,4\}] = \gamma_1\), \(U[\{2,3\}] = \gamma_2\), and
\[
F(\zeta) = w^6 - 3\beta w^4 + 2(\gamma_1+\gamma_2)w^3 - 3\det(U)\,\overline{\beta}\, w^2 + \det U.
\]
If \(\beta = 0\), then \(F\) is a polynomial in \(w^3\) and \(k \le 2\). If \(\gamma_1+\gamma_2 = 0\), then \(F\) is a polynomial in \(w^2\) and \(k \le 3\). If neither coefficient vanishes, then \(k = 6\), one has \(\Lambda = \Gamma + s^{-1}\mathbb{Z}\) with \(\lvert\Gamma\rvert = 6\), and all left endpoints are congruent modulo \(s\mathbb{Z}\).
\end{proposition}

\begin{proof}
Label the lengths \(\ell_1 = \ell_2 = s\) and \(\ell_3 = \ell_4 = 2s\). Theorem~\ref{thm:levelclass} places the levels at \(L/3 = 2s\), \(L/2 = 3s\), and \(2L/3 = 4s\). At level \(s\) the only subsets are the short singletons \(\{1\}\) and \(\{2\}\). The third coordinate of \(V_s\) is zero, so both coefficients vanish: \(u_{11} = u_{22} = 0\), and the complementary triple minors \(U[\{2,3,4\}]\) and \(U[\{1,3,4\}]\) vanish with them. Set \(\beta = u_{12}u_{21}\). The short principal minor is then \(U[\{1,2\}] = -\beta\).

At level \(4s\) the subsets are \(\{1,2,3\}\), \(\{1,2,4\}\) and \(\{3,4\}\), with coefficients \(-u_{44}\), \(-u_{33}\) and \(U[\{1,2\}] = -\beta\). The level identity equates \(-u_{33}-u_{44}\), \(-u_{44}-\beta\) and \(-u_{33}-\beta\), so \(u_{33} = u_{44} = \beta\), and each of the three coefficients equals \(-\beta\). The coefficient of \(w^4\) is therefore \(-3\beta\). At level \(2s\) the subsets are the complements \(\{4\}\), \(\{3\}\) and \(\{1,2\}\). The same identity equates their three coefficients, and the complementary-minor identity gives \(U[\{3,4\}] = \det(U)\,\overline{U[\{1,2\}]} = -\det(U)\,\overline{\beta}\). Each coefficient equals this value, so the coefficient of \(w^2\) is \(-3\det(U)\,\overline{\beta}\).

At level \(3s\) the subsets are the pairs \(\{1,3\}\), \(\{1,4\}\), \(\{2,3\}\) and \(\{2,4\}\), in two complementary couples. The level identity is \(c_T = c_{T^c}\), which reads \(U[\{1,3\}] = U[\{2,4\}]\) and \(U[\{1,4\}] = U[\{2,3\}]\). Call these common values \(\gamma_1\) and \(\gamma_2\). Each couple contributes twice, so the coefficient of \(w^3\) is \(2(\gamma_1+\gamma_2)\). The constant term is \(\det U\) and the leading term is \(w^6\). No other exponent occurs. This is the displayed polynomial.

If \(\beta = 0\), then \(F\) is a polynomial in \(w^3 = \ex{(L/2)\zeta}\), so Proposition~\ref{prop:periodsupport} gives a period \(2/L\) and \(k \le 2\). If \(\gamma_1+\gamma_2 = 0\), then \(F\) is a polynomial in \(w^2 = \ex{(L/3)\zeta}\), and \(k \le 3\). Suppose that both coefficients are nonzero. Every \(\ell_T\) lies in \(s\mathbb{Z} = (L/6)\mathbb{Z}\), so Proposition~\ref{prop:periodsupport} gives the period \(s^{-1}\). The frequencies \(4s\) and \(3s\) of \(F\) carry the nonzero coefficients \(-3\beta\) and \(2(\gamma_1+\gamma_2)\), so the converse in Proposition~\ref{prop:periodsupport} gives \(4sp \in \mathbb{Z}\) and \(3sp \in \mathbb{Z}\), hence \(sp \in \mathbb{Z}\), for every period \(p > 0\). Thus \(s^{-1}\) is the minimal period, \(h = s\), \(k = L/s = 6\), and \(\Lambda = \Gamma + s^{-1}\mathbb{Z}\) with \(\lvert\Gamma\rvert = 6\). Since every \(\ell_j\) lies in \(s\mathbb{Z}\), one has \(D_\ell(\lambda + s^{-1}) = D_\ell(\lambda)\). The kernel of \(M(\lambda+s^{-1})\) is therefore the kernel of \(M(\lambda)\), and Proposition~\ref{prop:secular} puts \(y_{\lambda+s^{-1}}\) in the line \(\mathbb{C} y_\lambda\). The coordinates give \(\ex{a_j/s}\) independent of \(j\), so \(a_j - a_1 \in s\mathbb{Z}\).
\end{proof}

\begin{proposition}
\label{prop:cosetbound}
Let \(\Omega\) be a spectral union of four intervals, with minimal period \(h^{-1}\) and with \(k = L/h\). Then \(k \le 4\), unless \(k = 6\) and, after translation and scaling, \(\Omega = A + [0,1)\), where \(A \subset \Z\) is a union of four integer runs of lengths \(2,2,1,1\) in some order and is a spectrum of the six-element set \(S = h\Gamma\), which has no nontrivial translation symmetry.
\end{proposition}

\begin{proof}
Theorem~\ref{thm:levelclass} places the lengths in one of its three cases. The third case gives \(k \le 4\). The first case gives \(k = 1\) by Theorem~\ref{thm:cell1112}. In the second case, Proposition~\ref{prop:cell1122} gives \(k \le 3\), or else \(k = 6\) and, after a translation and a scaling by \(1/s\), the set is \(A + [0,1)\) with \(A \subset \mathbb{Z}\) a union of four runs of lengths \(1,1,2,2\). Theorem~\ref{thm:fiber}, in the direction from a spectrum of \(\Omega\) to its fibers, says that \(A\) is a spectrum of \(S = s\Gamma\). A nontrivial translation symmetry of \(S\) would be a period of \(\Lambda\) strictly smaller than \(s^{-1}\), against minimality.
\end{proof}

Lemma~\ref{lem:sixpoint} excludes this six-point fiber, and Corollary~\ref{cor:k4} records the resulting bound \(k\le 4\).

\section{Lattice multitiles with at most five cosets}
\label{sec:four}

Let \(E\subset\R\) be a bounded measurable set with a spectrum \(\Gamma+h^{-1}\Z\) made of \(k\) cosets. By Theorem~\ref{thm:fiber}, almost every fiber of \(E\) is a spectrum of the \(k\)-point set \(S=h\Gamma\), and \(E\) is a \(k\)-fold multitile by \(h\Z\). It is enough to find one set \(B\subset\Z\) that is a tiling complement of almost every fiber at once. Kolountzakis and Matolcsi proved that a spectral set of at most five elements in a finite abelian group, or in \(\Z^d\), tiles \cite[Theorems~2.1 and~2.3]{km2006}. Applied to the fibers one at a time, their theorem gives a complement for each fiber separately. Dutkay and Haussermann proved that a common complement exists for \(k\le 5\), and hence that every bounded spectral set with such a spectrum tiles the real line \cite{dutkay-haussermann2015}. Shi's fiber lemma \cite[Lemma~7.1]{shi2026} gives the cases \(k\in\{2,3\}\) in the present notation.

We give an explicit construction of a common complement depending only on \(S\). For \(k\in\{2,3,5\}\) every complex Hadamard matrix of order \(k\) is equivalent to the Fourier matrix, elementarily for \(k\le 3\) and by a theorem of Haagerup \cite{haagerup1997} for \(k=5\). As in the proof of \cite[Theorem~2.1]{km2006}, every integer spectrum of \(S\) then reduces to a coset of one cyclic subgroup, and an interval of residues is a common complement (Theorem~\ref{thm:primecomplement}). For \(k=4\) the zero sums of four unimodular numbers have a continuous parameter, and the Fourier form is not available. Two facts replace it. Such a zero sum splits into two antipodal pairs (Lemma~\ref{lem:antipodal}), which makes the nonzero differences of \(S\) rational of even denominator (Lemma~\ref{lem:fournodulus}). After clearing the denominators, a dephased Hadamard matrix of order four has a column with entries \(\pm 1\). Consequently, modulo their least common multiple, every integer spectrum of \(S\) is a union of two cosets of a subgroup of order two (Theorem~\ref{thm:znfour}). One complement again serves every spectrum (Theorem~\ref{thm:complement}), and Theorem~\ref{thm:fiberfew} turns the common complement into a tiling of the real line.

\begin{lemma}
\label{lem:antipodal}
Let \(z_1,z_2,z_3,z_4\) be complex numbers of modulus \(1\) with \(z_1+z_2+z_3+z_4=0\). If \(z_1+z_2\neq 0\), then \(\{z_3,z_4\}=\{-z_1,-z_2\}\) as multisets. In every case the four numbers can be partitioned into two antipodal pairs.
\end{lemma}

\begin{proof}
Put \(w=z_1+z_2=-(z_3+z_4)\). If \(w=0\), then \(z_2=-z_1\) and \(z_4=-z_3\), which is a partition into two antipodal pairs. Suppose \(w\neq 0\). The conditions \(|z|=1\) and \(|w-z|=1\) describe the intersection of two distinct circles of radius \(1\), so there are at most two such \(z\). Both \(z_1\) and \(z_2\) lie on that intersection, and both \(z_3\) and \(z_4\) lie on the corresponding intersection for \(-w\), which is the reflection of the first through the origin. If \(z_1\neq z_2\), the first intersection is exactly \(\{z_1,z_2\}\), so \(\{z_3,z_4\}\subset\{-z_1,-z_2\}\). These two numbers are distinct: \(z_3=z_4\) would force \(|w|=2\), hence \(z_1=z_2\). If \(z_1=z_2\), then \(|w|=2\), the circles are tangent, and \(z_3=z_4=-z_1\). In both subcases, pairing \(z_1\) and \(z_2\) with their negatives among \(z_3\) and \(z_4\) partitions the four numbers into two antipodal pairs.
\end{proof}

\begin{figure}[ht]
\centering
\begin{overpic}[width=0.55\textwidth]{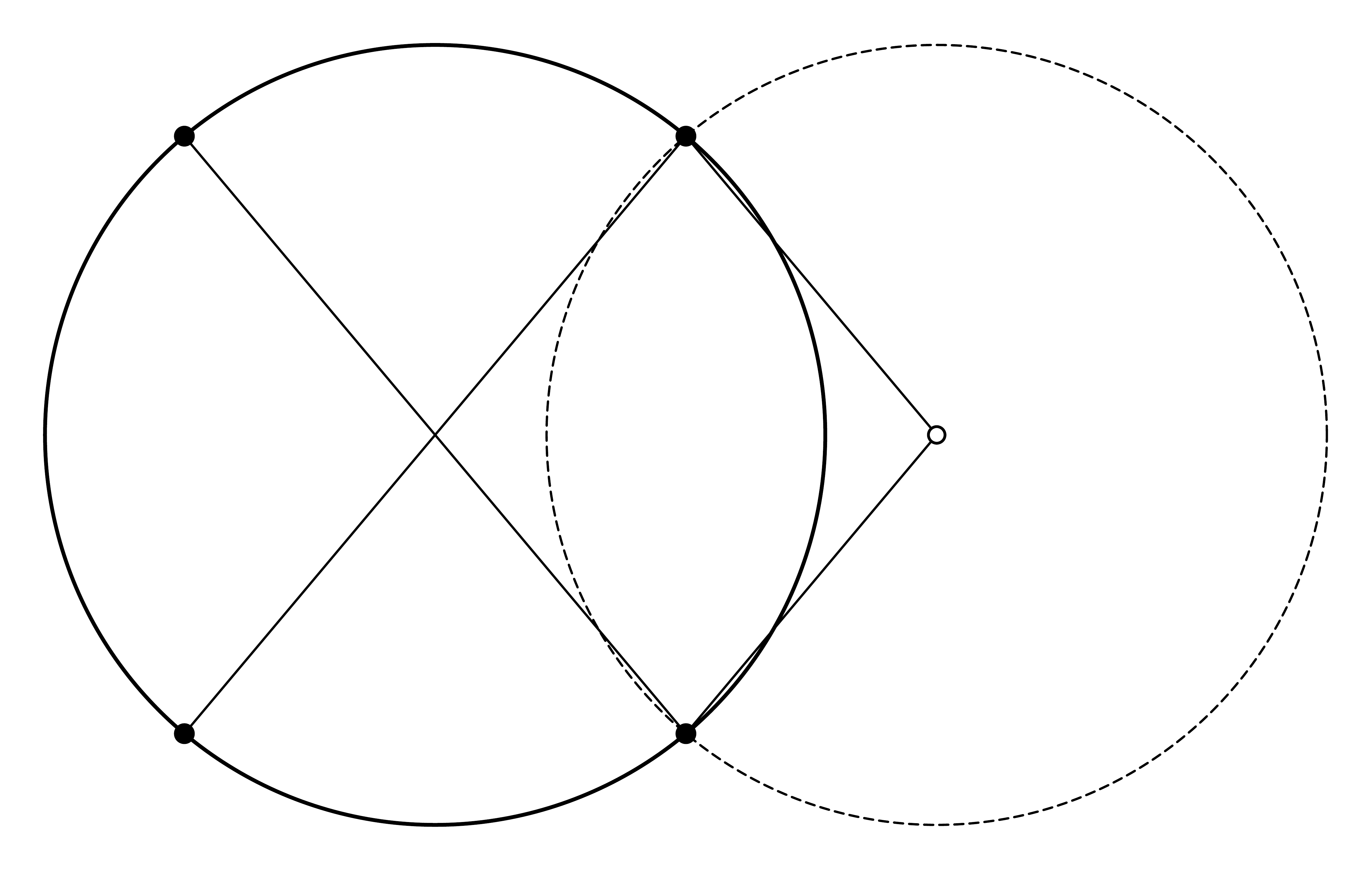}
\put(70.30,31.72){\makebox(0,0)[l]{$w$}}
\put(29.70,31.72){\makebox(0,0)[r]{$0$}}
\put(50.00,4.88){\makebox(0,0){$z_1$}}
\put(50.00,58.55){\makebox(0,0){$z_2$}}
\end{overpic}
\caption{Four points of modulus one with zero sum and \(z_1+z_2\neq 0\). The unit circle and its translate by \(w=z_1+z_2\), dashed, meet in \(z_1\) and \(z_2\). The rhombus with vertices \(0\), \(z_1\), \(w\), \(z_2\) arranges the four vectors head to tail, after interchanging \(z_3\) and \(z_4\) if necessary, and the two diameters join the antipodal pairs.}
\label{fig:antipodes}
\end{figure}

Figure~\ref{fig:antipodes} shows the two circles of the proof and the two pairs.

Write \(\ZN{N}=\Z/N\Z\) and \(\zroot{N}=e^{2\pi i/N}\) for each positive integer \(N\). For a finite set \(S'\subset\ZN{N}\), a set \(A_0\subset\ZN{N}\) is a spectrum of \(S'\) when \(\lvert A_0\rvert=\lvert S'\rvert\) and the matrix \(\bigl(\ex{la/N}\bigr)_{l\in S',\,a\in A_0}\) is a complex Hadamard matrix.

\begin{lemma}
\label{lem:reduction}
Let \(S\subset\R/\Z\) be a finite set with at least two elements and with \(S-S\subset\Q/\Z\). Let \(N\) be the least common multiple of the denominators in lowest terms of the nonzero elements of \(S-S\), and fix \(s_0\in S\). Then \(S-S\subset(1/N)\Z/\Z\), and \(S'=N(S-s_0)\subset\ZN{N}\) has \(\lvert S\rvert\) elements, contains \(0\), and satisfies \(\gcd(N,S')=1\). Reduction modulo \(N\) is injective on every spectrum \(A\subset\Z\) of \(S\), and the image is a spectrum of \(S'\). Conversely, a subset of \(\Z\) with \(\lvert S\rvert\) elements whose reduction modulo \(N\) is a spectrum of \(S'\) is a spectrum of \(S\).
\end{lemma}

\begin{proof}
The rule \(r/N+\Z\mapsto r\bmod N\) is an isomorphism of \((1/N)\Z/\Z\) with \(\ZN{N}\), so \(S'\) is a well-defined subset of \(\ZN{N}\) with \(\lvert S\rvert\) elements, and it contains \(0\). Since \(0\in S'\), the sets \(S'\) and \(S'-S'\) generate the same subgroup of \(\ZN{N}\), and \(\gcd(N,S')=\gcd(N,S'-S')\). The reduced denominator of \(r/N\) is \(N/\gcd(N,r)\). The least common multiple of these denominators, as \(r\) runs through \((S'-S')\setminus\{0\}\), equals \(N\) by the choice of \(N\). If \(d_r=\gcd(N,r)\), then \(\operatorname{lcm}_r(N/d_r)=N/\gcd_r d_r\), and \(\gcd_r d_r=\gcd(N,S'-S')\). Hence that gcd is \(1\).

The matrix \(\bigl(\ex{la/N}\bigr)_{l\in S',\,a\in A}\) is obtained from \(\bigl(\ex{sa}\bigr)_{s\in S,\,a\in A}\) by multiplying column \(a\) by \(\ex{-s_0 a}\), so its rows are orthogonal. Orthogonal columns are distinct, and the column indexed by \(a\) depends only on \(a\bmod N\). Reduction modulo \(N\) is therefore injective on \(A\), and the image is a spectrum of \(S'\). The converse is the same comparison of matrices.
\end{proof}

\begin{lemma}
\label{lem:fournodulus}
Let \(S\subset\R/\Z\) with \(\lvert S\rvert=4\), and suppose some \(A\subset\Z\) is a spectrum of \(S\). Then every \(\sigma\in(S-S)\setminus\{0\}\) is rational, and the denominator of \(\sigma\) in lowest terms is even. In particular the integer \(N\) of Lemma~\ref{lem:reduction} is even.
\end{lemma}

\begin{proof}
Let \(\sigma\in(S-S)\setminus\{0\}\) and fix \(a\in A\). Then \(\widehat A(\sigma)=0\). Lemma~\ref{lem:antipodal}, applied to the four numbers \(\ex{\sigma a''}\) with \(a''\in A\), pairs \(\ex{\sigma a}\) with its negative \(\ex{\sigma a'}\) for some \(a'\neq a\) in \(A\), so \(\sigma(a-a')\in\tfrac12+\Z\). Thus \(\sigma=(2j+1)/(2d)\) for \(d=a-a'\neq 0\) and some \(j\in\Z\). The greatest common divisor \(\gcd(2j+1,2d)\) equals \(\gcd(2j+1,d)\) and is odd, so the factor \(2\) remains in the reduced denominator. The least common multiple \(N\) of these finitely many denominators is therefore even.
\end{proof}

The notation \(A\oplus B=G\) means that every element of \(G\) has a unique representation \(a+b\), with \(a\in A\) and \(b\in B\).

\begin{lemma}
\label{lem:pullback}
Let \(\varphi\colon G\to G'\) be a homomorphism of abelian groups, let \(A\subset G\) be a finite set on which \(\varphi\) is injective, and let \(C'\subset G'\) satisfy \(\varphi(A)\oplus C'=G'\). Then \(A\oplus\varphi^{-1}(C')=G\).
\end{lemma}

\begin{proof}
Take \(g\in G\) and write \(\varphi(g)=\varphi(a)+c'\) with \(a\in A\) and \(c'\in C'\). Then \(c=g-a\) lies in \(\varphi^{-1}(C')\), and \(g=a+c\). If \(a_1+c_1=a_2+c_2\) with \(a_i\in A\) and \(c_i\in\varphi^{-1}(C')\), the images give two representations of the same element of \(G'\) as a sum from \(\varphi(A)\) and \(C'\). Uniqueness forces \(\varphi(a_1)=\varphi(a_2)\), hence \(a_1=a_2\) by injectivity on \(A\), and then \(c_1=c_2\). Surjectivity of \(\varphi\) is not required.
\end{proof}

Two complex Hadamard matrices are equivalent when one is obtained from the other by permuting rows, permuting columns, and multiplying rows and columns by unimodular constants. The Fourier matrix of order \(p\) is \(F_p=\bigl(\zroot{p}^{mj}\bigr)_{m,j\in\ZN{p}}\).

\begin{theorem}
\label{thm:primecomplement}
Let \(p\in\{2,3,5\}\), let \(S\subset\R/\Z\) have \(p\) elements, and suppose some \(A\subset\Z\) is a spectrum of \(S\). Then \(S-S\subset\Q/\Z\), and \(p\) divides the denominator in lowest terms of every nonzero element of \(S-S\). Let \(N\) be as in Lemma~\ref{lem:reduction}, and put \(K=N/p\). Every spectrum \(A\subset\Z\) of \(S\) reduces modulo \(N\) onto a coset of the subgroup \(K\ZN{N}\) of order \(p\), and \(A\oplus B=\Z\) for
\[
B=\{g\in\Z:\ g\bmod N\in[0,K)\}.
\]
\end{theorem}

For \(p\in\{2,3\}\) the coset structure and the complement \(B\) appear in the proof of Shi's fiber lemma \cite[Lemma~7.1]{shi2026}. The argument below treats the three primes together.

\begin{proof}
Every complex Hadamard matrix of order \(p\in\{2,3,5\}\) is equivalent to \(F_p\). For \(p=2\) this is immediate. For \(p=3\), multiply rows and columns so that the first row and the first column consist of ones. A unimodular row \((1,a,b)\) orthogonal to \((1,1,1)\) satisfies \(1+a+b=0\), hence \(\operatorname{Re}a=-1/2\) and \(\{a,b\}=\{\zroot{3},\zroot{3}^2\}\). After a permutation of the last two columns the second row is \((1,\zroot{3},\zroot{3}^2)\), and orthogonality to it forces the third row to be \((1,\zroot{3}^2,\zroot{3})\). For \(p=5\) the equivalence is a theorem of Haagerup \cite{haagerup1997}, which is also the input of \cite[Theorem~2.1]{km2006}.

Let \(A\) be a spectrum of \(S\), and fix \(s_0\in S\) and \(a_0\in A\). The matrix \(H=\bigl(\ex{(s-s_0)(a-a_0)}\bigr)_{s\in S,\,a\in A}\) is obtained from \(\bigl(\ex{sa}\bigr)\) by multiplying rows and columns by unimodular constants, so it is a complex Hadamard matrix, and its row \(s_0\) and its column \(a_0\) consist of ones. By the equivalence there are bijections \(\pi_1\colon S\to\ZN{p}\) and \(\pi_2\colon A\to\ZN{p}\) and unimodular functions \(\delta_1\), \(\delta_2\) with \(H_{s,a}=\delta_1(s)\,\zroot{p}^{\pi_1(s)\pi_2(a)}\,\delta_2(a)\). Since \(H_{s,a_0}=H_{s_0,a}=H_{s_0,a_0}=1\), the constants cancel in
\[
H_{s,a}=\frac{H_{s,a}\,H_{s_0,a_0}}{H_{s,a_0}\,H_{s_0,a}}=\zroot{p}^{\sigma(s)\tau(a)},
\]
where \(\sigma(s)=\pi_1(s)-\pi_1(s_0)\) and \(\tau(a)=\pi_2(a)-\pi_2(a_0)\). Thus \(\sigma\colon S\to\ZN{p}\) and \(\tau\colon A\to\ZN{p}\) are bijections with \(\sigma(s_0)=\tau(a_0)=0\), and
\begin{equation}
\label{eq:fourierform}
(s-s_0)(a-a_0)\equiv\sigma(s)\tau(a)/p\pmod{\Z},\qquad s\in S,\ a\in A.
\end{equation}

Let \(a_1\in A\) satisfy \(\tau(a_1)=1\), and put \(d=a_1-a_0\neq 0\). Then \(d(s-s_0)\equiv\sigma(s)/p\), so \(s-s_0\) is rational. For \(s\neq s_0\) the number \(d(s-s_0)\) has denominator exactly \(p\) in lowest terms, and that denominator divides the denominator of \(s-s_0\). The base point \(s_0\) was arbitrary, so \(p\) divides the denominator of every nonzero element of \(S-S\), and \(p\) divides \(N\).

Let \(S'=N(S-s_0)\) be as in Lemma~\ref{lem:reduction}. For \(\ell=N(s-s_0)\in S'\), relation \eqref{eq:fourierform} reads \(\ell(a-a_0)\equiv K\sigma(s)\tau(a)\pmod{N}\), so \(p\,\ell(a-a_0)\equiv 0\pmod N\) \(\forall \ell\in S'\). Since \(\gcd(N,S')=1\), there are integers \(c_\ell\) with \(\sum_{\ell\in S'}c_\ell\ell\equiv 1\pmod N\), and therefore \(p(a-a_0)\equiv 0\pmod N\). Thus \(a-a_0\bmod N\) lies in \(K\ZN{N}\) \(\forall a\in A\). Reduction modulo \(N\) is injective on \(A\) by Lemma~\ref{lem:reduction}, and \(A\) has \(p\) elements, so the reduction of \(A\) is \(a_0+K\ZN{N}\).

Every \(r\in\ZN{N}\) has exactly one representation \(r=a_0+jK+t\) with \(0\le j<p\) and \(0\le t<K\), so \((a_0+K\ZN{N})\oplus[0,K)=\ZN{N}\). Lemma~\ref{lem:pullback}, applied to reduction modulo \(N\), gives \(A\oplus B=\Z\). The integer \(N\), and with it \(B\), depends only on \(S\).
\end{proof}

The next statement is the structure forced by a single Hadamard matrix of order four whose entries are \(N\)th roots of unity. For \(r\in\ZN{N}\), the condition \(r\bmod N\in[0,M)\) means that the representative in \(\{0,\ldots,N-1\}\) is strictly less than \(M\). Parity in \(\ZN{N}\) is defined only for even \(N\).

\begin{theorem}
\label{thm:znfour}
Let \(N\) be even, \(M=N/2\), and \(S'\subset\ZN{N}\) with \(\lvert S'\rvert=4\), \(0\in S'\), and \(\gcd(N,S')=1\). Let \(A\subset\ZN{N}\) be a four-element spectrum of \(S'\). Then exactly one nonzero element of \(S'\) is even. Call it \(\ell_2\), and call the two odd elements \(\ell_3\) and \(\ell_4\). The set \(A\) is invariant under translation by \(M\), so
\[
A=\{a_0,\,a_0+M,\,a_2,\,a_2+M\}\quad\forall a_0\in A
\]
and a suitable \(a_2\in A\). The cross difference \(v=a_2-a_0\) satisfies
\[
\ell_2 v\equiv M\equiv(\ell_4-\ell_3)v\pmod{N}.
\]
Consequently \(\ell_2(A-a)\subset\{0,M\}\) \(\forall a\in A\), while neither \(\ell_3(A-a)\) nor \(\ell_4(A-a)\) is contained in \(\{0,M\}\). Conversely, if the nonzero elements of \(S'\) consist of one even element \(\ell_2\) and two odd elements \(\ell_3,\ell_4\), and if \(v\in\ZN{N}\) satisfies the displayed congruence, then \(\{a_0,a_0+M,a_0+v,a_0+v+M\}\) is a four-element spectrum of \(S'\) \(\forall a_0\). In particular \(\ell_2\) depends only on \(S'\). Moreover \(4\) divides \(N\).
\end{theorem}

\begin{proof}
Fix \(a_0\in A\) and write \(D=A-a_0=\{0,d_1,d_2,d_3\}\). The matrix \(H=\bigl(\ex{ld/N}\bigr)_{l\in S',\,d\in D}\) is a complex Hadamard matrix: it is the matrix of the spectral pair, multiplied on the left by a diagonal unitary matrix. Its row \(l=0\) and its column \(d=0\) consist of ones. Every other row is orthogonal to the row of ones, so its entries sum to zero and include the entry \(1\). Lemma~\ref{lem:antipodal} gives the multiset of entries as \(\{1,-1,u,-u\}\) for some \(u\) of modulus \(1\). The same holds for every column \(d\neq 0\). Call a row or a column real when every entry is \(\pm 1\). A non-real row or column then contains the entry \(-1\) exactly once.

Let \(B\) be the \(3\times 3\) block of \(H\) on the nonzero rows and columns, and suppose that no nontrivial column of \(H\) is real. Each of the three nontrivial columns of \(H\) begins with the entry \(1\), so each of its entries equal to \(-1\) lies in \(B\). Thus \(B\) contains exactly three entries equal to \(-1\). Each nonzero row of \(H\) also begins with \(1\), so its entries equal to \(-1\) lie in \(B\). Such a row contains the entry \(-1\) at least once, and exactly twice if it is real. If some nonzero row were real, \(B\) would contain at least \(2+1+1=4\) entries equal to \(-1\). Hence every nonzero row is non-real, and each row of \(B\) contains exactly one entry \(-1\). Permuting the nonzero rows, one may place these three entries on the diagonal of \(B\). The remaining two entries of each row are antipodal, and none of them is \(\pm 1\), so
\[
B=\begin{pmatrix}-1&u&-u\\ c&-1&-c\\ d&-d&-1\end{pmatrix}
\]
with \(u,c,d\notin\{\pm 1\}\). The first nontrivial column of \(H\) is \((1,-1,c,d)\). In this multiset the entry \(1\) is antipodal to \(-1\), because neither \(c\) nor \(d\) equals \(-1\), and Lemma~\ref{lem:antipodal} forces the other pair to be antipodal. Hence \(d=-c\). The second nontrivial column is then \((1,u,-1,c)\), and the same reasoning yields \(c=-u\). The second and third rows of \(H\) are therefore \((1,-1,u,-u)\) and \((1,-u,-1,u)\). Their Hermitian inner product is
\[
1+\overline{u}-u-\lvert u\rvert^2=\overline{u}-u,
\]
since \(\lvert u\rvert=1\). Orthogonality forces \(u=\overline{u}\), so \(u=\pm 1\), a contradiction. Some nontrivial column \(d^*\) is therefore real.

For that column, \(\ell d^*\in\{0,M\}\) \(\forall \ell\in S'\), and in particular \(2\ell d^*\equiv 0\pmod{N}\). Since \(\gcd(N,S')=1\), integers \(c_\ell\) exist with \(\sum_{\ell\in S'}c_\ell\ell\equiv 1\pmod{N}\). Multiplying by \(2d^*\) gives \(2d^*\equiv 0\pmod{N}\). The only solutions in \(\ZN{N}\) are \(0\) and \(M\), and \(d^*\neq 0\), so \(d^*=M\). Thus \(M\in D\), the real column is \(\bigl((-1)^\ell\bigr)_{\ell\in S'}\), and orthogonality to the column of ones yields \(\sum_{\ell\in S'}(-1)^\ell=0\). Exactly two elements of \(S'\) are even. One of them is \(0\), so exactly one nonzero element is even. This is the element \(\ell_2\).

The same argument at every base point shows \(a+M\in A\) \(\forall a\in A\). A finite set invariant under translation by \(M\) is a union of orbits of that translation, and here there are two orbits.

It remains to locate \(v\). For \(r\in\ZN{N}\), write \(\widehat A(r)=\sum_{a\in A}\ex{ra/N}\) for the cyclic Fourier sum. With \(A=\{a_0,a_0+M,a_0+v,a_0+v+M\}\),
\begin{equation}
\label{eq:fourfourier}
\widehat A(r)=(1+(-1)^r)\,\ex{ra_0/N}\,(1+\ex{rv/N}),\qquad r\in\ZN{N}.
\end{equation}
The right-hand side vanishes for every odd \(r\). Among the nonzero differences of \(S'=\{0,\ell_2,\ell_3,\ell_4\}\), the even ones are precisely \(\pm\ell_2\) and \(\pm(\ell_4-\ell_3)\), because \(\ell_2\) is even and \(\ell_3,\ell_4\) are odd. For even \(r\), vanishing is equivalent to \(rv\equiv M\pmod{N}\). Both \(\ell_2 v\equiv M\) and \((\ell_4-\ell_3)v\equiv M\) follow. Since \(\ell_2\) is even, \(\ell_2 M\equiv 0\), and therefore \(\ell_2(A-a)\subset\{0,M\}\) \(\forall a\in A\). If \(\ell_3 v\in\{0,M\}\), then \(\ell_4 v=\ell_3 v+M\in\{0,M\}\) as well, so \(vS'\subset\{0,M\}\). The congruence argument of the real column then forces \(v\in\{0,M\}\), contradicting \(\lvert A\rvert=4\). Thus neither \(\ell_3(A-a)\) nor \(\ell_4(A-a)\) lies in \(\{0,M\}\).

For the converse, the congruence \(\ell_2 v\equiv M\) excludes \(v\in\{0,M\}\), so the four elements are distinct. Formula~\eqref{eq:fourfourier} vanishes for odd \(r\), and it vanishes for \(r=\pm\ell_2\) and \(r=\pm(\ell_4-\ell_3)\) because each of those products with \(v\) is congruent to \(M\).

Finally \(\ell_2 v\equiv M\pmod{N}\) with \(\ell_2\) even. The representative \(M=N/2\) is a difference of two even integers, hence even, so \(4\) divides \(N\). The element \(\ell_2\) was read off as the unique even nonzero element of \(S'\), independently of \(A\).
\end{proof}

\begin{theorem}
\label{thm:complement}
Let \(N\), \(M\), and \(S'=\{0,\ell_2,\ell_3,\ell_4\}\) be as in Theorem~\ref{thm:znfour}, with \(\ell_2\) even and \(\ell_3,\ell_4\) odd, and suppose that \(S'\) has a spectrum in \(\ZN{N}\). Define
\[
B=\{g\in\Z:\ \ell_2 g\bmod N\in[0,M)\ \text{and}\ \ell_3 g\bmod N\in[0,M)\}.
\]
Then \(A\oplus B=\Z\) \(\forall A\subset\Z\) with \(\lvert A\rvert=4\) whose reduction modulo \(N\) is a spectrum of \(S'\). The set \(B\) is invariant under translation by \(N\), and \(\lvert B\cap[0,N)\rvert=N/4\). The same conclusions hold for the set defined with \(\ell_4\) in place of \(\ell_3\).
\end{theorem}

\begin{proof}
Let \(\varphi\colon\Z\to\ZN{N}\times\ZN{N}\) be the homomorphism \(\varphi(g)=(\ell_2 g,\ell_3 g)\), and write \(C'=[0,M)\times[0,M)\) inside \(\ZN{N}\times\ZN{N}\). Let \(A\) be as in the statement. By Theorem~\ref{thm:znfour}, the reduction of \(A\) is \(\{a_0,a_0+M,a_0+v,a_0+v+M\}\) with \(\ell_2 v\equiv M\pmod{N}\). Put \(t=\ell_3 v\bmod N\) and \(c_0=\varphi(a_0)\). Evenness of \(\ell_2\) gives \(\ell_2 M\equiv 0\), and oddness of \(\ell_3\) gives \(\ell_3 M\equiv M\), so
\[
\varphi(A)=c_0+\{(0,0),(0,M)\}+\{(0,0),(M,t)\}.
\]
These four points are distinct, so \(\varphi\) is injective on \(A\).

The sum \(\{(0,0),(0,M)\}+C'\) equals \([0,M)\times\ZN{N}\), because \([0,M)\) and \(M+[0,M)\) partition \(\ZN{N}\). Adding \(\{(0,0),(M,t)\}\) then fills \(\ZN{N}\times\ZN{N}\). The two sets have cardinalities \(4\) and \(M^2\), and \(4M^2=N^2\), so the sum is direct: \(\varphi(A)\oplus C'=\ZN{N}\times\ZN{N}\). Lemma~\ref{lem:pullback} yields \(A\oplus B=\Z\).

Translation by \(N\) does not change \(\varphi\), so \(B+N=B\). The same homomorphism on \(\ZN{N}\), restricted to the reduction \(A_0\) of \(A\), satisfies \(A_0\oplus(B\bmod N)=\ZN{N}\) by the pullback lemma in the finite group. Hence \(\lvert B\cap[0,N)\rvert=N/4\). Replacing \(\ell_3\) by \(\ell_4\) replaces \(t\) by \(t+M\), and the same four-point decomposition holds.
\end{proof}

\begin{theorem}
\label{thm:fiberfew}
Let \(h>0\), let \(k\in\{1,2,3,4,5\}\), and let \(E\subset\R\) be a bounded measurable set. Suppose \(\Lambda=\Gamma+h^{-1}\Z\) is a spectrum of \(E\), where \(\lvert\Gamma\rvert=k\) and the elements of \(\Gamma\) are pairwise distinct modulo \(h^{-1}\Z\), and put \(S=h\Gamma\subset\R/\Z\). Let \(B=\Z\) if \(k=1\), let \(B\) be the set of Theorem~\ref{thm:primecomplement} if \(k\in\{2,3,5\}\), and let \(B\) be the complement of Theorem~\ref{thm:complement} built from \(S\) if \(k=4\). Then
\[
\sum_{b\in B}1_E(t-hb)=1
\]
for almost every \(t\in\R\). Thus \(E\) tiles \(\R\) by the translation set \(hB\), which depends on \(E\) only through \(h\) and \(S\).
\end{theorem}

The existence of a common complement in all these cases is due to Dutkay and Haussermann \cite{dutkay-haussermann2015}. The statement above specifies the complement obtained from the constructions in this section.

\begin{proof}
Write \(A_x=\{m\in\Z:x+mh\in E\}\) for \(x\in[0,h)\). By Theorem~\ref{thm:fiber} there is a set \(X_0\subset[0,h)\) of full measure such that \(A_x\) is a spectrum of \(S\) \(\forall x\in X_0\). If \(k=1\), each such fiber is a single point, and \(A_x\oplus\Z=\Z\). If \(k\in\{2,3,5\}\), Theorem~\ref{thm:primecomplement} gives \(A_x\oplus B=\Z\) \(\forall x\in X_0\). If \(k=4\), Lemmas~\ref{lem:fournodulus} and~\ref{lem:reduction} apply to \(S\), and Theorem~\ref{thm:znfour}, applied to any one reduced fiber, shows that \(S'\) has a unique even nonzero element. The complement \(B\) of Theorem~\ref{thm:complement} is then determined by \(S\) and by the choice of base point and of one odd element, and \(A_x\oplus B=\Z\) \(\forall x\in X_0\).

Let \(t=x+hn\) with \(x\in[0,h)\) and \(n\in\Z\). For \(x\in X_0\) and \(b\in B\), the point \(t-hb\) lies in \(E\) iff \(n-b\in A_x\). The sum over \(b\in B\) therefore equals the number of representations of \(n\) in \(A_x\oplus B\), which is \(1\). The exceptional set \(([0,h)\setminus X_0)+h\Z\) is null.
\end{proof}

\begin{remark}
\label{rem:secondcomplement}
The complement of Theorem~\ref{thm:complement} is not the only set that tiles with every spectrum of \(S\). Let \(q_2\) and \(q_4\) be the orders of \(\ell_2/N\) and \((\ell_4-\ell_3)/N\) in \(\R/\Z\), and put \(m=\operatorname{lcm}(q_2,q_4)\). Both orders are even. Indeed \(\ell_2 v\equiv N/2\pmod{N}\) for the cross difference \(v\) of any spectrum, so \(\gcd(\ell_2,N)\) divides \(N/2\) and the reduced denominator \(q_2=N/\gcd(\ell_2,N)\) is even, and likewise for \(q_4\). Define
\[
B(S)=\{b\in\Z:b\bmod N\in[0,N/2)\ \text{and}\ b\bmod m\in[0,m/2)\}.
\]
This set depends on \(S\) only through the two integers \(m\) and \(N\).

The cross difference \(v\) of any spectrum satisfies \((\ell_2/N)v\in\tfrac12+\Z\) and \(((\ell_4-\ell_3)/N)v\in\tfrac12+\Z\). If \(x\in\R/\Z\) has even order \(q\), the integers \(d\) with \(xd\in\tfrac12+\Z\) are exactly \(q/2+q\Z\). In lowest terms \(x=p/q\) with \(p\) odd, and the condition reads \(pd\equiv q/2\pmod{q}\). The inverse of \(p\) modulo \(q\) is odd, so it fixes the class of \(q/2\). An element of \(q/2+q\Z\) has \(2\)-adic valuation \(v_2(q)-1\). The two classes containing \(v\) can meet only when \(v_2(q_2)=v_2(q_4)\). When the valuations are equal to an integer \(a\ge 1\), write \(q_2=2^a u\) and \(q_4=2^a w\) with \(u,w\) odd. Their common elements are \(2^{a-1}\operatorname{lcm}(u,w)\) times an odd integer, which is \(m/2+m\Z\). Thus \(v\equiv m/2\pmod{m}\).

Since \(\ell_2\) is even one has \(v_2(q_2)\le v_2(N)-1\), and \(v_2(m)=v_2(q_2)\), so \(m\) divides \(N/2\).

Number a spectrum so that its residues are \(a_0\), \(a_0+N/2\), \(a_2\), and \(a_2+N/2\), with \(a_2\equiv a_0+m/2\pmod{m}\). The translates by multiples of \(N\) may be absorbed into \(B(S)\), which is \(N\)-periodic, so it is enough to treat these four integers. The pair \(\{a_0,a_2\}\) tiles \(\Z\) with \(B_0=\{b:b\bmod m\in[0,m/2)\}\), because the two residues modulo \(m\) are \(0\) and \(m/2\) after a shift. Since \(m\) divides \(N/2\), the condition modulo \(m\) is unchanged by translation by \(N/2\), and \(B_0\) is the disjoint union of \(B(S)\) and \(B(S)+N/2\). The sum \(\{0,N/2\}\oplus B(S)\) is direct and equals \(B_0\). Composing the two direct sums gives \(A\oplus B(S)=\Z\). There are \((N/2)/m\cdot(m/2)=N/4\) elements of \(B(S)\) in each half-open interval of length \(N\).

The two complements differ in general. For \(N=12\) and \(S'=\{0,2,3,5\}\) one has \(\ell_2=2\), and either odd element may be used in Theorem~\ref{thm:complement}. The resulting set is \(\{0,1,8\}+12\Z\). Here \(q_2=q_4=6\), so \(m=6\) and \(B(S)=\{0,1,2\}+12\Z\). The spectra of \(S'\) are the three translates of \(\{0,3,6,9\}\), and each tiles \(\ZN{12}\) with either complement.
\end{remark}

\section{The six-point fiber}
\label{sec:six}

This section excludes the exceptional fiber of Proposition~\ref{prop:cosetbound}. The input is the rationality of finite spectra with integer frequencies. Once the differences of \(S\) are rational, a consecutive pair in the fiber gives a vanishing sum of six roots of unity. Of the three possible shapes, only the two with rotational symmetry are compatible with the orthogonality relations on the difference set of the fiber.

\begin{theorem}
\label{thm:rational}
Let \((A,S)\) be a finite spectral pair, in the sense of Definition~\ref{def:spectralpair}, with \(A\subset\Z\) and \(S\subset\R/\Z\) of the same cardinality. Then \(S-S\subset\Q/\Z\).
\end{theorem}

The theorem is due to Fu and Song \cite[Theorem~1.2]{fu-song2026}. Their statement takes a spectrum \(\Gamma\) normalized so that \(0\in\Gamma\) and concludes \(\Gamma\subset\Q/\Z\), and the form above is the same fact applied to \(\Gamma=S-s_0\). Their proof rests on a modulus rigidity theorem for Vandermonde systems that satisfy inverse-orthogonality relations \cite[Theorem~1.1]{fu-song2026}, proved with the positive definite kernel \(\sinh(tx)/\sinh(Tx)\), followed by Galois conjugation and Kronecker's theorem. The proof of Theorem~\ref{thm:main} uses Theorem~\ref{thm:rational} only through Lemma~\ref{lem:sixpoint}. Appendix~\ref{sec:rational} gives the proof. It follows Fu and Song, except that their modulus rigidity theorem is obtained as a special case of Lemma~\ref{lem:equalmoduli}, a two-sided form of it which the six-number theorem of Section~\ref{sec:sixnumbers} requires.

\begin{proposition}
\label{prop:sixroots}
A vanishing sum of six distinct roots of unity is a union of three antipodal pairs, or a union of two rotated equilateral triangles, or a rotation of the configuration
\[
\bigl\{\zroot{5},\zroot{5}^2,\zroot{5}^3,\zroot{5}^4,\zroot{6},\zroot{6}^5\bigr\}.
\]
\end{proposition}

\begin{figure}[ht]
\centering
\begin{overpic}[width=0.9\textwidth]{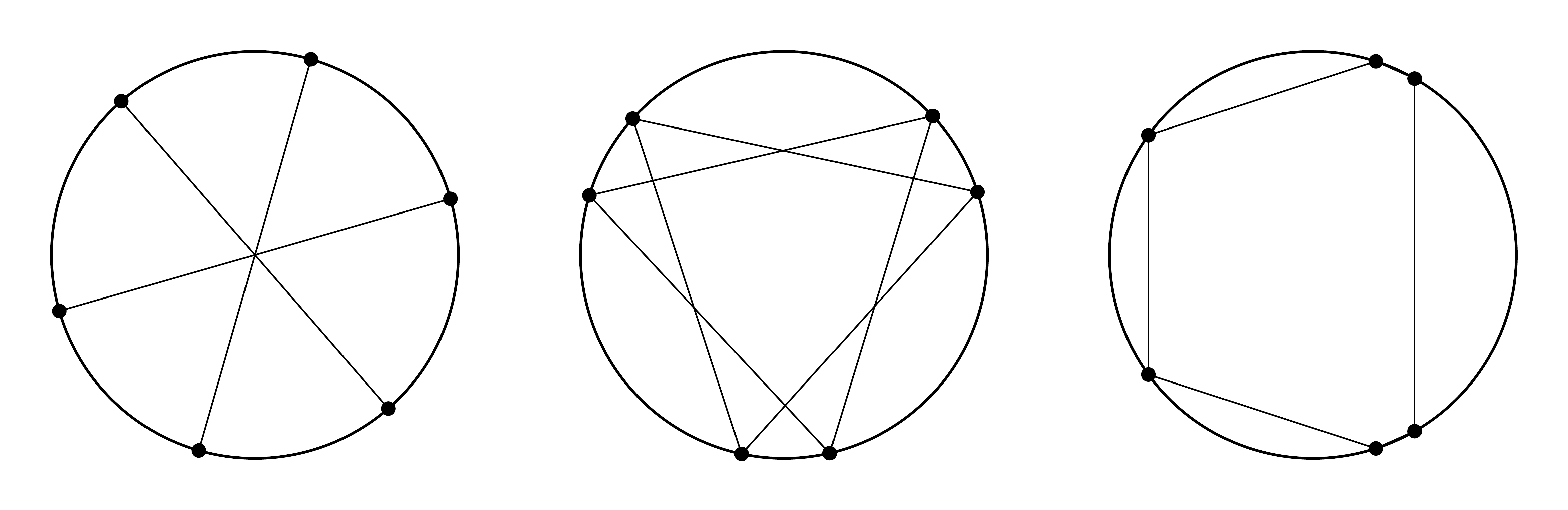}
\end{overpic}
\caption{The three configurations of Proposition~\ref{prop:sixroots}, with each minimal vanishing subsum drawn as the polygon on its roots: three antipodal pairs, two rotated equilateral triangles, and the six roots displayed in the statement, which form a single minimal relation. The first configuration is stable under rotation by a half turn, and the second under rotation by a third of a turn.}
\label{fig:sixconfigs}
\end{figure}

Figure~\ref{fig:sixconfigs} shows the three configurations. The first two are invariant under a nontrivial rotation. The third is excluded for spectral pairs by Lemma~\ref{lem:exotic} and the congruence argument in the proof of Lemma~\ref{lem:sixpoint}.

A relation \(\sum_i a_i\eta_i=0\) with positive integers \(a_i\) and distinct roots of unity \(\eta_i\) has weight \(\sum_i a_i\). It is minimal when no nontrivial subrelation exists. For each prime \(p\), the sum \(R_p\) of all \(p\)th roots of unity is a minimal relation of weight \(p\). The notation \((R_5:R_3)\) is that of Poonen and Rubinstein~\cite{poonen-rubinstein1998}: it is the relation obtained by subtracting a rotated copy of \(R_3\) from \(R_5\) and absorbing the minus signs into the roots, and its underlying set is the configuration displayed above. Mann's theorem, in the form of their Lemma~1, states that a minimal relation of weight \(k\) may be rotated so that every root is an \(m\)th root of unity, where \(m\) is the product of the primes at most \(k\)~\cite{mann1965,poonen-rubinstein1998}. For \(k\le 6\), this product divides \(30\), so every root may be taken to be a \(30\)th root of unity. Their Table~1, whose completeness is their Theorem~3, lists the minimal relations of weight at most six up to rotation: \(R_2\), \(R_3\), \(R_5\), and \((R_5:R_3)\), and no others. In particular there is no minimal relation of weight \(1\) or \(4\). Lam and Leung describe the weights that occur for roots of a fixed order~\cite{lam-leung2000}; the list used here is the table.

\begin{proof}
Any vanishing relation is a sum of minimal relations~\cite{poonen-rubinstein1998}. The six roots may therefore be partitioned into minimal vanishing sums of distinct roots of unity. A single root does not vanish. Two distinct unimodular numbers sum to zero precisely when they are antipodal. Three unimodular numbers \(1\), \(u\) and \(v\) sum to zero only when \(\{1,u,v\}=\{1,\zroot{3},\zroot{3}^2\}\): the identity \(|1+u|=1\) forces \(\operatorname{Re}u=-1/2\). Thus a minimal sum of weight three is a rotated triangle. Four distinct unimodular numbers with sum zero form two antipodal pairs, by Lemma~\ref{lem:antipodal}, so a relation of weight four on distinct roots is never minimal, in agreement with the table. A minimal sum of weight five is a rotated copy of \(R_5\), which cannot appear in a partition of six because the complementary weight would be \(1\). The only remaining partitions of six into parts from \(\{2,3,6\}\) are \(2+2+2\), \(3+3\) and \(6\).

Three antipodal pairs are stable under multiplication by \(-1\). Two rotated triangles are two cosets of the group of cube roots of unity, and that set is stable under multiplication by \(\zroot{3}\). A single minimal sum of weight six is, up to rotation, the configuration in the statement.
\end{proof}

\begin{lemma}
\label{lem:exotic}
Let \(T=\{\zroot{5},\zroot{5}^2,\zroot{5}^3,\zroot{5}^4,\zroot{6},\zroot{6}^5\}\) and let \(\rho\) be a root of unity. For \(t\in\Z\), the power sum \(\sum_{\eta\in\rho T}\eta^t\) vanishes iff \(5\nmid t\) and \(t\equiv\pm 1\pmod{6}\).
\end{lemma}

\begin{proof}
The factor \(\rho^t\) does not affect vanishing. The sum over \(T\) is \(\sum_{a=1}^4\zroot{5}^{at}+2\cos(\pi t/3)\). The sum of all fifth roots of unity vanishes, so the first term equals \(-1\) when \(5\nmid t\) and equals \(4\) when \(5\mid t\). The cosine equals \(1/2\) precisely when \(t\equiv\pm 1\pmod{6}\). The two contributions cancel iff both conditions hold, and \(4+2\cos(\pi t/3)\) is never zero.
\end{proof}

\begin{lemma}
\label{lem:sixpoint}
Let \((A,S)\) be a spectral pair with \(|A|=|S|=6\), \(A\subset\Z\) and \(S\subset\R/\Z\). Suppose \(A\) contains two consecutive integers. Then \(S+\tfrac12=S\) or \(S+\tfrac13=S\).
\end{lemma}

\begin{proof}
Theorem~\ref{thm:rational} gives \(S-S\subset\Q/\Z\). Translate \(S\) so that \(0\in S\). The six numbers \(\omega_s=e^{2\pi i s}\) are then distinct roots of unity. Two consecutive integers of \(A\) differ by \(1\), so \(\sum_{s\in S}\omega_s=\widehat{\muS}(1)=0\). Proposition~\ref{prop:sixroots} says that these roots form three antipodal pairs, two rotated triangles, or a rotation of the configuration in Lemma~\ref{lem:exotic}.

Three antipodal pairs are stable under \(\omega\mapsto-\omega\), so \(S+\tfrac12=S\). Two rotated triangles are stable under multiplication by \(\zroot{3}\), so \(S+\tfrac13=S\). In the remaining case, Lemma~\ref{lem:exotic} says that \(\widehat{\muS}(t)=0\) iff \(5\nmid t\) and \(t\equiv\pm 1\pmod{6}\). The pair being spectral, \(\widehat{\muS}\) vanishes on \((A-A)\setminus\{0\}\), by Definition~\ref{def:spectralpair}. Every nonzero difference of elements of \(A\) is therefore congruent to \(\pm 1\) modulo \(6\).

This is impossible as soon as \(A\) has three elements. If \(a>b>c\) lie in \(A\), both \(a-b\) and \(b-c\) are congruent to \(\pm 1\) modulo \(6\), so their sum \(a-c\) is congruent to \(0\) or to \(\pm 2\) modulo \(6\). None of these three residue classes is \(\pm 1\), so \(\widehat{\muS}(a-c)\neq 0\). The exotic configuration does not occur.
\end{proof}

The argument did not use the size of the gaps in \(A\), only one consecutive pair and the existence of three elements. No further restriction on the run structure of \(A\) is needed.

\begin{corollary}
\label{cor:k4}
Let \(\Omega\subset\R\) be a union of four pairwise disjoint non-adjacent bounded intervals, and let \(\Lambda\) be a spectrum of \(\Omega\), with minimal period \(h^{-1}\). Then \(k=\lvert\Omega\rvert/h\le 4\).
\end{corollary}

\begin{proof}
By Proposition~\ref{prop:cosetbound}, either \(k\le 4\), or else, up to translation and scaling, \(\Omega\) is a union of integer runs of lengths \(2\), \(2\), \(1\) and \(1\), with \(k=6\) and with \(S\) free of translation symmetry. In the second case Theorem~\ref{thm:fiber} says that the fiber \(A_x\) is a spectrum of \(S\) for almost every \(x\). That fiber is a six-element subset of \(\Z\) containing two consecutive integers. Lemma~\ref{lem:sixpoint} gives \(S+\tfrac12=S\) or \(S+\tfrac13=S\), contradicting the absence of a translation symmetry. The second case is empty, and \(k\le 4\).
\end{proof}

\section{Proof of the main theorem}
\label{sec:assembly}

\begin{figure}[ht]
\centering
\begin{overpic}[width=\textwidth]{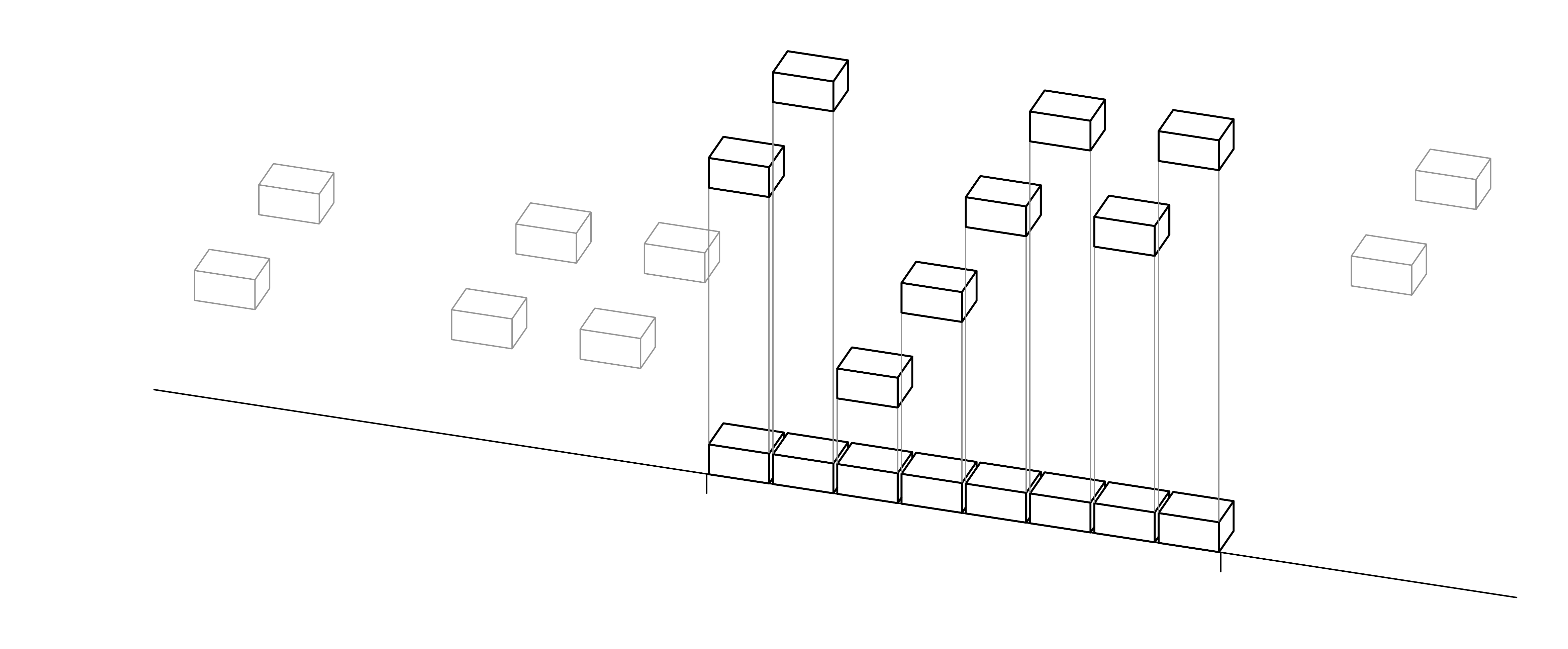}
\put(45.08,8.81){\makebox(0,0){\small$0$}}
\put(77.87,3.80){\makebox(0,0){\small$8$}}
\put(11.18,23.18){\makebox(0,0)[r]{$\Omega-8$}}
\put(15.28,28.65){\makebox(0,0)[r]{$\Omega-7$}}
\put(43.97,30.36){\makebox(0,0)[r]{$\Omega$}}
\put(48.07,35.82){\makebox(0,0)[r]{$\Omega+1$}}
\end{overpic}
\caption{The set \(\Omega=[0,1)\cup[4,5)\cup[6,7)\cup[10,11)\) of Theorem~\ref{thm:weak}, which is spectral with \(h=1\) and \(S=\{0,\tfrac18,\tfrac28,\tfrac38\}\), and its translates by \(-8\), \(-7\) and \(1\), raised to separate levels. For this \(S\) the complement is \(B=\{0,1\}+8\Z\). The pieces lying over \([0,8)\), between the marks \(0\) and \(8\) on the base line, fill one period of the line, and the parts outside it are gray.}
\label{fig:tiling}
\end{figure}

\begin{proof}[Proof of Theorem~\ref{thm:main}]
Let \(\Omega\subset\R\) be a bounded measurable set which agrees almost everywhere with a union of four bounded intervals, and suppose that \(\Omega\) is spectral. Null sets do not affect spectrality or tiling. After merging overlapping or adjacent intervals and discarding intervals of length zero, we may therefore take \(\Omega\) to be a disjoint union of at most four non-adjacent half-open intervals. If there are at most three components, the set tiles \(\R\) by translations, by Theorem~1.1 of Shi \cite{shi2026}. Assume from here on that
\[
\Omega=\bigsqcup_{j=1}^{4}[a_j,b_j),
\]
with four pairwise disjoint non-adjacent components. Let \(\Lambda\) be a spectrum, and write \(L=\lvert\Omega\rvert\).

By Theorem~\ref{thm:periodic}, \(\Lambda\) is periodic, every positive period lies in \(L^{-1}\Z_{>0}\), and there is a minimal positive period \(h^{-1}\) with \(h>0\). Write \(\Lambda=\Gamma+h^{-1}\Z\), where the elements of the finite set \(\Gamma\) are pairwise distinct modulo \(h^{-1}\Z\), and put \(k=\lvert\Gamma\rvert\) and \(S=h\Gamma\subset\R/\Z\). The set \(S\) has no nontrivial translation symmetry. Theorem~\ref{thm:fiber} gives \(L=kh\), so \(\Omega\) is a \(k\)-fold multitile by \(h\Z\), and for almost every \(x\in[0,h)\) the fiber \(A_x\) is a spectrum of \(S\).

Proposition~\ref{prop:cosetbound} gives \(k\le 4\), except possibly when, up to translation and scaling, \(\Omega\) is a union of integer runs of lengths \(2,2,1,1\), with \(k=6\) and with \(S\) free of translation symmetry. Lemma~\ref{lem:sixpoint} excludes that six-point fiber, and Corollary~\ref{cor:k4} records the resulting bound: \(k\le 4\).

Theorem~\ref{thm:fiberfew}, applied to \(E=\Omega\), now gives a tiling of \(\R\) by \(hB\). If \(k=1\), then \(h=L\), \(B=\Z\), and the tiling is by \(L\Z\). If \(k\in\{2,3\}\), the set \(B\) is that of Theorem~\ref{thm:primecomplement}. If \(k=4\), it is the complement of Theorem~\ref{thm:complement}.

The translation set in the case \(k=4\) is read off from \(S\). Lemmas~\ref{lem:reduction} and~\ref{lem:fournodulus} produce the even modulus \(N\) and the four-element set \(S'\subset\ZN{N}\). Theorem~\ref{thm:znfour} selects the unique even nonzero element \(\ell_2\in S'\) and either odd element \(\ell_3\). The complement of Theorem~\ref{thm:complement} is
\[
B=\{g\in\Z:\ell_2 g\bmod N\in[0,N/2)\ \text{and}\ \ell_3 g\bmod N\in[0,N/2)\}.
\]
It is invariant under translation by \(N\), it meets each period in \(N/4\) integers, and it depends on \(\Omega\) only through \(S\). Theorem~\ref{thm:fiberfew} gives the tiling identity
\[
\sum_{b\in B}1_\Omega(t-hb)=1
\]
for almost every \(t\in\R\). Figure~\ref{fig:tiling} shows one period of such a tiling.
\end{proof}

\section{Six unimodular numbers with four vanishing power sums}
\label{sec:sixnumbers}

The exceptional fiber considered in the proof of Corollary~\ref{cor:k4} has the form \(A=\{0,1,x,x+1,y,z\}\). If \(A\) is a spectrum of a six-point set \(S\), then \(\widehat{\muS}(n)=\sum_{s\in S}e^{2\pi i ns}=0\) at every nonzero difference \(n\in A-A\), in particular at \(n=1\), \(x-1\), \(x\) and \(x+1\). Lemma~\ref{lem:sixpoint} excluded such an \(S\) without translation symmetry by first proving that \(S-S\) is rational, which used the full orthogonality of the spectral pair. This section proves that the four power sums alone force the symmetry, with no rationality assumed and no other difference of \(A\) used. The statement concerns six complex numbers and no spectral pair, and it may be of independent interest. For algebraic solutions, the proof forces equality of the moduli at every place of the number field generated by the six numbers. Kronecker's theorem and the classification of vanishing sums of six roots of unity then give the symmetry, and the Nullstellensatz extends the conclusion to arbitrary complex solutions.

For nonzero complex numbers \(\omega_1,\ldots,\omega_6\) we write \(\ps{n}=\sum_{s=1}^6\omega_s^n\) for \(n\in\Z\) and \(\es{j}\) for the elementary symmetric functions. Negative exponents are allowed because the \(\omega_s\) are nonzero.

The half-turn and third-turn symmetries of six unimodular numbers with sum zero can be read from the coefficients of their polynomial.

\begin{lemma}
\label{lem:normalform}
Let \(\omega_s=e^{2\pi i s}\) for \(s\) in a six-element set \(S\subset\R/\Z\), and suppose \(\sum_{s\in S}\omega_s=0\). After a rotation of \(S\), the polynomial \(q(z)=\prod_{s\in S}(z-\omega_s)\) has the form
\[
q(z)=z^6+cz^4-rz^3+\overline{c}\,z^2+1
\]
with \(c\in\C\) and \(r\in\R\). The power sums \(\ps{n}=\sum_s\omega_s^n\) and the elementary symmetric functions \(\es{j}\) then satisfy \(\ps{2}=-2c\), \(\ps{3}=3r\) and \(\ps{5}=-5cr\), and \(S+\tfrac12=S\) or \(S+\tfrac13=S\) iff \(cr=0\).
\end{lemma}

\begin{proof}
The roots are unimodular, so \(\overline{\omega_s}=\omega_s^{-1}\). Rotating \(S\) multiplies every root by a fixed unimodular constant. The sum remains zero, and the constant may be chosen so that the product of the roots equals \(1\). For a monic polynomial with unimodular roots and constant term \(1\), the elementary symmetric functions satisfy \(\overline{\es{j}}=\es{6-j}\). The vanishing of the sum is \(\es{1}=0\), so \(\es{5}=0\), and \(\es{3}\) is real. Writing \(c=\es{2}\) and \(r=\es{3}\) gives the displayed shape, with \(\es{4}=\overline{c}\).

Newton's identities~\cite[Chapter~I, \S 2]{macdonald1995} give the three power sums. For \(1\le k\le 6\),
\[
\ps{k}-\es{1}\ps{k-1}+\es{2}\ps{k-2}-\cdots+(-1)^{k-1}\es{k-1}\ps{1}+(-1)^k k\es{k}=0,
\]
with \(\ps{0}\) omitted. The cases \(k=2\), \(k=3\) and \(k=5\), using \(\es{1}=\es{5}=0\), \(\es{2}=c\), \(\es{3}=r\) and \(\es{4}=\overline{c}\), are \(\ps{2}+2c=0\), \(\ps{3}-3r=0\) and \(\ps{5}+3cr+2rc=0\).

If \(r=0\), then \(q\) is a polynomial in \(z^2\), so its roots are stable under \(\omega\mapsto-\omega\) and \(S+\tfrac12=S\). If \(c=0\), then \(q\) is a polynomial in \(z^3\), so its roots are stable under multiplication by \(\zroot{3}\) and \(S+\tfrac13=S\). Conversely, \(S+\tfrac12=S\) means that the set of roots is stable under \(\omega\mapsto-\omega\). Then \(q(-z)=q(z)\), because both are monic of degree six with the same roots, and the coefficient of \(z^3\) forces \(r=0\). If \(S+\tfrac13=S\), the same comparison of \(q(\zroot{3}z)\) with \(q(z)\) forces the coefficients of \(z^4\) and \(z^2\) to vanish.
\end{proof}

\begin{theorem}
\label{thm:sixnumbers}
Let \(x\ge 2\) be an integer and let \(\omega_1,\ldots,\omega_6\) be nonzero complex numbers with \(\omega_1\cdots\omega_6=1\) whose power sums \(\ps{n}\) vanish for every exponent \(n\) in the set
\[
\{\pm 1,\ \pm(x-1),\ \pm x,\ \pm(x+1)\}.
\]
Then \(\es{3}=0\), or \(\es{2}=\es{4}=0\). In particular, six unimodular numbers with sum zero and \(\ps{x-1}=\ps{x}=\ps{x+1}=0\) satisfy \(\es{2}\es{3}=0\).
\end{theorem}

The second statement follows from the first. If \(|\omega_s|=1\) \(\forall s\), then \(\overline{\omega_s}=\omega_s^{-1}\), so \(\ps{-n}=\overline{\ps{n}}\), and the vanishing of \(\ps{1}\), \(\ps{x-1}\), \(\ps{x}\), \(\ps{x+1}\) gives the vanishing at the indicated positive and negative exponents. Replacing every \(\omega_s\) by \(\lambda\omega_s\), with \(\lambda\) unimodular and \(\lambda^6\,\omega_1\cdots\omega_6=1\), multiplies \(\ps{n}\) by \(\lambda^n\) and \(\es{j}\) by \(\lambda^j\), so it preserves the hypotheses and multiplies \(\es{2}\es{3}\) by \(\lambda^5\). The first statement then gives \(\es{3}=0\) or \(\es{2}=0\). Distinctness of the \(\omega_s\) is not needed anywhere in this section. The theorem is stated for \(x\ge 2\), which is what the proof needs, namely that \(0,1,x,x+1\) are four distinct integers and that \(x-1\ge 1\). The exceptional fiber considered in the proof of Corollary~\ref{cor:k4} has \(x\ge 3\), its two runs of length two being separated by at least one integer.

Throughout the section, \(W=\{x-1,x,x+1\}\) and \(-W=\{-x-1,-x,-x+1\}\). Three consecutive integers contain an even one and an odd one; we write \(n_e\) for the smallest even element of \(W\) and \(n_o\) for the largest odd element of \(W\), so that \(n_o\ge 3\) because \(x+1\ge 3\). We say that nonzero elements \(\omega_1,\ldots,\omega_6\) of a field satisfy the window condition at \(x\) when \(\omega_1\cdots\omega_6=1\) and \(\ps{n}=0\) \(\forall n\in\{1,-1\}\cup W\cup(-W)\). Replacing every \(\omega_s\) by \(\omega_s^{-1}\) preserves the window condition, since it exchanges \(\ps{n}\) with \(\ps{-n}\) and the exponent set is symmetric. The window condition fixes three of the elementary symmetric functions: \(\es{1}=\ps{1}=0\), \(\es{6}=1\), and \(\es{5}=\es{6}\sum_s\omega_s^{-1}=\ps{-1}=0\), so
\begin{equation}
\label{eq:windowshape}
\prod_{s=1}^6(z-\omega_s)=z^6+\es{2}z^4-\es{3}z^3+\es{4}z^2+1 .
\end{equation}

The archimedean step rests on the positive definite kernel of Fu and Song \cite[Lemma~4.1]{fu-song2026}.

\begin{lemma}
\label{lem:kernel}
Let \(0<t<T\). Define \(k(x)=\sinh(tx)/\sinh(Tx)\) for \(x\neq 0\), and \(k(0)=t/T\). Then \(k\) is positive definite on \(\R\): for every finite sequence of real numbers \(x_1,\ldots,x_n\), the matrix \(\bigl(k(x_p-x_q)\bigr)\) is positive semidefinite.
\end{lemma}

\begin{proof}
The value at the origin is the limit of the quotient, and \(k\) is even. Euler's product for the sine, in the form \(\sinh z=z\prod_{n\ge 1}\bigl(1+z^2/(n^2\pi^2)\bigr)\), gives
\[
k(x)=\frac{t}{T}\prod_{n\ge 1}f_n(x),
\qquad
f_n(x)=\frac{n^2\pi^2+t^2 x^2}{n^2\pi^2+T^2 x^2},
\]
at \(x=0\) as well as for \(x\neq 0\), since every factor \(f_n(0)\) equals \(1\). Each factor splits as a constant term plus a Cauchy kernel,
\[
f_n(x)=\Bigl(\frac{t}{T}\Bigr)^2+\Bigl(1-\Bigl(\frac{t}{T}\Bigr)^2\Bigr)\frac{1}{1+T^2 x^2/(n^2\pi^2)}.
\]
Both coefficients are positive.

The Gaussian kernel is positive definite. For \(g>0\) and a finite set of real points, the matrix \(\bigl(e^{2gab}\bigr)\) is an entrywise limit of sums \(\sum_{k\le K}(2g)^k v_k v_k^\top/k!\) with \(v_k=(a^k)_a\), each of which is positive semidefinite, and a limit of positive semidefinite matrices is positive semidefinite. The congruence \(M\mapsto DMD\) by the positive diagonal matrix \(D=\operatorname{diag}(e^{-ga^2})\) preserves that property and produces \(\bigl(e^{-g(a-b)^2}\bigr)\). The Cauchy kernel is a positive integral of Gaussians,
\[
\frac{1}{1+c^2 x^2}=\int_0^\infty e^{-s}e^{-s c^2 x^2}\,ds\qquad(c>0),
\]
so its matrices are positive semidefinite as well. Each factor matrix \((f_n(a-b))\) is a positive combination of the all-ones matrix and a Cauchy matrix, hence positive semidefinite.

The Schur product theorem~\cite[Theorem~7.5.3]{horn-johnson2013} closes the cone of positive semidefinite matrices under entrywise products. Every partial product \((t/T)\,f_1\cdots f_N\), evaluated on a finite set of differences, is therefore positive semidefinite. These partial products converge entrywise to \(k\), by the Euler product, and the limit is positive semidefinite.

The Fourier transform gives the same conclusion directly. One has
\[
k(x)=\int_\R\rho(\xi)\,\ex{x\xi}\,d\xi,\qquad \rho(\xi)=\frac{\pi\sin(\pi t/T)}{T\bigl(\cosh(2\pi^2\xi/T)+\cos(\pi t/T)\bigr)},
\]
where \(\rho\) is positive and integrable, and a positive combination of characters is positive definite.
\end{proof}

The next lemma is a two-sided form of the modulus rigidity theorem of Fu and Song \cite[Theorem~1.1]{fu-song2026}. The two eliminations in its proof give weighted counts of the nodes, each at least the number of exponents; since the weights assigned to each modulus class sum to one, the two counts add to the total number of nodes. The comparison with their theorem follows the proof.

\begin{lemma}
\label{lem:equalmoduli}
Let \(A'\subset\Z\) be a set of \(N\) integers and let \(D=(A'-A')\setminus\{0\}\). Let \(w_1,\ldots,w_d\) be nonzero complex numbers, not necessarily distinct, with \(d<2N\), and suppose \(\sum_{j=1}^d w_j^n=0\) \(\forall n\in D\). Then \(|w_1|=\cdots=|w_d|\).
\end{lemma}

\begin{figure}[ht]
\centering
\begin{overpic}[width=0.62\textwidth]{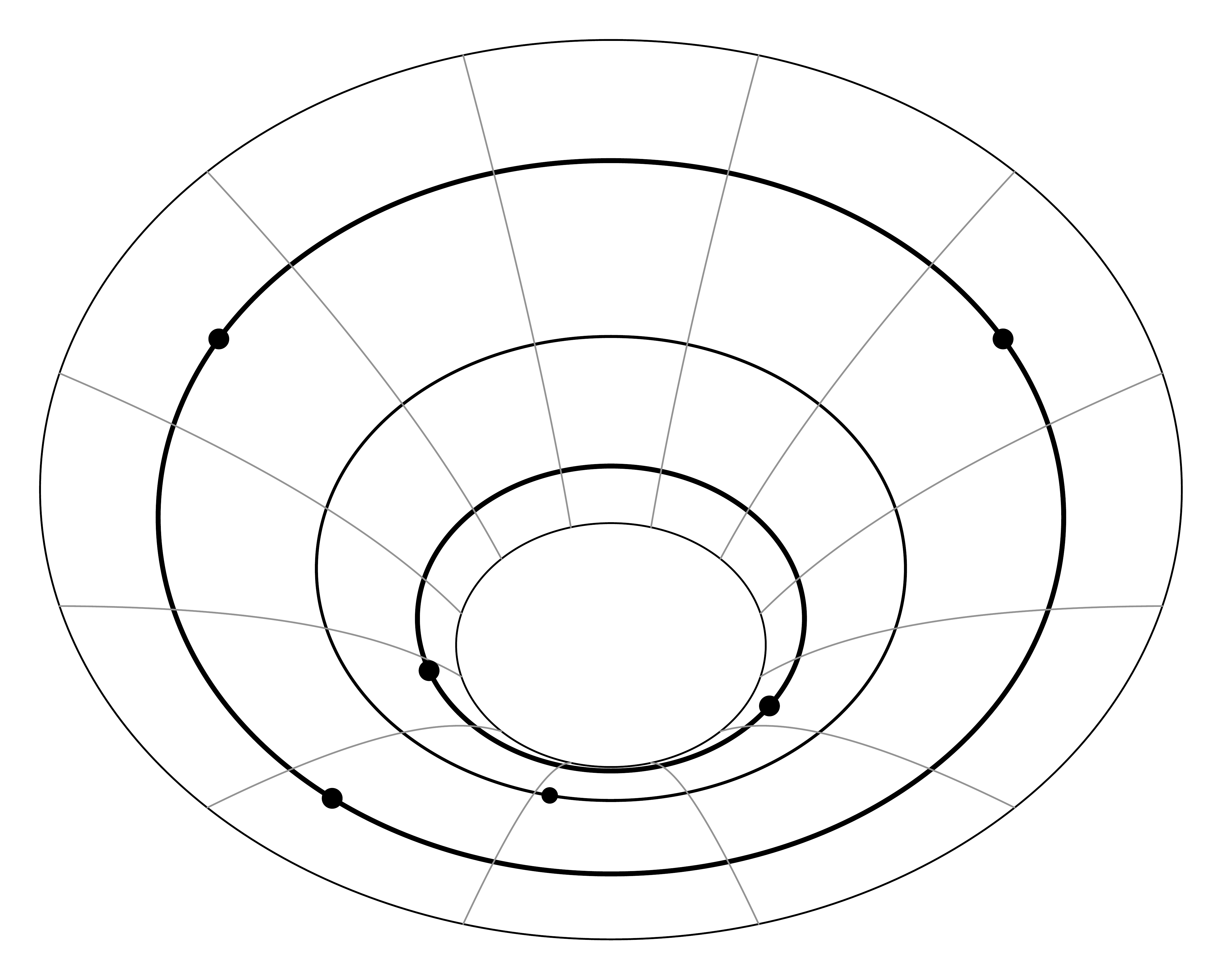}
\put(69.97,31.86){\makebox(0,0){\small$e^{t_1}$}}
\put(80.58,36.01){\makebox(0,0){\small$e^{t_2}$}}
\put(91.89,39.08){\makebox(0,0){\small$e^{T}$}}
\end{overpic}
\caption{The modulus classes of Lemma~\ref{lem:equalmoduli} on the graph of \(\log\lvert w\rvert\) over the plane. The class of modulus \(e^{t_r}\) is the circle of that radius at height \(t_r\). The innermost and the outermost modulus classes, at heights \(0\) and \(T\), are drawn heavy: each of the two eliminations~\eqref{eq:twoelim} expresses the power sums of one of them through those of the intermediate classes.}
\label{fig:moduli}
\end{figure}

\begin{proof}
Multiplying every \(w_j\) by a positive constant \(c\) multiplies \(\sum_j w_j^n\) by \(c^n\) and preserves the hypothesis, so we may assume that the least modulus is \(1\). Write the distinct moduli as \(e^{t_1}<\cdots<e^{t_s}\) with \(0=t_1<\cdots<t_s=T\), as in Figure~\ref{fig:moduli}. Let \(I_r\) be the class of indices of modulus \(e^{t_r}\), of size \(m_r\ge 1\), so that \(\sum_r m_r=d\). Write \(w_j=e^{t_r}u_j\) with \(|u_j|=1\) for \(j\in I_r\), and put \(U_r(n)=\sum_{j\in I_r}u_j^n\), so that \(U_r(-n)=\overline{U_r(n)}\). Suppose \(s\ge 2\); the aim is to show \(d\ge 2N\).

Fix \(n\in D\). The hypothesis at \(n\), grouped by classes, is the first of the two relations
\[
\sum_{r=1}^s e^{t_rn}\,U_r(n)=0,
\qquad
\sum_{r=1}^s e^{-t_rn}\,U_r(n)=0 .
\]
Since \(D=-D\), the hypothesis also holds at \(-n\), where it reads
\[
\sum_{r=1}^s e^{-t_rn}\,U_r(-n)=0 .
\]
Its complex conjugate is the second relation, because \(e^{-t_rn}\) is real and \(\overline{U_r(-n)}=U_r(n)\). Subtracting the second relation from the first, and subtracting \(e^{Tn}\) times the second from \(e^{-Tn}\) times the first, gives
\[
\sum_{r=1}^s 2\sinh(t_rn)\,U_r(n)=0,
\qquad
\sum_{r=1}^s 2\sinh((T-t_r)n)\,U_r(n)=0,
\]
where the second relation uses \(e^{(t_r-T)n}-e^{(T-t_r)n}=-2\sinh((T-t_r)n)\). In the first sum the term \(r=1\) vanishes because \(t_1=0\), and in the second the term \(r=s\) vanishes because \(t_s=T\). Since \(T>0\) and \(n\ne 0\), \(\sinh(Tn)\ne 0\), and dividing gives the two eliminations
\begin{equation}
\label{eq:twoelim}
\begin{gathered}
U_s(n)=-\sum_{r=2}^{s-1}\frac{\sinh(t_rn)}{\sinh(Tn)}\,U_r(n),\\
U_1(n)=-\sum_{r=2}^{s-1}\frac{\sinh((T-t_r)n)}{\sinh(Tn)}\,U_r(n),
\end{gathered}
\end{equation}
valid \(\forall n\in D\), with empty sums when \(s=2\). The innermost class has dropped out of the first, and the outermost out of the second.

For each class let \(E_r\) be the \(N\times m_r\) matrix with columns \(c_j=(u_j^a)_{a\in A'}\), \(j\in I_r\), and let \(G_r=E_rE_r^*=\sum_{j\in I_r}c_jc_j^*\). This is a positive semidefinite Hermitian matrix indexed by \(A'\), with diagonal entries \(m_r\) and with off-diagonal entries \((G_r)_{ab}=\sum_{j\in I_r}u_j^a\overline{u_j^b}=U_r(a-b)\), because \(\overline{u_j}=u_j^{-1}\). For \(2\le r\le s-1\) let \(C_r\) and \(C_r'\) be the matrices indexed by \(A'\) with entries \(k_{t_r}(a-b)\) and \(k_{T-t_r}(a-b)\), where \(k_t(x)=\sinh(tx)/\sinh(Tx)\) for \(x\ne 0\) and \(k_t(0)=t/T\). Since \(0<t_r<T\) and \(0<T-t_r<T\), Lemma~\ref{lem:kernel} says that \(C_r\) and \(C_r'\) are positive semidefinite. Their diagonal entries are \(t_r/T\) and \(1-t_r/T\).

Consider the matrix \(L=G_s+\sum_{r=2}^{s-1}C_r\circ G_r\), where \(\circ\) is the entrywise product. For \(a\ne b\) in \(A'\), the entry \(L_{ab}\) is \(U_s(n)+\sum_{r=2}^{s-1}k_{t_r}(n)U_r(n)\) with \(n=a-b\in D\), which is zero by the first elimination in~\eqref{eq:twoelim}. Every diagonal entry equals \(\lambda=m_s+\sum_{r=2}^{s-1}(t_r/T)\,m_r\). Hence \(L=\lambda I\), and
\[
\lambda I-G_s=\sum_{r=2}^{s-1}C_r\circ G_r .
\]
Each summand is positive semidefinite. Indeed \(C_r\circ G_r=\sum_{j\in I_r}C_r\circ(c_jc_j^*)\), and for any vector \(c\) the matrix \(C_r\circ(cc^*)\) has entries \((C_r)_{ab}\,c_a\overline{c_b}\), which are the entries of \(D_cC_rD_c^*\) with \(D_c\) the diagonal matrix with diagonal \(c\); a congruence preserves positive semidefiniteness. Therefore \(\lambda I-G_s\) is positive semidefinite. Now choose \(j_0\in I_s\) and test with \(v=c_{j_0}\), whose squared length is \(v^*v=\sum_{a\in A'}|u_{j_0}^a|^2=N\). On one side \(v^*(\lambda I)v=\lambda N\). On the other, \(v^*G_sv=v^*E_sE_s^*v=\sum_{j\in I_s}|c_j^*v|^2\ge|c_{j_0}^*c_{j_0}|^2=N^2\). Positive semidefiniteness of \(\lambda I-G_s\) gives \(\lambda N\ge N^2\), that is \(\lambda\ge N\).

The same argument with \(L'=G_1+\sum_{r=2}^{s-1}C_r'\circ G_r\) uses the second elimination in~\eqref{eq:twoelim}: the off-diagonal entries of \(L'\) vanish, the diagonal entries equal \(\lambda'=m_1+\sum_{r=2}^{s-1}(1-t_r/T)\,m_r\), so \(\lambda'I-G_1\) is positive semidefinite, and the test vector \(c_{j_1}\) with \(j_1\in I_1\) gives \(\lambda'\ge N\).

Adding the two inequalities,
\[
2N\le\lambda+\lambda'=m_1+m_s+\sum_{r=2}^{s-1}\Bigl(\frac{t_r}{T}+1-\frac{t_r}{T}\Bigr)m_r=\sum_{r=1}^sm_r=d,
\]
which contradicts \(d<2N\). Hence \(s=1\).
\end{proof}

\begin{figure}[ht]
\centering
\begin{overpic}[width=0.5\textwidth]{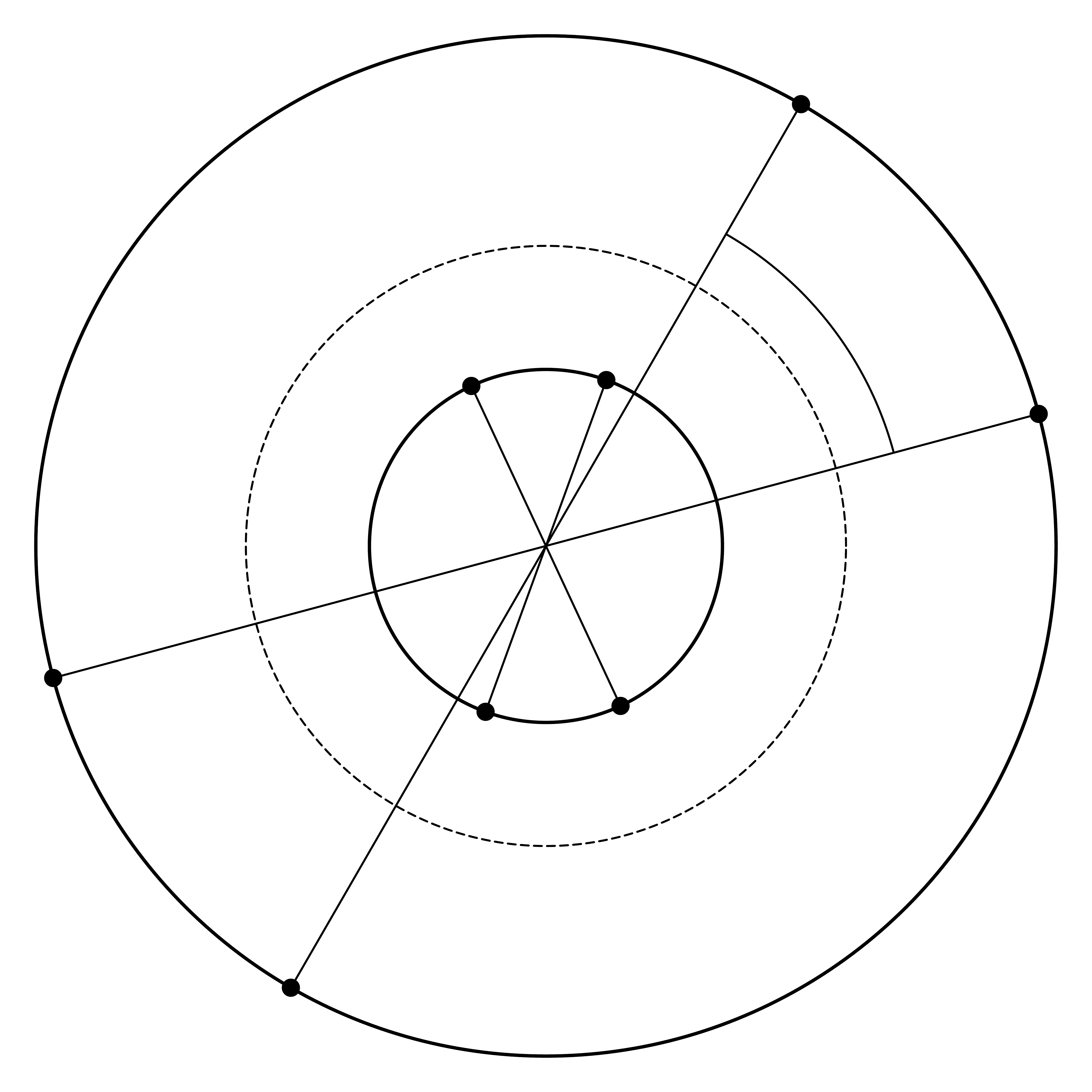}
\put(81.83,74.43){\makebox(0,0){\small$\pi/x$}}
\put(88.76,19.29){\makebox(0,0)[l]{\small$R$}}
\put(49.28,70.60){\makebox(0,0){\small$R^{-1}$}}
\end{overpic}
\caption{The eight numbers that show the bound \(d<2N\) of Lemma~\ref{lem:equalmoduli} to be sharp, for \(A'=\{0,1,x,x+1\}\) and \(x=4\). Each modulus class, on the circle of radius \(R\) or \(R^{-1}\), consists of two antipodal pairs, drawn as diameters and rotated against each other by \(\pi/x\). The dashed circle is the unit circle. For six numbers, fewer than \(2N=8\), Lemma~\ref{lem:equalmoduli} leaves a single circle, and Corollary~\ref{cor:unimodular} identifies it as the unit circle.}
\label{fig:twomoduli}
\end{figure}

For \(N=d\), Lemma~\ref{lem:equalmoduli} contains the modulus rigidity theorem of Fu and Song: by Lemma~\ref{lem:columns}, the inverse-orthogonality relations of a Vandermonde system with \(d\) nodes and \(d\) exponents give the vanishing of the power sums on the difference set of the exponents, and \(d<2d\). Appendix~\ref{sec:rational} uses it in that form. Their proof has, in addition, the row relations inside each class at its disposal, and it uses only the first of the two eliminations~\eqref{eq:twoelim}. Here only the power sums on the difference set are available, and the number of nodes may be as large as \(2N-1\). The loss is compensated by using both eliminations and adding the two bounds, and in the last step a congruence by a diagonal matrix takes the place of the Schur product theorem. The bound \(d<2N\) is sharp. For \(A'=\{0,1,x,x+1\}\) with \(x\) even, \(N=4\), and \(d=8\), the numbers \(\pm Ru\), \(\pm Rue^{i\pi/x}\), \(\pm R^{-1}u'\), \(\pm R^{-1}u'e^{i\pi/x}\), with \(R>1\) and \(u,u'\) unimodular, have two moduli, and their power sums vanish on \(D=\pm\{1,x-1,x,x+1\}\): each class is two antipodal pairs, so its power sum vanishes at odd \(n\), and at \(n=\pm x\) the factor \(1+e^{\pm i\pi}\) vanishes. Figure~\ref{fig:twomoduli} shows this configuration for \(x=4\).

\begin{corollary}
\label{cor:unimodular}
Let \(x\ge 2\). Nonzero complex numbers \(\omega_1,\ldots,\omega_6\) satisfying the window condition at \(x\) are unimodular.
\end{corollary}

\begin{proof}
Take \(A'=\{0,1,x,x+1\}\), a set of \(N=4\) integers because \(x\ge 2\). Its nonzero differences are \(\pm 1\), \(\pm(x-1)\), \(\pm x\), \(\pm(x+1)\), the exponents of the window condition. Since \(d=6<8=2N\), Lemma~\ref{lem:equalmoduli} gives a common modulus \(\rho\). Then \(1=|\omega_1\cdots\omega_6|=\rho^6\), so \(\rho=1\).
\end{proof}

The second step treats the non-archimedean places, and for it \(K\) is a field of characteristic \(0\) with a non-archimedean absolute value \(|\cdot|\), with valuation ring \(\mathcal O=\{a:|a|\le 1\}\), maximal ideal \(\mathfrak m=\{a:|a|<1\}\), residue field \(\kappa=\mathcal O/\mathfrak m\), and residue characteristic \(p\), which is \(0\) or a prime. A unit is an element of absolute value \(1\), and \(\bar a\in\kappa\) denotes the residue of \(a\in\mathcal O\). The ultrametric inequality \(|a+b|\le\max(|a|,|b|)\), with equality when \(|a|\ne|b|\), has the consequence used throughout: in a vanishing finite sum, the largest absolute value among the terms is attained by at least two terms, since a term of strictly largest absolute value would make the sum nonzero. Every integer has absolute value at most \(1\), and every nonzero integer has positive absolute value, because \(K\) has characteristic \(0\). The example to keep in mind is a number field \(K\) with the \(\mathfrak p\)-adic absolute value attached to a nonzero prime ideal \(\mathfrak p\) of its ring of integers \(\mathcal O_K\): the valuation ring is the localization of \(\mathcal O_K\) at \(\mathfrak p\), the residue field is \(\mathcal O_K/\mathfrak p\), and \(p\) is the rational prime below \(\mathfrak p\).

\begin{lemma}
\label{lem:supnorm}
Let \(u_1,\ldots,u_m\) be units of \(\mathcal O\), not necessarily distinct, and let \(\sigma_n=\sum_{i=1}^mu_i^n\) for \(n\in\Z\). Then \(\forall j\in\Z\),
\[
\max\bigl(|\sigma_j|,|\sigma_{j+1}|,\ldots,|\sigma_{j+m-1}|\bigr)=\sup_{n\in\Z}|\sigma_n| ,
\]
and this common value is at least \(|\sigma_0|=|m|>0\).
\end{lemma}

\begin{proof}
The \(u_i\) are the roots of \(f(T)=\prod_i(T-u_i)=T^m-\epsilon_1T^{m-1}+\cdots+(-1)^m\epsilon_m\), whose coefficients \(\epsilon_i\) are polynomials with integer coefficients in the \(u_i\), hence lie in \(\mathcal O\), and whose constant term \(\pm\epsilon_m=\pm u_1\cdots u_m\) is a unit. Multiplying \(f(u_i)=0\) by \(u_i^n\) and summing over \(i\) gives, \(\forall n\in\Z\),
\[
\sigma_{n+m}=\epsilon_1\sigma_{n+m-1}-\epsilon_2\sigma_{n+m-2}+\cdots-(-1)^m\epsilon_m\sigma_n .
\]
Put \(N_j=\max(|\sigma_j|,\ldots,|\sigma_{j+m-1}|)\). The recurrence at \(n=j\) expresses \(\sigma_{j+m}\) as a combination of \(\sigma_j,\ldots,\sigma_{j+m-1}\) with coefficients of absolute value at most \(1\), so the ultrametric inequality gives \(|\sigma_{j+m}|\le N_j\), hence \(N_{j+1}\le N_j\). The same recurrence, solved for \(\sigma_j\), expresses \(\sigma_j\) as \(\epsilon_m^{-1}\) times a combination of \(\sigma_{j+1},\ldots,\sigma_{j+m}\) with coefficients of absolute value at most \(1\); since \(|\epsilon_m^{-1}|=1\), this gives \(|\sigma_j|\le N_{j+1}\), hence \(N_j\le N_{j+1}\). Thus \(N_j\) does not depend on \(j\). Every \(|\sigma_n|\) is at most \(N_n\), so the common value is the supremum, and it is attained. Finally \(\sigma_0=m\), and \(|m|>0\) because \(m\ne 0\) in \(K\).
\end{proof}

\begin{lemma}
\label{lem:transfer}
Let \(x\ge 2\) and let \(\omega_1,\ldots,\omega_6\in K^\times\) satisfy \(\ps{n}=0\) \(\forall n\in W\cup(-W)\). Suppose the absolute values \(|\omega_s|\) take exactly two values \(\rho_+>\rho_-\), and let \(I_+\) and \(I_-\) be the sets of indices where they are attained. Then \(I_+\) or \(I_-\) has at least four elements.
\end{lemma}

\begin{proof}
Suppose both classes have at most three elements. Choose \(s_+\in I_+\) and \(s_-\in I_-\), put \(\alpha=\omega_{s_+}\) and \(\beta=\omega_{s_-}\), and write \(\omega_s=\alpha u_s\) for \(s\in I_+\) and \(\omega_s=\beta v_s\) for \(s\in I_-\). The \(u_s\) and \(v_s\) are units. Let \(S(n)=\sum_{s\in I_+}u_s^n\), \(T(n)=\sum_{s\in I_-}v_s^n\), and \(\theta=|\beta|/|\alpha|=\rho_-/\rho_+\), so that \(0<\theta<1\). For \(n\in W\cup(-W)\) the hypothesis reads \(\alpha^nS(n)+\beta^nT(n)=0\), so \(|\alpha|^n|S(n)|=|\beta|^n|T(n)|\), that is
\begin{equation}
\label{eq:classratio}
|S(n)|=\theta^n|T(n)| .
\end{equation}
Let \(M_S=\sup_n|S(n)|\) and \(M_T=\sup_n|T(n)|\). By Lemma~\ref{lem:supnorm} both are finite and positive, and each is the maximum over any run of consecutive indices whose length is the number of elements of the class, hence also the maximum over any run of three consecutive indices, which contains such a run.

For \(n\in W\) one has \(n\ge x-1\ge 1\), so \(\theta^n\le\theta^{x-1}\) because \(\theta<1\), and~\eqref{eq:classratio} gives \(|S(n)|\le\theta^{x-1}M_T\). Since \(W\) is a run of three consecutive integers, \(M_S=\max_{n\in W}|S(n)|\le\theta^{x-1}M_T\). For \(n\in W\), the exponent \(-n\) lies in \(-W\), and~\eqref{eq:classratio} at \(-n\) reads \(|S(-n)|=\theta^{-n}|T(-n)|\), that is \(|T(-n)|=\theta^n|S(-n)|\le\theta^{x-1}M_S\). Since \(-W\) is a run of three consecutive integers, \(M_T=\max_{n\in W}|T(-n)|\le\theta^{x-1}M_S\). Combining, \(M_S\le\theta^{2(x-1)}M_S\). Since \(M_S>0\), this gives \(\theta^{2(x-1)}\ge 1\), which is false for \(0<\theta<1\) and \(x\ge 2\).
\end{proof}

The lemma uses only the two windows \(W\) and \(-W\), not the exponents \(\pm 1\), and it holds in every residue characteristic. It is the ultrametric counterpart of adding the two bounds in Lemma~\ref{lem:equalmoduli}: each window bounds the supremum of one class by a contraction of the supremum of the other.

\begin{lemma}
\label{lem:residues}
Let \(x\ge 2\) and let \(\omega_1,\ldots,\omega_6\in K^\times\) satisfy \(\ps{n}=0\) \(\forall n\in W\). Let \(\rho_+=\max_s|\omega_s|\) and \(I=\{s:|\omega_s|=\rho_+\}\), and suppose that \(I\) has at most three elements and is not all of \(\{1,\ldots,6\}\). Fix \(s_0\in I\) and put \(u_s=\omega_s/\omega_{s_0}\) for \(s\in I\). Then the \(u_s\) are units, and for every \(b\in\kappa\) the number of indices \(s\in I\) with \(\bar u_s=b\) is zero in \(\kappa\), that is, divisible by \(p\) when \(p\) is a prime. Consequently \(p\in\{2,3\}\); if \(p=3\) then \(I\) has exactly three elements and \(\bar u_s=1\) \(\forall s\in I\); if \(p=2\) then \(I\) has exactly two elements and \(\bar u_s=1\) for both.

The same holds with \(\rho_+\) replaced by \(\rho_-=\min_s|\omega_s|\) and \(I\) by \(\{s:|\omega_s|=\rho_-\}\), when \(\ps{n}=0\) \(\forall n\in-W\).
\end{lemma}

\begin{proof}
Fix \(n\in W\), so that \(n\ge 1\). For \(s\notin I\) one has \(|\omega_s^n|=|\omega_s|^n<\rho_+^n\), and the set of such \(s\) is nonempty, so by the ultrametric inequality \(\bigl|\sum_{s\notin I}\omega_s^n\bigr|<\rho_+^n\). From \(\ps{n}=0\),
\[
\Bigl|\sum_{s\in I}\omega_s^n\Bigr|=\Bigl|\sum_{s\notin I}\omega_s^n\Bigr|<\rho_+^n .
\]
Dividing by \(\omega_{s_0}^n\), of absolute value \(\rho_+^n\), gives \(\bigl|\sum_{s\in I}u_s^n\bigr|<1\). Each \(u_s\) has absolute value \(\rho_+/\rho_+=1\), so it is a unit, and the inequality says \(\sum_{s\in I}\bar u_s^{\,n}=0\) in \(\kappa\), for \(n=x-1,x,x+1\).

Let \(b_1,\ldots,b_t\) be the distinct elements of \(\kappa\) among the \(\bar u_s\), \(s\in I\), and let \(\mu_1,\ldots,\mu_t\) be their multiplicities, positive integers with \(\mu_1+\cdots+\mu_t=|I|\le 3\). Each \(b_j\) is nonzero, as the residue of a unit. The three vanishing residue sums read
\[
\sum_{j=1}^t\bigl(\mu_jb_j^{x-1}\bigr)\,b_j^{\,i}=0\qquad(i=0,1,2).
\]
Since \(t\le 3\), the first \(t\) of these equations form a linear system in the unknowns \(\mu_jb_j^{x-1}\) whose matrix \((b_j^{\,i})_{0\le i<t,\,1\le j\le t}\) is a Vandermonde matrix in the distinct elements \(b_j\), hence invertible over \(\kappa\). Therefore \(\mu_jb_j^{x-1}=0\) in \(\kappa\) \(\forall j\), and since \(b_j\ne 0\), the integer \(\mu_j\) is zero in \(\kappa\), which is impossible if \(p=0\). If \(p>0\), it says that \(p\) divides \(\mu_j\), and if \(p\ge 5\) this contradicts \(1\le\mu_j\le 3\). If \(p=3\), each \(\mu_j\) is a positive multiple of \(3\) and their sum is at most \(3\), so \(t=1\) and \(\mu_1=3\): all three residues coincide, and they equal \(\bar u_{s_0}=1\). If \(p=2\), each \(\mu_j\) is a positive even integer and their sum is at most \(3\), so \(t=1\) and \(\mu_1=2\): the two residues coincide and equal \(\bar u_{s_0}=1\).

For the minimum, apply the above to \(\omega_s^{-1}\), whose power sums at \(n\in W\) are \(\ps{-n}\), and whose largest absolute value \(\rho_-^{-1}\) is attained exactly on \(\{s:|\omega_s|=\rho_-\}\). The ratios formed there are the inverses \(u_s^{-1}\) of the ratios \(u_s=\omega_s/\omega_{s_0}\), and a unit has residue \(1\) iff its inverse does, so the conclusion holds for the \(u_s\) as stated.
\end{proof}

\begin{lemma}
\label{lem:classes}
Let \(x\ge 2\) and let \(\omega_1,\ldots,\omega_6\in K^\times\) satisfy the window condition at \(x\). Suppose the \(\omega_s\) are not all units. Let \(\rho_+\) and \(\rho_-\) be the largest and the smallest absolute value and \(I_+\), \(I_-\) the sets of indices where they are attained. Then \(\rho_+>1>\rho_-\), the sets \(I_+\) and \(I_-\) are disjoint, and \(p\in\{2,3\}\). If \(p=3\), then \(I_+\) and \(I_-\) have three elements each, so the absolute values take exactly two values. If \(p=2\), then at least one of \(I_+\) and \(I_-\) has at most three elements, and each of them that has at most three elements is a pair \(\{s_0,s_1\}\) with \(\omega_{s_1}=\omega_{s_0}w\) for a unit \(w\) with \(\bar w=1\).
\end{lemma}

\begin{proof}
If \(\rho_+\le 1\), then every \(|\omega_s|\le 1\) and \(\prod_s|\omega_s|=1\) force every \(|\omega_s|=1\), against the assumption; so \(\rho_+>1\), and symmetrically \(\rho_-<1\). Hence \(\rho_+\ne\rho_-\) and \(I_+\cap I_-=\emptyset\), so \(|I_+|+|I_-|\le 6\) and one of the two sets has at most three elements. Neither is all of \(\{1,\ldots,6\}\). Lemma~\ref{lem:residues}, applied to the set with at most three elements, with the window \(W\) for \(I_+\) and \(-W\) for \(I_-\), gives \(p\in\{2,3\}\). If \(p=3\), that set has three elements, so the other has at most three as well and Lemma~\ref{lem:residues} makes it a set of three elements too; the two sets exhaust the indices. If \(p=2\), Lemma~\ref{lem:residues} says that each of \(I_+\), \(I_-\) with at most three elements is a pair \(\{s_0,s_1\}\) with \(\bar u_{s_1}=1\), where \(u_{s_1}=\omega_{s_1}/\omega_{s_0}\); put \(w=u_{s_1}\).
\end{proof}

\begin{proposition}
\label{prop:charthree}
Let \(x\ge 2\) and suppose that \(p=3\). Then elements \(\omega_1,\ldots,\omega_6\in K^\times\) satisfying the window condition at \(x\) are units.
\end{proposition}

\begin{proof}
Otherwise Lemma~\ref{lem:classes} gives exactly two absolute values, attained on two classes of three elements each, which Lemma~\ref{lem:transfer} excludes.
\end{proof}

For \(p\notin\{2,3\}\) the conclusion of Proposition~\ref{prop:charthree} is part of Lemma~\ref{lem:classes}. The case \(p=2\) requires the exponents \(\pm 1\) and a finer analysis. Two facts about a pair are used repeatedly. Let \(y\) be a unit and \(n\in\Z\). First, \(1+y^{-n}=y^{-n}(1+y^n)\), so \(|1+y^{-n}|=|1+y^n|\). Second, if \(p=2\) and \(\bar y=1\), then for every odd \(n\ge 1\),
\begin{equation}
\label{eq:oddpair}
|1+y^n|=|1+y|,
\end{equation}
because \(1+y^n=(1+y)(1-y+y^2-\cdots+y^{n-1})\) and the second factor has residue \(1-1+1-\cdots+1=\bar n=1\) in \(\kappa\), the number of terms being odd; so it is a unit. By the first fact,~\eqref{eq:oddpair} holds for every odd \(n\in\Z\).

\begin{lemma}
\label{lem:twopatterns}
Let \(x\ge 2\), suppose \(p=2\), and let \(\omega_1,\ldots,\omega_6\in K^\times\) satisfy the window condition at \(x\) and not be all units. List the distinct absolute values in decreasing order and the sizes of the corresponding classes in the same order. Then, after possibly replacing every \(\omega_s\) by \(\omega_s^{-1}\), the list of sizes is \((2,4)\), \((2,2,2)\) or \((2,1,1,2)\). In every extreme class of size two, the two numbers can be written as \(a\) and \(aw\), where \(w\) is a unit with \(\bar w=1\).
\end{lemma}

\begin{proof}
By Lemma~\ref{lem:classes}, one of the two extreme classes has at most three elements and is then a pair of the stated form. The other extreme class has at most four elements, since the two are disjoint, and if it has at most three elements then Lemma~\ref{lem:residues} makes it a pair; so it has two or four elements. If it has four, the two classes exhaust the indices and the sizes are \((2,4)\) or \((4,2)\); inversion exchanges the two orders. If it has two, it is a pair of the stated form by Lemma~\ref{lem:classes}, and the two remaining indices carry absolute values strictly between \(\rho_-\) and \(\rho_+\), either equal, giving \((2,2,2)\), or distinct, giving \((2,1,1,2)\).
\end{proof}

\begin{proposition}
\label{prop:smallpatterns}
In the setting of Lemma~\ref{lem:twopatterns}, the patterns \((2,1,1,2)\) and \((2,2,2)\) do not occur.
\end{proposition}

\begin{proof}
Pattern \((2,1,1,2)\). Write the six numbers as \(\alpha,\alpha w_A\), then \(\beta\), then \(\gamma\), then \(\delta,\delta w_D\), with absolute values \(\rho_1>\rho_2>\rho_3>\rho_4\) and \(\bar w_A=\bar w_D=1\), and put \(\lambda_A=|1+w_A|\). Let \(n\in\{1,n_o\}\); both are odd exponents of the window condition. The four terms of \(\ps{n}=\alpha^n(1+w_A^n)+\beta^n+\gamma^n+\delta^n(1+w_D^n)\) have absolute values \(\rho_1^n\lambda_A\) by~\eqref{eq:oddpair}, \(\rho_2^n\), \(\rho_3^n\), and \(\rho_4^n|1+w_D^n|\le\rho_4^n\). Among the last three, \(\rho_2^n\) is strictly the largest. Since \(\ps{n}=0\), the largest absolute value is attained twice, which is possible only if \(\rho_1^n\lambda_A=\rho_2^n\): if \(\rho_1^n\lambda_A\) were larger, the first term would be a strict maximum, and if it were smaller, the second term would be. Thus \(\rho_1\lambda_A=\rho_2\) and \(\rho_1^{n_o}\lambda_A=\rho_2^{n_o}\). Dividing the second relation by the first gives \(\rho_1^{n_o-1}=\rho_2^{n_o-1}\) with \(n_o-1\ge 2\), so \(\rho_1=\rho_2\), a contradiction.

Pattern \((2,2,2)\). Write the six numbers as \(\alpha,\alpha w_A\), then \(\mu,\mu w_B\), then \(\delta,\delta w_C\), with absolute values \(\rho_1>\rho_2>\rho_3\) and \(\bar w_A=\bar w_C=1\). Put \(\lambda_A=|1+w_A|\), \(\lambda_B=|1+w_B|\), \(\lambda_C=|1+w_C|\), \(\theta_1=\rho_2/\rho_1<1\) and \(\theta_2=\rho_3/\rho_2<1\). For odd \(n\), the \(A\) term \(\alpha^n(1+w_A^n)\) of \(\ps{n}\) has absolute value \(\rho_1^n\lambda_A\) and the \(C\) term has absolute value \(\rho_3^n\lambda_C\), by~\eqref{eq:oddpair}; the \(B\) term has absolute value \(\rho_2^n|1+w_B^n|\le\rho_2^n\).

First, \(\ps{n_o}=0\). The \(B\) term has absolute value at most \(\rho_2^{n_o}\) and the \(C\) term at most \(\rho_3^{n_o}<\rho_2^{n_o}\). If \(\rho_1^{n_o}\lambda_A>\rho_2^{n_o}\), the \(A\) term is a strict maximum, which is impossible. Hence \(\lambda_A\le\theta_1^{n_o}\). Likewise \(\ps{-n_o}=0\), where the \(A\) term has absolute value \(\rho_1^{-n_o}\lambda_A\le\rho_1^{-n_o}<\rho_2^{-n_o}\), the \(B\) term at most \(\rho_2^{-n_o}\), and the \(C\) term \(\rho_3^{-n_o}\lambda_C\); if \(\rho_3^{-n_o}\lambda_C>\rho_2^{-n_o}\) the \(C\) term is a strict maximum. Hence \(\lambda_C\le\theta_2^{n_o}\).

Second, \(\ps{1}=0\). The \(A\) term has absolute value \(\rho_1\lambda_A\le\rho_1\theta_1^{n_o}=\rho_2\theta_1^{n_o-1}<\rho_2\), and the \(C\) term has absolute value \(\rho_3\lambda_C\le\rho_3<\rho_2\). If \(\lambda_B=1\), the \(B\) term, of absolute value \(\rho_2\), is a strict maximum, which is impossible. Hence \(\lambda_B<1\), that is, \(\overline{w_B}=-1=1\) in \(\kappa\), and~\eqref{eq:oddpair} applies to \(w_B\) as well: for odd \(n\), the \(B\) term has absolute value \(\rho_2^n\lambda_B\). For odd \(n\) write \(x_A(n)=\rho_1^n\lambda_A\), \(x_B(n)=\rho_2^n\lambda_B\), \(x_C(n)=\rho_3^n\lambda_C\) for the three absolute values; by the first fact about pairs, these formulas hold for negative odd \(n\) as well.

Third, the case where some \(\lambda\) vanishes, that is, some \(w\) equals \(-1\). If all three vanish, the six numbers are \(\pm\alpha,\pm\mu,\pm\delta\), and \(\ps{n_e}=2(\alpha^{n_e}+\mu^{n_e}+\delta^{n_e})\); the three summands have the distinct absolute values \(\rho_1^{n_e}>\rho_2^{n_e}>\rho_3^{n_e}\), so their sum has absolute value \(\rho_1^{n_e}>0\), and \(|2|>0\) because \(K\) has characteristic \(0\); thus \(\ps{n_e}\ne 0\), a contradiction. If exactly two vanish, \(\ps{1}\) has one nonzero term and cannot vanish. If exactly one vanishes, say for the class \(i\), then at \(n=1\) and at \(n=n_o\) the two remaining terms, both nonzero, must have equal absolute values: \(\rho_j\lambda_j=\rho_k\lambda_k\) and \(\rho_j^{n_o}\lambda_j=\rho_k^{n_o}\lambda_k\) for the other two classes \(j\ne k\). Dividing gives \((\rho_j/\rho_k)^{n_o-1}=1\), so \(\rho_j=\rho_k\), a contradiction.

Fourth, the case where all three \(\lambda\) are positive. For odd \(n\) and \(i\in\{A,B,C\}\) put \(f_i(n)=\log x_i(n)\), so that \(f_A(n)=n\log\rho_1+\log\lambda_A\), \(f_B(n)=n\log\rho_2+\log\lambda_B\) and \(f_C(n)=n\log\rho_3+\log\lambda_C\). These are affine functions of a real variable \(n\) with distinct slopes. Their maximum \(F\) is convex and piecewise affine, and each slope occurs on a single interval, so \(F\) changes slope at no more than two points. If two of the \(f_i\) attain \(F\) at some \(n\), then \(F\) changes slope at \(n\), because the two have different slopes. Now \(\ps{n}=0\) at the three odd exponents \(-1\), \(1\) and \(n_o\), so at each of them the largest of \(x_A(n)\), \(x_B(n)\), \(x_C(n)\) is attained at least twice, and \(F\) changes slope at three distinct points. This is a contradiction, and all cases are excluded.
\end{proof}

\begin{lemma}
\label{lem:newtonparity}
Let \(e_1,e_2,e_3,e_4\) be the elementary symmetric functions of four indeterminates and \(p_m\) their power sums. For every odd \(m\ge 1\), \(p_m\) lies in the ideal generated by \(e_1\) and \(e_3\) in \(\Z[e_1,e_2,e_3,e_4]\). Consequently, if \(v_1,\ldots,v_4\) are units of \(\mathcal O\) and \(T_n=\sum_{j=1}^4v_j^n\), then \(|T_m|\le\max(|T_1|,|T_{-1}|)\) for every odd integer \(m\).
\end{lemma}

\begin{proof}
Newton's identities~\cite[Chapter~I, \S 2]{macdonald1995} give \(p_1=e_1\), \(p_2=e_1p_1-2e_2\), \(p_3=e_1p_2-e_2p_1+3e_3=e_1^3-3e_1e_2+3e_3\), and for \(m\ge 5\),
\[
p_m=e_1p_{m-1}-e_2p_{m-2}+e_3p_{m-3}-e_4p_{m-4}.
\]
Thus \(p_1\) and \(p_3\) lie in the ideal \((e_1,e_3)\). For odd \(m\ge 5\), the first and third terms on the right lie in \((e_1,e_3)\) because of their factors \(e_1\) and \(e_3\), and the second and fourth lie there by induction, \(m-2\) and \(m-4\) being odd and positive. This proves the first assertion.

For the units, write \(p_m=Fe_1+Ge_3\) with \(F,G\in\Z[e_1,\ldots,e_4]\), and evaluate at \(v_1,\ldots,v_4\). The values of \(e_1,\ldots,e_4\) lie in \(\mathcal O\), hence so do those of \(F\) and \(G\), and the ultrametric inequality gives \(|T_m|\le\max(|e_1(v)|,|e_3(v)|)\) for odd \(m\ge 1\). Here \(e_1(v)=T_1\), and \(e_3(v)=e_4(v)\sum_jv_j^{-1}=e_4(v)\,T_{-1}\) with \(|e_4(v)|=|v_1v_2v_3v_4|=1\), so \(|e_3(v)|=|T_{-1}|\). For odd \(m\le -1\), apply the case of positive exponents to the units \(v_j^{-1}\), whose power sum at \(-m\) is \(T_m\) and whose power sums at \(1\) and \(-1\) are \(T_{-1}\) and \(T_1\).
\end{proof}

The ideal membership over \(\Z\) gives a uniform bound for all odd power sums of four units in terms of \(T_1\) and \(T_{-1}\). If both sums vanish, then \(e_1=e_3=0\), so the quartic is even and the four units form two antipodal pairs.

\begin{proposition}
\label{prop:pattern24}
In the setting of Lemma~\ref{lem:twopatterns}, the pattern \((2,4)\) does not occur.
\end{proposition}

\begin{proof}
In the pattern \((2,4)\) the pair is an extreme class, so by Lemma~\ref{lem:twopatterns} its two numbers can be written as \(a\) and \(au\), where \(u\) is a unit with \(\bar u=1\). Replacing every \(\omega_s\) by its inverse if necessary, which preserves the window condition and the property \(\bar u=1\), we may assume that the pair carries the smaller absolute value. Write the quadruple as \(bv_1,\ldots,bv_4\) with \(v_j\) units, put \(\rho_1=|a|<\rho_2=|b|\), \(\theta=\rho_1/\rho_2<1\), \(a(n)=|1+u^n|\) and \(T_n=\sum_jv_j^n\). For \(n\in\{1\}\cup W\) the relation \(\ps{n}=0\) reads \(a^n(1+u^n)=-b^nT_n\), so \(\rho_1^na(n)=\rho_2^n|T_n|\), that is \(|T_n|=\theta^na(n)\); and \(\ps{-n}=0\) reads \(a^{-n}(1+u^{-n})=-b^{-n}T_{-n}\), so \(|T_{-n}|=\theta^{-n}a(-n)=\theta^{-n}a(n)\), by the first fact about pairs.

Put \(a_o=a(1)\). By~\eqref{eq:oddpair}, \(a(n_o)=a_o\). Thus \(|T_1|=\theta a_o\), \(|T_{-1}|=\theta^{-1}a_o\) and \(|T_{-n_o}|=\theta^{-n_o}a_o\). Lemma~\ref{lem:newtonparity} at the odd exponent \(-n_o\) gives
\[
\theta^{-n_o}a_o=|T_{-n_o}|\le\max(\theta a_o,\theta^{-1}a_o)=\theta^{-1}a_o .
\]
If \(a_o>0\), dividing by \(a_o\) gives \(\theta^{1-n_o}\le 1\), which is false because \(\theta<1\) and \(1-n_o<0\). Hence \(a_o=0\), that is, \(u=-1\). This is consistent with \(\bar u=1\), since \(p=2\). The preceding contradiction also excludes \(u=1\), where \(a_o=|2|>0\).

With \(u=-1\), the relations at \(n=\pm 1\) give \(T_1=0\) and \(T_{-1}=0\). Hence \(e_1(v)=T_1=0\) and \(e_3(v)=e_4(v)T_{-1}=0\), and the \(v_j\) are the roots of \(z^4+e_2(v)z^2+e_4(v)\). This polynomial is even, so the multiset of its roots is invariant under \(z\mapsto-z\): the \(v_j\) are \(s,-s,t,-t\) for some \(s,t\in K^\times\), with \(e_4(v)=s^2t^2\), so \(|st|=1\). For every integer \(k\), \(T_{2k}=2(s^{2k}+t^{2k})\) and \(T_{-2k}=2(s^{-2k}+t^{-2k})=T_{2k}/(st)^{2k}\), so \(|T_{-2k}|=|T_{2k}|\). At the even exponent \(n_e\in W\) one has \(a(n_e)=|1+(-1)^{n_e}|=|2|\), and the two relations give
\[
\theta^{n_e}|2|=|T_{n_e}|=|T_{-n_e}|=\theta^{-n_e}|2| .
\]
Since \(|2|>0\), this gives \(\theta^{2n_e}=1\), so \(\theta=1\), contradicting \(\theta<1\).
\end{proof}

\begin{proposition}
\label{prop:nonarch}
Let \(x\ge 2\), and let \(K\) be a field of characteristic \(0\) with a non-archi\-medean absolute value. Elements \(\omega_1,\ldots,\omega_6\in K^\times\) satisfying the window condition at \(x\) have absolute value \(1\).
\end{proposition}

\begin{proof}
Suppose not. Lemma~\ref{lem:classes} gives \(p\in\{2,3\}\). Proposition~\ref{prop:charthree} excludes \(p=3\). For \(p=2\), Lemma~\ref{lem:twopatterns} leaves the three patterns, which Propositions~\ref{prop:smallpatterns} and~\ref{prop:pattern24} exclude.
\end{proof}

The two local steps combine, through Kronecker's theorem, to show that algebraic solutions of the window condition are roots of unity, and the classification of vanishing sums of roots of unity then gives the symmetry.

\begin{lemma}
\label{lem:conjugates}
Let \(x\ge 2\) and let \(\omega_1,\ldots,\omega_6\) be algebraic numbers satisfying the window condition at \(x\). Then every \(\omega_s\) is a root of unity.
\end{lemma}

\begin{proof}
Let \(K=\Q(\omega_1,\ldots,\omega_6)\), a number field. The window condition consists of polynomial identities with integer coefficients in the \(\omega_s\) and their inverses, so it is preserved by every field homomorphism from \(K\) into a field.

Let \(\mathfrak p\) be a nonzero prime ideal of the ring of integers \(\mathcal O_K\). The \(\mathfrak p\)-adic absolute value on \(K\) is non-archimedean and \(K\) has characteristic \(0\), so Proposition~\ref{prop:nonarch} gives \(|\omega_s|_{\mathfrak p}=1\), that is, \(\omega_s\) is a unit in the localization of \(\mathcal O_K\) at \(\mathfrak p\). An element of \(K\) that lies in the localization of \(\mathcal O_K\) at every nonzero prime ideal lies in \(\mathcal O_K\)~\cite[Chapter~I, \S 11]{neukirch1999}. Hence every \(\omega_s\) is an algebraic integer.

Let \(\beta\) be an algebraic conjugate of \(\omega_s\). The isomorphism \(\Q(\omega_s)\to\Q(\beta)\) sending \(\omega_s\) to \(\beta\) extends to a field embedding \(\tau\colon K\to\C\), since \(K\) is a finite extension of \(\Q(\omega_s)\) and \(\C\) is algebraically closed~\cite[Chapter~V, Theorem~2.8]{lang2002}. The numbers \(\tau(\omega_1),\ldots,\tau(\omega_6)\) are nonzero complex numbers satisfying the window condition at \(x\), so Corollary~\ref{cor:unimodular} gives \(|\tau(\omega_s)|=1\), that is, \(|\beta|=1\). Every conjugate of the algebraic integer \(\omega_s\) has modulus \(1\), and Kronecker's theorem~\cite{kronecker1857} says that \(\omega_s\) is a root of unity.
\end{proof}

\begin{lemma}
\label{lem:rootsofunity}
Let \(x\ge 2\) and let \(\omega_1,\ldots,\omega_6\) be roots of unity, not necessarily distinct, with \(\ps{1}=\ps{x-1}=\ps{x}=\ps{x+1}=0\). Then the multiset \(\{\omega_1,\ldots,\omega_6\}\) is invariant under \(z\mapsto-z\), in which case \(\es{3}=0\), or under \(z\mapsto\zroot{3}z\), in which case \(\es{2}=\es{4}=0\).
\end{lemma}

\begin{proof}
In the terminology recalled before Proposition~\ref{prop:sixroots}, the identity \(\ps{1}=0\) is a vanishing relation of weight six among the distinct roots of unity that occur among the \(\omega_s\), the coefficient of each being its multiplicity. Every relation is a sum of minimal relations~\cite{poonen-rubinstein1998}: if the relation is not minimal, it has a nontrivial proper subrelation, whose coefficients are at most the given ones, and subtracting it leaves a relation of smaller weight; induction on the weight gives a decomposition of the multiset into sub-multisets, each of which sums to zero and is minimal. The minimal relations of weight at most six are, up to rotation, \(R_2\), \(R_3\), \(R_5\) and \((R_5:R_3)\), as recalled there from Mann's theorem and from Poonen and Rubinstein's Theorem~3~\cite{mann1965,poonen-rubinstein1998}; each of them is a sum of distinct roots with coefficients \(1\). Since there is no minimal relation of weight \(1\) or \(4\), the weight six splits as \(2+2+2\), \(3+3\) or \(6\).

If the split is \(6\), the multiset is \(\rho T\) for a root of unity \(\rho\) and the set \(T\) of Lemma~\ref{lem:exotic}, and that lemma says that \(\ps{t}=\sum_{\eta\in\rho T}\eta^t\) vanishes only when \(t\equiv\pm 1\pmod 6\). One of \(x-1\), \(x\), \(x+1\) is a positive multiple of \(3\), and it is congruent to \(0\) or \(3\) modulo \(6\); its power sum is nonzero, against the hypothesis. If the split is \(2+2+2\), the multiset is \(\{a_1,-a_1,a_2,-a_2,a_3,-a_3\}\), invariant under \(z\mapsto-z\). If the split is \(3+3\), it is \(\{a_1,a_1\zroot{3},a_1\zroot{3}^2,a_2,a_2\zroot{3},a_2\zroot{3}^2\}\), invariant under \(z\mapsto\zroot{3}z\). For distinct roots, the splits \(2+2+2\), \(3+3\) and \(6\) are, in this order, the three configurations of Figure~\ref{fig:sixconfigs}.

Replacing every \(\omega_s\) by \(c\,\omega_s\) multiplies \(\es{j}\) by \(c^j\). If the multiset is invariant under \(z\mapsto-z\), then \(\es{3}=(-1)^3\es{3}\), so \(\es{3}=0\). If it is invariant under \(z\mapsto\zroot{3}z\), then \(\es{2}=\zroot{3}^2\es{2}\) and \(\es{4}=\zroot{3}^4\es{4}=\zroot{3}\es{4}\), and since \(\zroot{3}^2\ne 1\ne\zroot{3}\), both vanish.
\end{proof}

It remains to pass from algebraic solutions to complex ones. The desired coefficient identities are polynomial, so the following consequence of the Nullstellensatz applies.

\begin{lemma}
\label{lem:descent}
Let \(I\) be an ideal of \(\Q[X_1,\ldots,X_k]\) and \(f\in\Q[X_1,\ldots,X_k]\). If \(f\) vanishes at every common zero of \(I\) in \(\overline{\Q}^{\,k}\), then \(f\) vanishes at every common zero of \(I\) in \(\C^k\).
\end{lemma}

\begin{proof}
By Hilbert's Nullstellensatz over the algebraically closed field \(\overline\Q\)~\cite[Chapter~IX]{lang2002}, some power \(f^M\) lies in the ideal generated by \(I\) in the ring \(\overline\Q[X_1,\ldots,X_k]\), say \(f^M=\sum_ig_ih_i\) with \(g_i\in I\) and \(h_i\in\overline\Q[X_1,\ldots,X_k]\). Choose a basis \((b_\lambda)\) of \(\overline\Q\) as a vector space over \(\Q\) with one basis element \(b_0=1\). Expanding the coefficients of each \(h_i\) in this basis gives \(h_i=\sum_\lambda h_{i\lambda}b_\lambda\) with \(h_{i\lambda}\in\Q[X_1,\ldots,X_k]\), finitely many of them nonzero, and the expansion of a polynomial in \(\overline\Q[X_1,\ldots,X_k]\) as \(\sum_\lambda q_\lambda b_\lambda\) with \(q_\lambda\in\Q[X_1,\ldots,X_k]\) is unique, coefficient by coefficient. Now \(f^M=\sum_\lambda\bigl(\sum_ig_ih_{i\lambda}\bigr)b_\lambda\), while the left side has coefficients in \(\Q\) and so has only the term \(\lambda=0\) in its expansion. Comparing the two expansions gives \(f^M=\sum_ig_ih_{i0}\in I\). Hence \(f^M\), and with it \(f\), vanishes at every common zero of \(I\) in \(\C^k\).
\end{proof}

\begin{proof}[Proof of Theorem~\ref{thm:sixnumbers}]
Let \(\omega_1,\ldots,\omega_6\) be nonzero complex numbers satisfying the window condition at \(x\), and let \(c=\es{2}\), \(r=\es{3}\), \(d=\es{4}\), so that \(\prod_s(z-\omega_s)=z^6+cz^4-rz^3+dz^2+1\) by~\eqref{eq:windowshape}. For \(n\ge 1\), Newton's identities express \(\ps{n}\) as a polynomial with integer coefficients in \(\es{1},\ldots,\es{6}\); substituting \(\es{1}=\es{5}=0\) and \(\es{6}=1\) gives a polynomial \(P_n\in\Z[c,r,d]\) with \(\ps{n}=P_n(c,r,d)\). For \(n\ge 1\), the power sum \(\ps{-n}\) is the \(n\)th power sum of the inverses \(\omega_s^{-1}\), whose elementary symmetric functions are \(\es{6-j}/\es{6}\), namely \(0,d,r,c,0,1\); hence \(\ps{-n}=P_n(d,r,c)\). Let \(I\subset\Q[c,r,d]\) be the ideal generated by the six polynomials \(P_n(c,r,d)\) and \(P_n(d,r,c)\) for \(n\in W\). The point \((c,r,d)\) is a common zero of \(I\) in \(\C^3\).

Conversely, let \((c_0,r_0,d_0)\in\overline\Q^{\,3}\) be a common zero of \(I\), and let \(\omega_1',\ldots,\omega_6'\) be the roots of \(z^6+c_0z^4-r_0z^3+d_0z^2+1\), with multiplicity. They are algebraic, nonzero because the constant term is \(1\), with product \(1\), and their elementary symmetric functions are \(0,c_0,r_0,d_0,0,1\), so \(\ps{1}=\es{1}=0\), \(\ps{-1}=\es{5}/\es{6}=0\), and \(\ps{\pm n}=0\) for \(n\in W\) by the expressions for \(\ps{n}\) and \(\ps{-n}\) above. Thus the \(\omega_s'\) satisfy the window condition at \(x\). By Lemma~\ref{lem:conjugates} they are roots of unity, and by Lemma~\ref{lem:rootsofunity} either \(r_0=0\) or \(c_0=d_0=0\). In both cases \(r_0c_0=0\) and \(r_0d_0=0\).

The polynomials \(rc\) and \(rd\) therefore vanish at every common zero of \(I\) in \(\overline\Q^{\,3}\), and Lemma~\ref{lem:descent} says that they vanish at every common zero of \(I\) in \(\C^3\), in particular at \((c,r,d)\). So \(\es{3}\es{2}=0\) and \(\es{3}\es{4}=0\), which is the assertion that \(\es{3}=0\) or \(\es{2}=\es{4}=0\). The second statement of the theorem was derived from the first after its statement.
\end{proof}

\begin{corollary}
\label{cor:sixrootsofunity}
Let \(x\ge 2\), let \(\omega_1,\ldots,\omega_6\) be nonzero complex numbers satisfying the window condition at \(x\), and let \(n_3\) be the element of \(W\) divisible by \(3\). If \(\es{3}=0\), then \(\omega_s^{6n_e}=1\) \(\forall s\). If \(\es{2}=\es{4}=0\), then \(\omega_s^{4n_3}=1\) \(\forall s\). In particular every \(\omega_s\) is a root of unity of order at most \(6(x+1)\), and the ratios of six unimodular numbers with sum zero and \(\ps{x-1}=\ps{x}=\ps{x+1}=0\) are roots of unity.
\end{corollary}

\begin{proof}
By Theorem~\ref{thm:sixnumbers}, \(\es{3}=0\) or \(\es{2}=\es{4}=0\). Suppose first that \(\es{3}=0\). By~\eqref{eq:windowshape} the polynomial \(\prod_s(z-\omega_s)=z^6+\es{2}z^4+\es{4}z^2+1\) is even, so the multiset of the \(\omega_s\) is \(\{\pm\alpha,\pm\beta,\pm\gamma\}\) for some nonzero \(\alpha\), \(\beta\), \(\gamma\), and \(1=\prod_s\omega_s=-(\alpha\beta\gamma)^2\). Put \(y_1=\alpha^{n_e}\), \(y_2=\beta^{n_e}\) and \(y_3=\gamma^{n_e}\). Since \(n_e\) is even, \(\ps{n_e}=2(y_1+y_2+y_3)\) and \(\ps{-n_e}=2(y_1^{-1}+y_2^{-1}+y_3^{-1})\), and both vanish. The elementary symmetric functions of \(y_1,y_2,y_3\) are therefore \(0\), \(0\) and \(y_1y_2y_3\), so every \(y_i\) satisfies \(y_i^3=y_1y_2y_3=\bigl((\alpha\beta\gamma)^2\bigr)^{n_e/2}=(-1)^{n_e/2}\). Hence \(\alpha^{6n_e}=y_1^6=1\), and in the same way \(\beta^{6n_e}=\gamma^{6n_e}=1\). Since \(6n_e\) is even, the same holds for \(-\alpha\), \(-\beta\) and \(-\gamma\).

Suppose now that \(\es{2}=\es{4}=0\). Then \(\prod_s(z-\omega_s)=z^6-\es{3}z^3+1\), which factors as \((z^3-\tau_1)(z^3-\tau_2)\) with \(\tau_1\tau_2=1\). The multiset of the \(\omega_s\) is therefore \(\{\alpha,\alpha\zroot{3},\alpha\zroot{3}^2,\beta,\beta\zroot{3},\beta\zroot{3}^2\}\) with \(\alpha^3=\tau_1\), \(\beta^3=\tau_2\), and \((\alpha\beta)^3=1\). Since \(3\) divides \(n_3\), one has \(\ps{n_3}=3(\alpha^{n_3}+\beta^{n_3})\), which vanishes, so \(\beta^{n_3}=-\alpha^{n_3}\). Multiplying by \(\alpha^{n_3}\) gives \(-\alpha^{2n_3}=(\alpha\beta)^{n_3}\), which is \(1\) because \((\alpha\beta)^3=1\) and \(3\) divides \(n_3\). Hence \(\alpha^{4n_3}=1\), then \(\beta^{4n_3}=\alpha^{4n_3}=1\), and \(\zroot{3}^{4n_3}=1\).

The orders divide \(6n_e\) or \(4n_3\), and both are at most \(6(x+1)\). For six unimodular numbers with sum zero and \(\ps{x-1}=\ps{x}=\ps{x+1}=0\), the rotation that makes their product \(1\) produces numbers satisfying the window condition, as shown after the statement of Theorem~\ref{thm:sixnumbers}. The rotated numbers are roots of unity, and a rotation does not change the ratios.
\end{proof}

With Lemma~\ref{lem:normalform}, Theorem~\ref{thm:sixnumbers} and Corollary~\ref{cor:sixrootsofunity} give Theorem~\ref{thm:introsix} of the introduction.

\begin{remark}
\label{rem:sixnumbers}
Let \(S\subset\R/\Z\) have six elements and let \(\omega_s=e^{2\pi is}\) for \(s\in S\). If \(\widehat{\muS}(n)=\sum_s\omega_s^n\) vanishes at \(n=1\), \(x-1\), \(x\), \(x+1\) for some integer \(x\ge 2\), then Theorem~\ref{thm:sixnumbers} gives \(\es{2}\es{3}=0\), and Lemma~\ref{lem:normalform}, after the rotation that makes the product of the \(\omega_s\) equal to \(1\), gives \(S+\tfrac12=S\) or \(S+\tfrac13=S\). The exceptional fiber \(A=\{0,1,x,x+1,y,z\}\) considered in the proof of Corollary~\ref{cor:k4} has \(1\), \(x-1\), \(x\) and \(x+1\) among its differences, so for that fiber the conclusion of Lemma~\ref{lem:sixpoint} follows from Theorem~\ref{thm:sixnumbers} without Theorem~\ref{thm:rational}: the four power sums carry the whole obstruction, and the remaining differences of \(A\), as well as the sizes of the gaps, play no role. Corollary~\ref{cor:sixrootsofunity} gives, moreover, \(S-S\subset\Q/\Z\), the conclusion of Theorem~\ref{thm:rational} for this \(S\), from four of its orthogonality relations.
\end{remark}

\section{A weak tiling that is not a convex combination of tilings}
\label{sec:weak}

\begin{definition}
\label{def:weaktile}
A bounded measurable set \(\Omega\subset\R\) weakly tiles its complement when there exists a positive locally finite measure \(\nu\) on \(\R\) such that \(1_\Omega*\nu=1_{\Omega^c}\) almost everywhere. The measure \(\nu\) is a weak tiling measure for \(\Omega\). If \(\nu\) is a sum of distinct unit point masses, the copies of \(\Omega\) are disjoint and fill the complement, and the weak tiling is a proper tiling of the complement.
\end{definition}

The definition is that of Lev and Matolcsi \cite{lev-matolcsi2022}. The same authors prove that every bounded spectral set weakly tiles its complement \cite[Theorem~1.5]{lev-matolcsi2022}.

Kolountzakis, Lev and Matolcsi ask two questions about a finite union of intervals \cite{klm2025}. Problem~4.2 asks whether a weak tiling of the complement forces a proper tiling of the real line, for a union of three or more intervals. Problem~4.3 asks whether every weak tiling measure of a finite union of intervals is a convex combination of proper tilings of the complement. An affirmative answer to the second question would imply an affirmative answer to the first. Theorem~\ref{thm:weak} answers the second question in the negative for four intervals. The set in the statement tiles the real line, and it is spectral, so the example does not bear on Problem~4.2.

For a pure point measure write \(\mu(t)\) for the mass at \(t\). Differentiating the tiling equation gives the following balance between the masses at left and right endpoints.

\begin{lemma}
\label{lem:balance}
Let \(\Omega=\bigcup_{j=1}^{4}(a_j,b_j)\) be a union of four pairwise disjoint bounded open intervals, and let \(\mu\) be a positive pure point measure, locally finite, with \(1_\Omega*\mu=1\) almost everywhere. Then \(\forall y\in\R\),
\[
\sum_{j=1}^{4}\mu(y-a_j)=\sum_{j=1}^{4}\mu(y-b_j).
\]
Conversely, if \(\mu\) is such a measure, \(\mu(\{0\})=1\), the displayed identity holds for every \(y\), and \(\Omega\cap(\Omega+t)\) is null for every atom \(t\neq 0\), then \(1_\Omega*\mu=1\) almost everywhere.
\end{lemma}

\begin{proof}
A locally finite measure may have infinitely many atoms in a bounded window, so the argument is made with distributions. The convolution \(F=1_\Omega*\mu\) is a locally integrable function, given at every point by \(F(x)=\sum_t\mu(t)\,1_\Omega(x-t)\), because \(1_\Omega\) is bounded with compact support and \(\mu\) is locally finite. Since \(1_\Omega\) has compact support, its distributional derivative may be convolved with \(\mu\), and
\[
F'=1_\Omega'*\mu=\sum_{j=1}^{4}\bigl(\delta_{a_j}-\delta_{b_j}\bigr)*\mu.
\]
This is a locally finite signed pure point measure, and its mass at \(y\) is \(\sum_j\mu(y-a_j)-\sum_j\mu(y-b_j)\). If \(F=1\) almost everywhere, then \(F'=0\), and every mass vanishes, which is the displayed identity. Conversely, if every mass vanishes, then \(F'=0\), so \(F\) agrees almost everywhere with a constant. The measure \(\mu\) has countably many atoms, and \(\Omega\cap(\Omega+t)\) is null for each atom \(t\neq 0\). For almost every \(x\in\Omega\), therefore, only the atom at \(0\) contributes to \(F(x)\), and \(F(x)=\mu(\{0\})=1\). The constant is \(1\).
\end{proof}

\begin{theorem}
\label{thm:weak}
Let \(\Omega=(0,1)\cup(4,5)\cup(6,7)\cup(10,11)\). Define \(\mu\colon\Z\to[0,1]\) by \(\mu(n)=1\) when \(n\equiv 0\pmod{8}\), by \(\mu(n)=1/2\) when \(n\equiv 7,9,11\pmod{12}\), and by \(\mu(n)=0\) otherwise, and put \(\nu=\sum_{n\neq 0}\mu(n)\,\delta_n\). Then \(\Omega\) is spectral, \(\nu\) is a weak tiling measure for \(\Omega\), and \(\nu\) is not a convex combination of proper tilings of the complement.
\end{theorem}

\begin{proof}
Write \(A=\{0,4,6,10\}\) and \(S=\{0,1/8,2/8,3/8\}\). For each \(r\in\{1,2,3\}\),
\[
\sum_{a\in A}\ex{ra/8}=0.
\]
The three sums are \(1-1-i+i\), \(1+1-1-1\), and \(1-1+i-i\). Thus \(A\) is a spectrum of \(S\). For \(h=1\) and for almost every \(x\in[0,1)\), the fiber of \(\Omega\) is exactly \(A\). The converse half of Theorem~\ref{thm:fiber} says that \(S+\Z\) is a spectrum of \(\Omega\).

The two rules defining \(\mu\) are disjoint, because every integer congruent to \(7\), \(9\), or \(11\) modulo \(12\) is odd, while every multiple of \(8\) is even. Let \(y\in\Z\). The residues of \(A\) modulo \(8\) are \(\{0,2,4,6\}\), each once. If \(y\) is even, the four integers \(y-a\) are even, the atoms of mass \(1/2\) do not contribute, and exactly one of the four is divisible by \(8\). The sum is \(1\). If \(y\) is odd, only the atoms of mass \(1/2\) contribute. Modulo \(12\), the set \(A\) is the union of the pairs \(\{0,6\}\) and \(\{4,10\}\). Translation by \(6\) interchanges \(\{7,9,11\}\) with its complement \(\{1,3,5\}\) among the odd residues, so exactly one value of \(y-a\) from each pair belongs to \(\{7,9,11\}\) modulo \(12\). The sum is again \(1\). Thus \(1_A*\mu=1\) on \(\Z\).

Identify \(\mu\) with the pure point measure \(\sum_{n\in\Z}\mu(n)\delta_n\). For non-integral \(x\), with \(m=\lfloor x\rfloor\), one has \(x-n\in\Omega\) iff \(m-n\in A\). Indeed each component of \(\Omega\) is a unit interval with integer endpoints, and the fractional part of \(x\) selects exactly the integer \(n=m-a\). Therefore
\[
(1_\Omega*\mu)(x)=\sum_{a\in A}\mu(m-a)=1.
\]
The equality \(1_\Omega*\mu=1\) holds almost everywhere. Since \(\mu(0)=1\) and \(\nu=\mu-\delta_0\) is positive and locally finite,
\[
1_\Omega*\nu=1_{\Omega^c}
\]
almost everywhere.

It remains to classify the proper tilings and compare them with \(\nu\). A proper tiling of the complement is \(\sum_{t\in T'}\delta_t\) for a set \(T'\) such that the translates \(\Omega+t\), \(t\in T'\), partition \(\Omega^c\) up to null sets. Then \(0\notin T'\), and adjoining the original copy of \(\Omega\) gives a tiling of \(\R\) by \(\Omega\) with translation set \(T=T'\cup\{0\}\), which contains \(0\). The gaps of \(\Omega\) are \((1,4)\), \((5,6)\), and \((7,10)\), of integer lengths \(3\), \(1\), and \(3\), and the components of \(\Omega\) occupy \((0,1)\), \((4,5)\), \((6,7)\), and \((10,11)\). In a tiling of \(\R\) by translates of \(\Omega\) whose translation set contains \(0\), every tile component is an open interval of length \(1\), disjoint from \(\Omega\) or equal to one of its components, and no such interval can cross an endpoint of \(\Omega\). The gap \((1,4)\) is therefore partitioned into intervals of length \(1\). The only such partition is \((1,2)\cup(2,3)\cup(3,4)\). The same matching of endpoints propagates through every later gap and through the ray \((-\infty,0)\). The tiling partition of \(\R\) is exactly \(\{(n,n+1):n\in\Z\}\). Every translation in the tiling set is an integer, because it carries some component of \(\Omega\) onto some \((n,n+1)\). Since \(\Omega=A+(0,1)\), a set \(T\subset\Z\) is the translation set of a tiling of \(\R\) by \(\Omega\) exactly when \(A\oplus T=\Z\).

Let \(T\subset\Z\) contain \(0\) and satisfy \(A\oplus T=\Z\). Then \(T\) contains no nonzero element of \(A-A=\{0,\pm 2,\pm 4,\pm 6,\pm 10\}\). Starting from \(8k\in T\), the integer \(8k+2\) has only one admissible partner in \(T\): the four candidates are \(8k+2\), \(8k-2\), \(8k-4\), and \(8k-8\), and the first three differ from \(8k\) by an element of \(A-A\). Thus \(8k-8\in T\). The integer \(8k+8\) likewise has only the admissible partner \(8k+8\), the other three candidates \(8k+4\), \(8k+2\), and \(8k-2\) being excluded by the same differences. By induction \(8\Z\subset T\). An even integer outside \(8\Z\) has residue \(2\), \(4\), or \(6\) modulo \(8\) and therefore differs from an element of \(8\Z\) by an element of \(A-A\). Hence \(T\cap 2\Z=8\Z\).

If \(b\in T\) is odd, the same four-candidate count at \(b+2\) forces \(b-8\in T\), and the count at \(b+8\) forces \(b+8\in T\). Thus \(T\) is a union of cosets of \(8\Z\). The set \(A\) meets every even residue class modulo \(8\) once, so \(A\oplus(c+8\Z)\) is the set of all odd integers, each represented once, for every odd \(c\). Two distinct odd cosets would represent each odd integer twice. Therefore \(T=8\Z\cup(c+8\Z)\) for a unique \(c\in\{1,3,5,7\}\). Each of these four sets does tile \(\Z\) with \(A\). A proper tiling measure of the complement is \(\nu_c=\sum_{t\in T\setminus\{0\}}\delta_t\) for one of these \(T\). Any convex combination of the \(\nu_c\), and any integral of them against a probability measure on \(\{1,3,5,7\}\), assigns one and the same mass to every nonzero point of a fixed coset \(c+8\Z\). This is incompatible with \(\nu\): one has \(\mu(7)=1/2\) and \(\mu(15)=0\), while \(15\equiv 7\pmod{8}\). There are only four proper tilings, so their convex combinations form a compact set of measures, which is closed in the vague topology (against compactly supported continuous functions). The measure \(\nu\) is therefore not a vague limit of convex combinations of proper tilings either.
\end{proof}

\section{Remarks}
\label{sec:remarks}

\begin{remark}
\label{rem:nintervals}
For a spectral union of any finite number \(n\) of intervals, the boundary unitary, the secular determinant and the uniqueness principle are available, and Section~\ref{sec:boundary} is written for every \(n\). The cofactor identity of Theorem~\ref{thm:cofactor} has the same proof for every \(n\). A monomial boundary unitary forces the set to be a fundamental domain of the lattice generated by its measure, and equal component lengths produce an \(n\)-fold lattice multitile by the lattice \((L/n)\Z\), where \(L=|\Omega|\) \cite[Section~10]{shi2026}. Nothing in those conclusions uses a restriction on \(n\).

The same is true of the fiber reformulation, Theorem~\ref{thm:fiber}, which holds for every bounded measurable set. The passage from a lattice multitile to a translational tiling when the spectrum has \(k\le 5\) cosets is the result of Dutkay and Haussermann \cite{dutkay-haussermann2015}, formulated here as Theorem~\ref{thm:fiberfew}. For a general spectrum with six or more cosets that passage remains open. The present classification does not establish the bound \(k\le n\) for arbitrary \(n\), with \(k\) taken in the minimal-period representation. In particular the method does not by itself give the analog of Theorem~\ref{thm:main} for \(n > 4\).
\end{remark}

\appendix
\section{Proof of the rationality theorem}
\label{sec:rational}

This appendix proves Theorem~\ref{thm:rational}. The argument is that of Fu and Song \cite[Sections~4 and~5]{fu-song2026}. Their modulus rigidity theorem \cite[Theorem~1.1]{fu-song2026} is Lemma~\ref{lem:rigidity} below, which is obtained here as the case \(N=d\) of Lemma~\ref{lem:equalmoduli}, and their positive definite kernel \cite[Lemma~4.1]{fu-song2026} is Lemma~\ref{lem:kernel}.

A finite spectral pair with integer frequencies gives a complex Hadamard matrix of generalized Vandermonde type, with entries \(z_s^a\) for \(z_s=e^{2\pi i s}\) and \(a\in A\). After translating \(S\) so that \(0\in S\), one node is \(1\). The proof below shows that the nodes are algebraic integers; they are not yet known to be roots of unity. A field automorphism preserves the inverse-orthogonality relations but may change the moduli of the nodes. Lemma~\ref{lem:rigidity} shows that all conjugate nodes have the same modulus, hence modulus \(1\), allowing Kronecker's theorem to be applied.

Throughout this appendix, \(d\) denotes the common cardinality. The proof of Theorem~\ref{thm:rational} combines the modulus rigidity of Lemma~\ref{lem:rigidity} with a Galois argument.

\begin{lemma}
\label{lem:columns}
Let \(w_1,\ldots,w_d\) be nonzero complex numbers and let \(A\subset\Z\) have \(d\) elements. Suppose
\[
\sum_{a\in A}\Bigl(\frac{w_j}{w_i}\Bigr)^a=d\,\delta_{ij}
\qquad(1\le i,j\le d).
\]
Then
\[
\sum_{j=1}^d w_j^{a-b}=d\,\delta_{ab}
\qquad(a,b\in A).
\]
\end{lemma}

\begin{proof}
Index rows by \(A\) and columns by \(\{1,\ldots,d\}\), and let \(M\) be the matrix with entries \(M_{aj}=w_j^a\). Let \(M^\vee\) be the matrix with entries \(M^\vee_{aj}=w_j^{-a}\). The \((i,j)\) entry of \((M^\vee)^\top M\) is \(\sum_{a\in A}w_i^{-a}w_j^a\), so the hypothesis reads \((M^\vee)^\top M=dI\). Since the matrices are square, the reverse product also satisfies \(M(M^\vee)^\top=dI\). The \((a,b)\) entry of the last identity is \(\sum_j w_j^{a-b}\).
\end{proof}

\begin{lemma}
\label{lem:rigidity}
Let \(w_1,\ldots,w_d\) be nonzero complex numbers and let \(A\subset\Z\) have \(d\) elements. Suppose \(\sum_{a\in A}(w_j/w_i)^a=d\,\delta_{ij}\) \(\forall 1\le i,j\le d\). Then \(|w_1|=\cdots=|w_d|\).
\end{lemma}

\begin{proof}
Lemma~\ref{lem:columns} gives \(\sum_{j=1}^d w_j^{a-b}=0\) \(\forall a\neq b\) in \(A\), that is, the power sums of \(w_1,\ldots,w_d\) vanish on \((A-A)\setminus\{0\}\). Lemma~\ref{lem:equalmoduli}, applied with \(A'=A\) and \(N=d\), gives the conclusion, since \(d<2d\).
\end{proof}

\begin{proof}[Proof of Theorem~\ref{thm:rational}]
The case \(d=1\) is the assertion \(\{0\}\subset\Q/\Z\) after a translation of \(S\). Assume \(d\ge 2\), and fix \(s_0\in S\). The translate \(\Gamma=S-s_0\) satisfies \(0\in\Gamma\), has \(d\) elements, and \(\widehat A\) vanishes on \(\Gamma-\Gamma\) except at the origin. It is enough to prove \(\Gamma\subset\Q/\Z\), since then \(S-S=\Gamma-\Gamma\subset\Q/\Z\). The set \(S\) itself need not meet \(\Q/\Z\).

Put \(c=\min A\) and replace \(A\) by \(A-c\). Then \(\min A=0\), and \(\widehat A\) changes by the unimodular factor \(e^{-2\pi i c\,\cdot}\), so the zero set is unchanged. The mask polynomial \(\mask{A}=\sum_{a\in A}X^a\in\Z[X]\) is monic of degree \(\max A\ge 1\) and has constant term \(1\). Write \(z_j=e^{2\pi i\gamma_j}\) for the elements \(\gamma_j\) of \(\Gamma\), with \(z_0=1\). The spectral hypothesis, including the diagonal, is
\[
\sum_{a\in A}\Bigl(\frac{z_j}{z_i}\Bigr)^a=d\,\delta_{ij}.
\]
The case \(i=0\) says that \(\mask{A}\) vanishes at \(z_j\) \(\forall j\neq 0\). The same polynomial takes the value \(d\) at \(z_0=1\), so \(1\) is not a root. Each \(z_j\) with \(j\neq 0\) is an algebraic integer, and so is \(z_0=1\).

Let \(K\) be the splitting field over \(\Q\) of the product of the minimal polynomials of \(z_0,\ldots,z_{d-1}\). Then \(K/\Q\) is finite and normal, and \(K\) contains every \(z_j\) and every algebraic conjugate of every \(z_j\). If \(\beta\) is a conjugate of \(z_j\), the assignment \(f(z_j)\mapsto f(\beta)\) is a field embedding of \(\Q(z_j)\) into \(K\). It extends to an embedding of \(K\) into an algebraic closure of \(\Q\). The image of \(K\) lies in \(K\), because \(K\) is normal, and an injective endomorphism of a finite-dimensional vector space is an automorphism. Thus \(\beta=\sigma(z_j)\) for some \(\sigma\in\operatorname{Gal}(K/\Q)\).

Fix \(\sigma\in\operatorname{Gal}(K/\Q)\) and set \(w_j=\sigma(z_j)\). Each \(w_j\) is nonzero, and \(w_0=\sigma(1)=1\). The displayed relations are polynomial identities with integer coefficients in the nodes and their inverses, so they pass to the conjugates \(w_j\). Lemma~\ref{lem:rigidity} gives \(|w_0|=\cdots=|w_{d-1}|\). Since \(w_0=1\), every \(|\sigma(z_j)|\) equals \(1\). Every algebraic conjugate of every \(z_j\) therefore has modulus \(1\).

Kronecker's theorem~\cite{kronecker1857} says that an algebraic integer whose every algebraic conjugate has modulus \(1\) is a root of unity. Hence \(z_j^{N_j}=1\) for some positive integer \(N_j\), so \(\gamma_j\in\Q/\Z\). If \(N\) is a common multiple of the \(N_j\), then \(\Gamma\subset N^{-1}\Z/\Z\).
\end{proof}

\section*{Acknowledgments}
I thank everyone who contributed useful discussions or took the time to read the paper and suggest edits, including Zi-jian Song, whose encouragement in our brief interactions and substantial contributions to the field I greatly appreciate.

\section*{Funding}
No funding was received for this research.
\clearpage
\bibliographystyle{amsplainfull}
\makeatletter\def\@biblabel#1{\@defaultbiblabelstyle{#1}}\def\bibsetup{}\makeatother
\bibliography{refs}
\end{document}